\documentclass[12pt]{article}

\usepackage{rotating}
\usepackage{graphicx}													
\usepackage{amsmath}
\usepackage{amsthm}						
\usepackage{amsfonts}								
\usepackage{amssymb}
\usepackage{tabu}
\usepackage{floatpag}

\newcommand\al{\alpha}
\newcommand\bt{\beta}
\newcommand\gm{\gamma}
\newcommand\dl{\delta}
\newcommand\lm{\lambda}
\newcommand\sg{\sigma}
\newcommand\Gm{\Gamma}
\newcommand\Dl{\Delta}
\newcommand\Sg{\Sigma}

\newcommand\F{\mathbb{F}}

\newcommand\cF{{\cal F}}

\newcommand\cJ{{\cal J}}
\newcommand\cM{{\cal M}}
\newcommand\cP{{\cal P}}

\newcommand\ad{\mathrm{ad}}

\newcommand\Aut{\mathrm{Aut}\,}

\newcommand\la{\langle}
\newcommand\ra{\rangle}
\newcommand\dla{\la\!\la}
\newcommand\dra{\ra\!\ra}
\newcommand\supp{\mathrm{supp}}

\newcommand{\ha}{\frac{1}{2}}
\newcommand{\qu}{\frac{1}{4}}
\newcommand{\thi}{\frac{1}{32}}

\newtheorem{lemma}{Lemma}[section]

\newtheorem{theorem}[lemma]{Theorem}
\newtheorem{definition}[lemma]{Definition}
\newtheorem{corollary}[lemma]{Corollary}
\newtheorem{proposition}[lemma]{Proposition}
\newtheorem{question}[lemma]{Question}

\renewcommand{\arraystretch}{1.5}

\title{Double axes and subalgebras of Monster type in Matsuo algebras}
\author{Alexey Galt, Vijay Joshi, Andrey Mamontov,\\
Sergey Shpectorov, Alexey Staroletov}

\begin{document}
\maketitle
\newcommand{\Addresses}{{
		\bigskip\noindent
		\footnotesize
		Alexey~Galt, \textsc{Sobolev Institute of Mathematics, Novosibirsk, Russia;}\\\nopagebreak
		\textsc{Novosibirsk State University, Novosibirsk, Russia;}\\\nopagebreak
		\textit{E-mail address: } \texttt{galt84@gmail.com}
		
		\medskip\noindent
		Vijay~Joshi, \textsc{School of Mathematics, University of Birmingham,
			Watson Building, Edgbaston, Birmingham B15 2TT, United Kingdom;}
		\\\nopagebreak
		\textit{E-mail address: } \texttt{vijay6d@gmail.com}
		
		\medskip\noindent
		Andrey~Mamontov, \textsc{Sobolev Institute of Mathematics, Novosibirsk, Russia;}\\\nopagebreak
		\textsc{Novosibirsk State University, Novosibirsk, Russia;}\\\nopagebreak
		\textit{E-mail address: } \texttt{andreysmamontov@gmail.com}
		
		\medskip\noindent
		Sergey~Shpectorov, \textsc{School of Mathematics, University of Birmingham,
		Watson Building, Edgbaston, Birmingham B15 2TT, United Kingdom;}
		\\\nopagebreak
		\textit{E-mail address: } \texttt{s.shpectorov@bham.ac.uk}

		\medskip\noindent
		Alexey~Staroletov, \textsc{Sobolev Institute of Mathematics, Novosibirsk, Russia;}\\\nopagebreak
		\textsc{Novosibirsk State University, Novosibirsk, Russia;}\\\nopagebreak
		\textit{E-mail address: } \texttt{staroletov@math.nsc.ru}
		
		\medskip
}}

\begin{abstract}
Axial algebras are a class of commutative non-associative algebras 
generated by idempotents, called axes, with adjoint action semi-simple 
and satisfying a prescribed fusion law. Axial algebras were introduced 
by Hall, Rehren and Shpectorov \cite{hrs,hrs1} as a broad generalization 
of Majorana algebras of Ivanov, whose axioms were derived from the 
properties of the Griess algebra for the Monster sporadic simple group. 
The class of axial algebras of Monster type includes Majorana algebras for 
the Monster and many other sporadic simple groups, Jordan algebras for 
classical and some exceptional simple groups, and Matsuo algebras 
corresponding to $3$-transposition groups. Thus, axial algebras of Monster 
type unify several strands in the theory of finite simple groups.

It is shown here that double axes, i.e., sums of two orthogonal axes in a 
Matsuo algebra, satisfy the fusion law of Monster type $(2\eta,\eta)$. 
Primitive subalgebras generated by two single or double axes are 
completely classified and $3$-generated primitive subalgebras are 
classified in one of the three cases. These classifications further lead 
to the general flip construction outputting a rich variety of axial 
algebras of Monster type. An application of the flip construction to the 
case of Matsuo algebras related to the symmetric groups results in three 
new explicit infinite series of such algebras.     
\end{abstract}

\section{Introduction}
\label{Introduction}

Axial algebras are commutative non-associative algebras generated by a set
of special idempotents, called axes, satisfying a prescribed fusion law. 
For example, the $196,884$-dimensional real unital Griess algebra $V$ is 
generated by a set of idempotents called $2A$-axes that satisfy the fusion 
law $\cM(\qu,\thi)$ shown in Table \ref{Monster}. 
\begin{table}[h]
\begin{center}
\begin{tabular}{|c||c|c|c|c|}
\hline
$\ast$&$1$&$0$&$\qu$&$\thi$\\
\hline\hline
$1$&$1$&&$\qu$&$\thi$\\
\hline
$0$&&$0$&$\qu$&$\thi$\\
\hline
$\qu$&$\qu$&$\qu$&$1,0$&$\thi$\\
\hline
$\thi$&$\thi$&$\thi$&$\thi$&$1,0,\qu$\\
\hline
\end{tabular}
\end{center}
\caption{Fusion law $\cM(\qu,\thi)$}\label{Monster}
\end{table}
The meaning of this is that the adjoint action $\ad_a:v\mapsto av$ of any 
$2A$-axis $a\in V$ is semi-simple (i.e., diagonalizable), with all 
eigenvalues in the set $\{1,0,\qu,\thi\}$ and the entries of the table 
prescribe how eigenvectors of $\ad_a$ multiply in $V$. For example, the 
entry $\thi\ast 0=\{\thi\}$ means that if $u$ is a $\thi$-eigenvector and $v$ is a 
$0$-eigenvector then $uv$ is a $\thi$-eigenvector; the entry 
$\qu\ast\qu=\{1,0\}$ means that the product of a $\qu$-eigenvector and a second 
$\qu$-eigenvector is a sum of a $1$-eigenvector and a $0$-eigenvector, and so 
on. 

Historically, the first class of algebras whose axioms included a fusion law 
were the Majorana algebras of Ivanov \cite{i}. These are real axial algebras 
satisfying the fusion law $\cM(\qu,\thi)$ plus certain additional axioms. The 
Griess algebra $V$ is the principal example of a Majorana algebra, and in fact, 
the axioms of Majorana algebras were derived from the properties of the Griess 
algebra exploited in the celebrated theorem of Sakuma \cite{s}. Note that 
Sakuma's theorem was proven in the context of vertex operator algebras (VOAs), 
where the Griess algebra arises as the weight-two component of the Moonshine 
VOA $V^\natural$, whose automorphism group is the sporadic simple group $M$, 
the Monster. 

There are further classes of algebras that exhibit axial behavior. Idempotents 
in associative algebras satisfy the simplest fusion law involving only 
$\{1,0\}$. More interestingly, in Jordan algebras, every idempotent leads to 
the so-called Peirce decomposition of the algebra, and this property can be 
restated in axial terms as the fusion law $\cJ(\ha)$ in Table \ref{Jordan}.    
\begin{table}[h]
\begin{center}
\begin{tabular}{|c||c|c|c|}
\hline
$\ast$&$1$&$0$&$\ha$\\
\hline\hline
$1$&$1$&&$\ha$\\
\hline
$0$&&$0$&$\ha$\\
\hline
$\ha$&$\ha$&$\ha$&$1,0$\\
\hline
\end{tabular}
\end{center}
\caption{Fusion law $\cJ(\ha)$}\label{Jordan}
\end{table}
Motivated by these and other examples, Hall, Rehren, and Shpectorov
\cite{hrs,hrs1} introduced axial algebras as a broad generalization of Majorana 
and Jordan algebras, allowing an arbitrary field, an arbitrary fusion law, 
and getting rid of all additional ``non-axial'' axioms. Focussing on a 
particular fusion law and/or adding additional axioms defines a subclass of 
axial algebras. In this paper, we work with two subclasses: algebras of Jordan 
type $\eta$ satisfying the fusion law $\cJ(\eta)$ as in Table \ref{Jordan}, but 
with an arbitrary $\eta\neq 1,0$ in place of $\ha$, and algebras of Monster 
type $(\al,\bt)$ for the fusion law $\cM(\al,\bt)$ as in Table \ref{Monster}, 
but with $\qu$ and $\thi$ substituted with arbitrary $\al$ and $\bt$ distinct 
from $1$, $0$ and from each other. The only other axiom that we require is 
primitivity, which refers to primitivity of each generating axis $a$ of the 
algebra $A$, understood as the condition that the $1$-eigenspace $A_1(a)$ of 
$\ad_a$ be $1$-dimensional spanned by $a$. (Note that since $a$ is an 
idempotent, it is a $1$-eigenvector for $\ad_a$.)

We note that the class of algebras of Jordan type is a subclass of algebras of 
Monster type, corresponding to the situation where the $\al$-eigenspace is 
trivial and $\bt=\eta$. 

We have already mentioned that the automorphism group of the Griess algebra 
is the Monster $M$. Jordan algebras also admit large automorphism groups, among 
which we find all classical groups and some exceptional groups of Lie type.
The apparent relation between axial algebras of Monster type and (simple) 
groups is for a good reason: the fusion laws $\cJ(\eta)$ and $\cM(\al,\bt)$ 
are $C_2$-graded and this results in the axial algebra $A$ admitting, for each 
axis $a$, an automorphism $\tau_a$ of order two, called the Miyamoto involution 
of $a$. The Miyamoto involutions generate the Miyamoto group of $A$, which is a 
subgroup of $\Aut(A)$. Naturally, the Miyamoto group of the Griess algebra $V$ 
is the Monster group and, similarly, for the Jordan algebras, the Miyamoto 
groups are the corresponding groups of Lie type. We stress that it is not just 
the known examples that are highly symmetric: all algebras of Monster type, 
whether known or not, admit significant automorphism groups. Hence the study 
of algebras of Monster type is highly interesting and relevant not only for 
specialists in non-associative algebras, but also for group theorists, because 
these algebras provide a platform to study sporadic simple groups alongside the 
groups of Lie type, giving the hope of a possible unified theory of (finite) 
simple groups.

Ideally, the aim of the theory is a complete classification of algebras of 
Monster type. For the smaller class of algebras of Jordan type, it is shown in 
\cite{hrs1} (with a correction in \cite{hss}) that, for $\eta\neq\ha$, algebras 
of Jordan type are exactly the Matsuo algebras, corresponding to 
$3$-transposition groups, or their factor algebras. The class of Matsuo 
algebras, playing an important r\^ole in this paper, was first introduced by 
Matsuo \cite{m} and then reintroduced and generalized by Hall, Rehren, and 
Shpectorov in \cite{hrs1}. The case $\eta=\ha$ is still open and it is 
significantly more complex because Jordan algebras appear in this case. 
However, the axial-style characterization of Jordan algebras by Segev \cite{se} 
and the recent result of Gorshkov and Staroletov \cite{gs}, bounding 
the dimension of a $3$-generated algebra of Jordan type, give hope that the 
classification of algebras of Jordan type can soon be completed. 

For the algebras of Monster type, classification results are currently few. 
Sakuma's Theorem was reproved by Ivanov, Pasechnik, Seress, and Shpectorov in 
\cite{ipss} in the context of Majorana algebras. Next, Hall, Rehren, and 
Shpectorov \cite{hrs} generalized the theorem by removing several assumptions 
specific to Majorana algebras leaving only the fusion law $\cM(\qu,\thi)$, 
primitivity, and the existence of a Frobenius form, a non-zero bilinear form 
associating with the algebra product. Later Rehren \cite{r,r1} attempted to 
remove the Frobenius form condition and also generalize to the case of arbitrary 
$2$-generated primitive algebras of Monster type $(\al,\bt)$. He managed to 
obtain an upper bound of $8$ for the dimension of such an algebra, modulo 
the extra assumptions that $\al\neq 2\bt,4\bt$. In a recent preprint 
\cite{fms}, Franchi, Mainardis, and Shpectorov make the next step, extending 
this bound, $8$, to one of the exceptional cases, $\al=2\bt$, and finding the 
real culprit, an infinite-dimensional $2$-generated algebras 
$HW$\footnote{$HW$ stands for ``highwater algebra'' because the algebra was 
discovered in Venice during the disastrous floods in November 2019.}, in the 
second exceptional case where $\al=4\bt$. They also remove the Frobenius 
form condition in the case of the law $\cM(\qu,\thi)$, thus achieving the 
ultimate classification in this key case. The general classification of 
$2$-generated algebras of Monster type remains an open and highly attractive 
problem.

Beyond the $2$-generated case, the results have been limited so far to finding 
subalgebras of the Griess algebra and computing the algebras of type $(\qu,\thi)$ 
for concrete Miyamoto groups by hand or using available GAP and MAGMA programs 
\cite{ipss,ish,is,i1,s,d,i2,i3,fim,w,pw,ms}. (See also the interesting survey 
\cite{fim1}.) While these computations lead to many new interesting examples of 
algebras, the current situation is somewhat paradoxical. Intuitively, it is 
clear that the class of algebras of Monster type should be significantly richer 
than its subclass of algebras of Jordan type. However, all (finitely many!) 
known examples outside of this subclass have been of type $(\qu,\thi)$. Clearly, 
the class of algebras of Monster type $(\al,\bt)$ should be suitably populated 
with examples before one can expect more general classification results.

In this paper we make a step towards remedying this situation. Namely, 
we introduce the flip construction that results in a rich variety of primitive 
algebras of Monster type $(2\eta,\eta)$, where $\eta$ is any scalar outside of 
$\{0,\ha,1\}$. It is interesting that $\al=2\eta$ and $\bt=\eta$ is one of the 
two exceptional situations ($\al=2\bt$) arising in the work of Rehren \cite{r}.
The flip construction starts with a Matsuo algebra $M=M_\eta(\Gm)$, where $\Gm$ 
is the Fischer space of an arbitrary $3$-transposition group. Recall that a 
$3$-transposition group is a pair $(G,D)$ where $G$ is a group generated by a 
normal subset $D$ of involutions satisfying the condition that $|cd|\leq 3$ for 
all $c,d\in D$. The Fischer space $\Gm$ of $(G,D)$ is the point-line geometry 
that has the involutions from $D$ as its points and, as lines, all triples 
$\{c,d,e\}\subseteq D$ satisfying $c^d=e=d^c$. Note that points $c$ and $d$ 
are collinear in $\Gm$ if and only if they do not commute as elements of $G$.

The Matsuo algebra $M_\eta(\Gm)$, defined over a field $\F$ of characteristic 
not equal to $2$, has $D$ as its basis, with multiplication defined on this 
basis as follows:
$$c\cdot d=\left\{
\begin{array}{rl}
c,&\mbox{if }c=d;\\
0,&\mbox{if $c$ and $d$ are non-collinear};\\
\frac{\eta}{2}(c+d-e),&\mbox{if $c$ and $d$ lie on a line, $\{c,d,e\}$}.
\end{array}
\right.$$
Here the parameter $\eta$ is an arbitrary element of $\F$, not equal to $1$ or 
$0$. (For an extended discussion of these concepts, see Section \ref{Matsuo}.) 

The elements of the basis $D$ of $M=M_\eta(\Gm)$ are the primitive axes of 
this algebra, satisfying the fusion law $\cJ(\eta)$ in Table \ref{Jordan 
type}. Note that non-collinear axes $a$ and $b$ are orthogonal in $M$, i.e., 
$a\cdot b=0$. In this case, $x=a+b$ is again an idempotent in $M$. We call 
the elements $x$ \emph{double axes}, as opposed to single axes, which are 
the elements of $D$.

Our first result concerns the fusion law that a double axis satisfies.

\begin{theorem}
If $\eta\neq\ha$ then the double axis $x=a+b$ satisfies the fusion law 
$\cM(2\eta,\eta)$:
\begin{table}[h]
\begin{center}
\begin{tabular}{|c||c|c|c|c|}
\hline
$\diamond$&$1$&$0$&$2\eta$&$\eta$\\
\hline\hline
$1$&$1$&&$2\eta$&$\eta$\\
\hline
$0$&&$0$&$2\eta$&$\eta$\\
\hline
$2\eta$&$2\eta$&$2\eta$&$1+0$&$\eta$\\
\hline
$\eta$&$\eta$&$\eta$&$\eta$&$1+0+2\eta$\\
\hline
\end{tabular}
\end{center}
\end{table} 
\end{theorem}

(This is proved in Section \ref{double} as Proposition \ref{double fusion 
rules}.)

That is, double axes are of Monster type $(2\eta,\eta)$. Note that single axes 
also satisfy the above law, with the $2\eta$-eigenspace trivial. Thus, any 
subalgebra in a Matsuo algebra generated by a collection of single and double 
axes is an algebra of Monster type $(2\eta,\eta)$. The pinch of salt in this 
statement is the absence of the word ``primitive''. Indeed, the double axis 
$x=a+b$ is not primitive in $M$, because its $1$-eigenspace is $2$-dimensional, 
equal to $\la a,b\ra$. However, $x$ can be primitive in a proper subalgebra of 
$M$, as long as this subalgebra does not contain $a$ and $b$. The question is 
therefore: can we find a large family of subalgebras, where the generating 
double axes are primitive?

We first look at the $2$-generated case, that is, at subalgebras generated by 
two axes. If both axes are single, say, $a,b\in D$, then $\dla 
a,b\dra$\footnote{We use double brackets for algebra generation, to distinguish 
it from the linear span of the vectors.} is well-known. If $a$ and $b$ are 
non-collinear then $\dla a,b\dra\cong\F^2$, the algebra known in the axial lingo 
as $2B$. If $a$ and $b$ are collinear then $\dla a,b\dra$ is a $3$-dimensional 
Matsuo algebra known as $3C(\eta)$.

So we assume that at least one of the two generating axes is a double axis. 
Accordingly, we consider two cases: either (A) the algebra is generated by a 
single axis $a$ and a double axis $b+c$; or (B) it is generated by two double 
axes, $a+b$ and $c+d$. The exact configuration is described by the diagram on 
the support set $\{a,b,c\}$ or $\{a,b,c,d\}$, the diagram being the graph in 
which two points from the support set are connected by an edge when they are 
collinear in the Fischer space $\Gm$. The possible diagrams for cases (A) and 
(B) are given in Figures \ref{Type A} and \ref{Type B}. Each diagram defines a 
Coxeter group, of which the $3$-transposition group $G$ (the minimal $G$, 
generated just by the support involutions) is a factor group. All these Coxeter 
groups are well-known, and so it is not difficult to list all possible $G$ and 
then build the corresponding $2$-generated subalgebras. The number of subcases 
is eliminated by Theorem \ref{symmetry} stating that the subalgebra cannot be 
primitive unless the diagram admits the flip symmetry. By the \emph{flip} we 
mean the permutation of the support set switching two points within each double 
axis. This eliminates diagrams $A_2$, $B_2$, $B_3$, and $B_5$. The remaining 
symmetric diagrams all lead to primitive $2$-generated subalgebras.

We can summarize this discussion as the following result.

\begin{theorem} \label{main 2-generated}
Every $2$-generated primitive subalgebra is one of the following: $2B\cong\F^2$, 
$3C(\eta)$, $3C(2\eta)$, the new $4$-dimensional algebra $Q_2(\eta)$, the new 
$5$-dimensional algebra shown in Table $\ref{t:B6(2)}$, or the new 
$8$-dimensional algebra in Table $\ref{t:B6(3)}$.
\end{theorem}

It is interesting to compare this with the list of Sakuma algebras arising 
for the Monster fusion law. Similarities between the two lists are quite obvious 
and the main difference is the absence here of an analogue of the Sakuma algebra 
$5A$ with the Miyamoto group $D_{10}$. Of course, the algebras above are the 
ones that live inside Matsuo algebras. Hypothetically, other $2$-generated 
primitive algebras of Monster type $(2\eta,\eta)$ might exist. Hence we are 
asking the following: Is the list in Theorem \ref{main 2-generated} the complete 
list of $2$-generated primitive algebras of Monster type $(2\eta,\eta)$?

We also attempt to determine $3$-generated primitive subalgebras of Matsuo 
algebras. Here again the case of three single axes is well understood. Hence one 
needs to deal only with the following cases: (C) two single axes and one double 
axis; (D) one single axis and two double axes; and (E) three double axes. The 
flip-symmetric diagrams in each of these cases are shown in Figures \ref{Type 
C}, \ref{Type D}, and \ref{Type E}. In this paper we complete the case (C). 
Some of the algebras arising here decompose as direct sums of smaller algebras. 
In the following theorem we only list the indecomposable ones.

\begin{theorem} \label{main 3-generated}
Each $3$-generated primitive indecomposable subalgebra in case (C) is one of the 
following: the new $9$-dimensional algebra $Q^3(\eta)$, the new $8$-dimensional 
algebra $2Q_2(\eta)$, the new $12$-dimensional algebra $3Q_2(\eta)$, or the new 
$24$-dimensional algebra shown in Table \ref{t:algebra C6-24}.
\end{theorem}

The classification of cases (D) and (E) is an ongoing project. For (D), the 
support set has cardinality $5$. So here one can use the classification of 
$5$-generated $3$-transposition groups by Hall and Soicher \cite{hs}. 
Unfortunately, there is no classification of $6$-generated $3$-transposition 
groups, which is what one would need for (E).

It is remarkable that all algebras found so far for flip-symmetric diagrams turn 
out to be primitive. This suggests that perhaps the flip symmetry condition is 
not only necessary for primitivity, but also sufficient. Developing this idea 
further, we propose in this paper the general flip construction. Consider an 
arbitrary $3$-transposition group $(G,D)$ and an involutive automorphism $\tau$ 
of $G$ such that $D^\tau=D$. Then $\tau$ acts on the Fischer space $\Gm$ and, 
consequently, also on the Matsuo algebra $M=M_\eta(\Gm)$. Let 
$H=\la\tau\ra\leq\Aut(M)$. If $\tau$ flips two orthogonal Matsuo axes $a$ and 
$b$ then the fixed subalgebra $M_H$ of $M$ contains the double axis $a+b$, but 
not $a$ and $b$. Hence $x=a+b$ is guaranteed to be primitive in $M_H$. (This is 
Proposition \ref{primitive}.) Hence we define the \emph{flip algebra} $A_\tau$ 
as the subalgebra of $M_H$ generated by all single and (flipped) double axes 
contained in $M_H$. This is always a primitive axial algebra of Monster type 
$(2\eta,\eta)$. 

There exist great many $3$-transposition groups and each Fischer space admits 
many classes of involutive automorphisms. Thus, the flip construction does 
indeed lead to a rich variety of algebras of Monster type. The work of 
understanding these algebras is just starting. In this paper, we completely 
analyze the simplest case of the symmetric group $G=S_n$. The special case of 
the fixed-point-free involution $\tau\in S_n$ (when $n$ is even, $n=2k$) gives 
the first proper infinite series $Q_k(\eta)$ of algebras of Monster type 
$(2\eta,\eta)$.

\begin{theorem} \label{main Qk}
Suppose that $n=2k$ and $\tau=(1,2)(3,4)\cdots(2k-1,2k)$. Then the flip algebra 
$Q_k(\eta)=A_\tau$ is of dimension $k^2$, spanned by its $k$ single axes and 
$k^2-k$ double axes.
\end{theorem}

We call the algebra $Q_k(\eta)$ the \emph{$k^2$-algebra}.

In the general case, the flip algebra decomposes as the direct sum of two 
smaller algebras.

\begin{theorem} \label{main general flip}
If $n=2k+m$, with $k,m\neq 0$, and $\tau=(1,2)(3,4)\cdots(2k-1,2k)$ then 
$A_\tau=Q_k(\eta)\oplus M_\eta(\Gm')$, where the first summand arises on the 
subset $\{1,2,\ldots,2k\}$ and the second, Matsuo summand comes from the fixed 
subset $\{2k+1,\ldots,n\}$ of $\tau$.
\end{theorem}

Note that in Theorem \ref{main Qk} the flip algebra $Q_k(\eta)$ coincides with 
the full fixed subalgebra $M_H$, but this is not the case in Theorem \ref{main 
general flip}.

We also investigate the properties of the algebra $Q_k(\eta)$. In particular, 
we show the following using the structure theory developed by Khasraw, McInroy, 
and Shpectorov \cite{kms}.

\begin{theorem}
The algebra $Q_k(\eta)$ is simple unless $\eta=-\frac{1}{2(k-1)}$ or 
$-\frac{1}{k-2}$. For the two exceptional values of $\eta$, the algebra has a 
unique maximal nontrivial ideal (the radical), the factor over which is a simple 
algebra.
\end{theorem}

The dimension of the radical for the special values of $\eta$ is (almost always) 
$1$ and $k-1$, respectively. 

Towards the end of the paper, we also introduce two further infinite series of 
algebras, $2Q_k(\eta)$ of dimension $2k^2$ and $3Q_k(\eta)$ of dimension $3k^2$. 
These arise for $3$-transposition groups of the form $2^{n-1}:S_n$ and 
$3^{n-1}:S_n$ and special choices of $\tau$. For these groups, we do not analyze 
all possible flips $\tau$. We also do not determine the special values of 
$\eta$, for which $2Q_k(\eta)$ and $3Q_k(\eta)$ are not simple. All these 
questions are left for another paper. 

The idea of double axes came from GAP experiments performed in Novosibirsk. Some 
of the above theorems appeared in the PhD thesis of Joshi \cite{j} at the 
University of Birmingham. This includes the fusion law obeyed by double axes and 
the classification of the $2$-generated subalgebras. The infinite series 
$Q_k(\eta)$ was discovered simultaneously in Birmingham and Novosibirsk (likely 
within an hour of each other) as part of ongoing collaboration. During nearly 
two years of work on this paper, several mathematicians having access to our 
drafts have started looking at different cases of $3$-transposition groups. 
Joshi's PhD thesis includes, in addition to the symmetric case, also an almost 
complete treatment of flips for $G=Sp_{2n}(2)$. Alsaeedi \cite{al} constructed 
an infinite series $Q^k(\eta)$, dual to our $Q_k(\eta)$, as a subalgebra of the 
Matsuo algebra for $G={}^+\Omega_{k+1}^+(3)$. Shi \cite{shi} provided a treatment of 
$G=O_{2n}^{\pm}(2)$ similar to Joshi's treatment of $Sp_{2n}(2)$. The case of 
$G=SU_n(2)$ was looked at by Hoffman, Rodrigues, and Shpectorov \cite{hrs2}. 
Finally, all flip algebras for the sporadic Fischer groups and triality groups 
were computed by the 3A (the Novosibirsk part of the team) using GAP. All these 
results, with the exception of \cite{al}, are still awaiting publication. 

Of course, all this work only concerns the almost simple $3$-transposition 
groups, so it is just the tip of the iceberg. The variety of flip algebras is 
enormous. In a sense, the flip construction is similar to the construction of 
twisted groups of Lie type. There, as well, every automorphism defines a fixed 
subgroup, but some automorphisms are more interesting than others, leading to 
new simple groups. In the same way, for flip algebras, the question is which 
of them are the key algebras, the true building blocks, and which are the 
composites. 

Let us now discuss the contents of the present paper. In Section 2, we 
introduce the basics of axial algebras. In Section 3, we discuss the  
class of Matsuo algebras corresponding to 3-transposition groups, or in 
geometric terms, to Fischer spaces. In Section 4, we prove that double axes in 
Matsuo algebras satisfy the fusion law of Monster type $(2\eta,\eta)$. In 
Section 5, we classify all $2$-generated primitive subalgebras generated by 
single and double axes. In all but two cases we can find the suitable 
configuration inside a symmetric group, and so the calculations in these cases 
are particularly uncomplicated. In the remaining two cases, we operate, instead, 
in terms of the Fischer spaces for the groups $2^3:S_4$ and $3^3:S_4$. In 
Section 6, we outline the classification of $3$-generated primitive subalgebras 
and complete the case (C). In Section 7, we introduce the flip construction and 
show that it always leads to a primitive algebra. In Section 8, we analyze the 
case of the symmetric group group $S_n$ and develop the properties of the new 
series of algebras $Q_k(\eta)$. Finally, in Section 9, we construct two further 
infinite series, $2Q_k(\eta)$ and $3Q_k(\eta)$. 

\bigskip
{\bf Acknowledgement:} This work was partly supported by Mathematical Center 
in Akademgorodok under agreement No. 075-15-2019-1675 with the Ministry of Science 
and Higher Education of the Russian Federation.
\section{Axial algebras}

\label{Axial algebras}

In this section we provide the necessary background on axial algebras.
Note that algebras are not assumed to be associative. Hence an algebra is just
a vector space with a bilinear product.

\subsection{Axes and axial algebras}
\label{basics}

We start with basic definitions. 

\begin{definition}
A \emph{fusion law} $\cF$ over a field $\F$ consists of a (finite) set $X\subseteq\F$ 
and a symmetric (commutative) product
$$\ast:X\times X\to 2^X,$$
where $2^X$ denotes the set of all subsets of $X$.
\end{definition}

Let $A$ be a commutative algebra over $\F$. For $a\in A$, we denote by 
$\ad_a$ the adjoint map $A\to A$ defined by $u\mapsto au$. For $\lm\in\F$, 
let $A_\lm(a)$ denote the $\lm$-eigenspace of $\ad_a$; that is, 
$A_\lm(a)=\{u\in A\mid au=\lm u\}$. Note that this notation makes sense 
even if $\lm$ is not an eigenvalue of $\ad_a$; in this case, we simply have 
$A_\lm(a)=0$. For $L\subset\F$, set $A_L(a):=\oplus_{\lm\in L}A_\lm(a)$. 

Let $\cF=(X,\ast)$ be a fusion law over $\F$. 

\begin{definition}
A non-zero idempotent $a\in A$ is an \emph{($\cF$-)axis} if 
\begin{enumerate}
	\item[(a)] $A=A_X(a)$; and 
	\item[(b)] $A_\lm(a)A_\mu(a)\subseteq A_{\lm\ast\mu}(a)$ for all 
	           $\lm,\mu\in X$.
\end{enumerate}
\end{definition}

Condition (a) means that $\ad_a$ is semi-simple and all its eigenvalues 
are in $X$. Note that, since $a$ is an idempotent, $1$ is an eigenvalue 
of $\ad_a$ and $a\in A_1(a)$. Hence we always assume that $1\in X$. 

\begin{definition}
An axis $a$ is \emph{primitive} if $A_1(a)$ is $1$-dimensional; that is, 
$A_1(a)=\langle a\rangle$.
\end{definition}

We will mostly deal with primitive axes. In this case, $A_1(a)A_\lm(a)=
\langle a\rangle A_\lm(a)=A_\lm(a)$ if $\lm\neq 0$ and $A_1(a)A_0(a)=0$.
In view of this, we can assume that $1\ast 0=\emptyset$ (assuming that 
$0\in X$) and $1\ast\lm=\{\lm\}$ for all $0\neq\lm\in X$. 

\begin{definition}
An algebra $A$ over $\F$ is an \emph{($\cF$-)axial} algebra if it is 
generated as algebra by a set of $\cF$-axes. The algebra $A$ is a 
\emph{primitive} \emph{($\cF$-)axial} algebra if it is generated by a set 
of primitive $\cF$-axes.
\end{definition}

In principle, an axial algebra should be formally defined as the pair 
consisting of the algebra and the set of generating axes. However, in 
practice we will just talk about the algebra and the set of generating 
axes will be assumed. 

Table \ref{M} shows the fusion law $\cM$ that is the focus of this 
paper.
\begin{table}[h]
\begin{center}
\begin{tabular}{|c||c|c|c|c|}
\hline
$\ast$&$1$&$0$&$\al$&$\bt$\\
\hline\hline
$1$&$1$&&$\al$&$\bt$\\
\hline
$0$&&$0$&$\al$&$\bt$\\
\hline
$\al$&$\al$&$\al$&$1,0$&$\bt$\\
\hline
$\bt$&$\bt$&$\bt$&$\bt$&$1,0,\al$\\
\hline
\end{tabular}
\end{center}
\caption{Fusion law $\cM(\al,\bt)$}\label{M}
\end{table}
Here $\al,\bt\in\F$ are arbitrary numbers distinct from each other and from 
$1$ and $0$. Each cell of the table lists the elements of the corresponding 
set $\lm\ast\mu$. For example, $1\ast 0=\emptyset$ and $\al\ast\al=\{1,0\}$.
Primitive axial algebras with this fusion law are called \emph{algebras of 
Monster type $(\al,\bt)$}. 

Often an axial algebra admits a bilinear form that associates with the 
algebra product.

\begin{definition}
A \emph{Frobenius form} on an axial algebra $A$ over a field $\F$ is a 
non-zero bilinear form $(\cdot,\cdot)$ on $A$ such that 
$$(uv,w)=(u,vw)$$
for all $u,v,w\in A$.
\end{definition}

The decomposition $A=\oplus_{\lm\in X}A_\lm(a)$ of an axial algebra $A$ 
with respect to an axis $a\in A$ is an orthogonal decomposition, that is, 
$$(A_\lm(a),A_\mu(a))=0$$
for all $\lm,\mu\in X$, $\lm\neq\mu$. 

We say that an axis $a\in A$ is \emph{singular} if $(a,a)=0$ and $a$ is 
\emph{non-singular} otherwise. A primitive axis lies in the radical
$$A^\perp=\{u\in A\mid (u,v)=0\mbox{ for all }v\in A\}$$
if and only if it is singular. 

\subsection{Ideals and simplicity}
\label{ideals}

Methods to find all ideals and, consequently, to check whether $A$ is simple 
were developed in \cite{kms}. Recall that every axial algebra comes with a set 
of generating axes. According to \cite{kms}, all ideals in an axial algebra $A$ 
can be classified into two types: (1) ideals not containing any of the 
generating axes, and (2) ideals containing generating axes. The first type of 
ideals are controlled by the radical of the algebra.

\begin{definition}
The \emph{radical} $R(A)$ of a primitive axial algebra $A$ is the unique maximal 
ideal not containing any of the generating axes of $A$.
\end{definition}

If $A$ admits a Frobenius form then it is easy to see that its radical $A^\perp$ 
is an ideal of $A$. Recall that a primitive axis $a\in A$ is contained in 
$A^\perp$ if and only if $a$ is singular. 

\begin{theorem}[\cite{kms}] \label{perp}
If $A$ is a primitive axial algebra admitting a Frobenius form with respect to 
which none of the generating axes of $A$ is singular then the radical $R(A)$ 
coincides with the radical $A^\perp$ of the Frobenius form.
\end{theorem} 

In particular, if the Frobenius form is non-degenerate (i.e., $A^\perp=0$) then 
$A$ has no non-zero ideals of the first kind. For the ideals of the second type, 
\cite{kms} suggest the following construction. Suppose $a\in A$ is a primitive 
axis. Since $A=\oplus_{\lm\in X}A_\lm(a)$, we can write an arbitrary $u\in A$ 
as $u=\sum_{\lm\in X}u_\lm$, where $u_\lm\in A_\lm(a)$ for each $\lm$. Since $a$ 
is primitive, $u_1=\al a$ for some $\al\in\F$. We call $u_1$ the \emph{projection} 
of $u$ onto $a$. 

\begin{definition} \label{projection graph}
Suppose $A$ is a primitive axial algebra. The \emph{projection graph} $\Dl$ of $A$ 
has as vertices all generating axes of $A$ and $\Dl$ has a directed edge from $b$ 
to $a$ if $b$ has a non-zero projection $b_1$ onto $a$.
\end{definition} 

It can be shown that if $b$ is contained in an ideal $J$ then every component 
$b_\lm$ is contained in $J$. In particular, if $b_1=\al a\neq 0$ then $a\in J$. 
It follows that if $\Dl$ is strongly connected (any vertex can be reached from any 
other vertex via a directed path) then every ideal $J$ of the second type contains 
all generating axes, and so $J=A$.

When $A$ admits a Frobenius form, the definition of $\Dl$ simplifies.

\begin{theorem}[\cite{kms}] \label{unoriented}
Suppose $A$ is a primitive axial algebra admitting a Frobenius form with all 
generating axes non-singular. Then the projection graph $\Dl$ is a simple 
(undirected) graph and $a$ and $b$ are adjacent in $\Dl$ if and only if 
$(a,b)\neq 0$.
\end{theorem}

In particular, if $\Dl$ is connected then $A$ has no proper ideals of the second 
type. The final result of this section summarizes our discussion.

\begin{theorem} \label{summary}
Suppose that a primitive axial algebra admits a Frobenius form and all generating 
axes are non-singular with respect to the form. If the projection graph $\Dl$ of 
$A$ is connected then $A$ is simple if and only if the Frobenius form has zero 
radical.  
\end{theorem}

\subsection{Grading and the Miyamoto group}
\label{Miyamoto}

When the fusion law $\cF$ is graded by an abelian group, every axis in an axial
algebra $A$ leads to automorphisms of $A$.

\begin{definition}
Suppose $\cF=(X,\ast)$ is a fusion law. A \emph{grading} of $\cF$ by a group 
$T$ is a partition $\{X_t|t\in T\}$ of the set $X$ (where we allow parts $X_t$ 
to be empty for some $t\in T$), such that, for all $\lm,\mu\in X$, if $\lm\in X_t$ 
and $\mu\in X_{t'}$ then $\lm\ast\mu$ is contained in $X_{tt'}$. 
\end{definition}

See \cite{dpsv} for a slightly different, categorical definition.

The fusion law $\cM(\al,\bt)$ in Table \ref{M} is graded by the group $C_2=\{1,-1\}$. 
Here $X_1=\{1,0,\al\}$ and $X_{-1}=\{\bt\}$.

Suppose $\cF$ is a fusion law graded by a group $T$ and $A$ is an $\cF$-axial algebra. 
For an axis $a\in A$ and a linear character $\xi$ of $T$, define the linear map 
$\tau_a(\xi): A\to A$ as follows: on each part $A_{X_t}(a)=\oplus_{\lm\in X_t}A_\lm(a)$, 
the map $\tau_a(\xi)$ acts as scalar $\xi(t)$. (That is, $\tau_a(\xi)(u)=\xi(t)u$ for all 
$u\in A_{X_t}(a)$.) It is easy to see that $\tau_a(\xi)$ is an automorphism of $A$ 
and the map $\xi\mapsto\tau_a(\xi)$ is a homomorphism from the group $T^\ast$ of linear 
characters of $T$ to $\Aut(A)$. The image of this homomorphism, 
$T_a=\{\tau_a(\xi)\mid\xi\in T^\ast\}$ is called the \emph{axial subgroup} of $\Aut(A)$ 
corresponding to the axis $a$ and the subgroup of $\Aut(A)$ generated by the axial 
subgroups $T_a$ for all generating axes $a\in A$ is called the \emph{Miyamoto group} of 
$A$.

When $T=C_2$ (for example, in the case of the fusion law $\cM(\al,\bt)$), the group 
$T^\ast$ is of order two (assuming that the characteristic of $\F$ is not $2$). Then 
the only non-identity character $\xi$ of $T$ produces the element $\tau_a=\tau_a(\xi)$, 
known as the \emph{Miyamoto involution}. Here the axial subgroup $T_a=\langle\tau_a\rangle$ 
is of order two (or one if $A_{X_{-1}}(a)=0$) and the Miyamoto group is simply the 
subgroup of $\Aut(A)$ generated by all Miyamoto involutions $\tau_a$.

\section{Matsuo algebras}

\label{Matsuo}

In this section we introduce the family of axial algebras called 
the \emph{Matsuo algebras}. 

Recall that a \emph{$3$-transposition group} is a pair $(G,D)$, 
where $G$ is a group and $D$ is a normal subset (union of conjugacy 
classes) of $G$ such that: 
\begin{enumerate}
\item[(a)] $D$ generates $G$; 
\item[(b)] every $d\in D$ is an involution (that is, $|d|=2$); and 
\item[(c)] for $c,d\in D$, $|cd|\leq 3$.
\end{enumerate}
In many cases, $D$ is unique for a given $G$, so it is common to talk 
about the $3$-transposition group $G$, instead of $(G,D)$.

An example of $3$-transposition group is given by $G=S_n$, the symmetric 
group on $n$ symbols, and $D=(1,2)^G$, the class of transpositions 
($2$-cycles). Finite $3$-transposition groups were classified, under some
restrictions, by Fischer \cite{f} and, in complete generality, by Cuypers 
and Hall \cite{ch}, who used geometric methods.

Suppose $(G,D)$ is a $3$-transposition group. The \emph{Fischer space} of 
$(G,D)$ is a point-line geometry $\Gm=\Gm(G,D)$, whose point set is $D$ 
and where points $c$ and $d$ are collinear if and only if $|cd|=3$. 
Furthermore, any two collinear points $c$ and $d$ lie in a 
unique common line, which consists of $c$, $d$, and the third point 
$e=c^d=d^c$. Note that $c$, $d$, and $e$ are the three involutions in the 
subgroup $S_3$ that any two of them generate. Thus, lines of $\Gm$ are in 
a bijection with the subgroups $S_3$ generated by the elements of $D$. 

It is easy to see that the connected components of the Fischer space 
$\Gm$ coincide with the conjugacy classes of $G$ contained in $D$. 
In particular, the Fischer space is connected if and only if $D$ is a 
single conjugacy class of $G$. We will also say in this case that $(G,D)$ 
is \emph{connected}. 

The $3$-transposition group $(G,D)$ can be recovered from $\Gm$ up to the 
center of $G$, which, clearly, acts trivially on $\Gm$. The further 
discussion will be mainly in terms of the Fischer space $\Gm$, but we will 
keep using $D$ for the point set of $\Gm$. 

Let us associate an algebra with the Fischer space $\Gm$. Select a field 
$\F$ of characteristic not equal to $2$ and let $\eta\in\F$, $\eta\neq 
0,1$. 

\begin{definition}
The \emph{Matsuo algebra} $M_\eta(\Gm)$ over $\F$, corresponding to $\Gm$ 
and $\eta$, has the point set $D$ as its basis. Multiplication is defined 
on the basis as follows:
$$c\cdot d=\left\{
\begin{array}{rl}
c,&\mbox{if }c=d;\\
0,&\mbox{if $c$ and $d$ are non-collinear};\\
\frac{\eta}{2}(c+d-e),&\mbox{if $c$ and $d$ lie on a line, $\{c,d,e\}$}.
\end{array}
\right.$$ 
\end{definition}

Here we use the dot for the algebra product to distinguish it from 
the multiplication in the $3$-transposition group $G$. In what follows 
we skip the dot as long as this causes no confusion. 

The following small cases will appear prominently in the remainder of the 
paper, so we need special names for them, coming from \cite{hrs} and 
generalizing the names of the Norton-Sakuma algebras \cite{ipss}. First of 
all, if the Fischer space $\Gm$ consists of a single point, say $c$, then 
$M\cong\F$ is referred to as the algebra $1A$. If $\Gm$ consists of two 
non-collinear points then $M\cong\F^2$ and it is referred to as the algebra 
$2B$. Finally, if $\Gm$ consists of three points forming a line then $M$ is 
$3$-dimensional called $3C(\eta)$. (The value of $\eta$ is irrelevant for $1A$ 
and $2B$.) 

Let us record now some properties of Matsuo algebras.

\begin{proposition}
If $\Gm$ is disconnected with components $\Gm_i$, $i\in I$, then 
$M_\eta(\Gm)=\oplus_{i\in I}M_\eta(\Gm_i)$.
\end{proposition}

Included in this statement is the property that 
$M_\eta(\Gm_i)M_\eta(\Gm_j)=0$ for all $i,j\in I$, $i\neq j$.

Clearly, each $c\in D$ is an idempotent in $M_\eta(\Gm)$. Furthermore, it 
can be shown that $c$ satisfies the fusion law $\cJ(\eta)$ in Table 
\ref{Jordan type}.
\begin{table}[h]
\begin{center}
\begin{tabular}{|c||c|c|c|}
\hline
$\ast$&$1$&$0$&$\eta$\\
\hline\hline
$1$&$1$&&$\eta$\\
\hline
$0$&&$0$&$\eta$\\
\hline
$\eta$&$\eta$&$\eta$&$1,0$\\
\hline
\end{tabular}
\end{center}
\caption{Fusion law $\cJ(\eta)$}\label{Jordan type}
\end{table}
This means that Matsuo algebras are \emph{algebras of Jordan type $\eta$}. 
This class of axial algebras was introduced in \cite{hrs}, where a 
partial classification, for $\eta\neq\frac{1}{2}$, was achieved. (See also a 
correction and further details in \cite{hss}.)

\begin{proposition}[\cite{hrs,hss}] \label{characterization}
Any algebra of Jordan type $\eta\neq\frac{1}{2}$ is either a Matsuo 
algebra or a factor algebra of a Matsuo algebra.
\end{proposition} 

The Matsuo algebra admits a Frobenius form $(\cdot,\cdot)$ such that 
$(c,c)=1$ for each axis $c\in D$. This form is given on the basis $D$ 
by the following:
$$
(c,d)=\left\{
\begin{array}{rl}
1,&\mbox{ if }c=d;\\
0,&\mbox{ if }|cd|=2;\\
\frac{\eta}{2},&\mbox{ if }|cd|=3.
\end{array}\right.
$$

The fusion law $\cJ(\eta)$ is $C_2$-graded, with $X_1=\{1,0\}$ and 
$X_{-1}=\{\eta\}$. The action of the Miyamoto involution $\tau_c$ on 
$M_\eta(\Gm)$ agrees with the action of the $3$-transposition $c$ on $D$ 
by conjugation. Correspondingly, the Miyamoto group of $M_\eta(\Gm)$ is 
isomorphic to the factor group $G/Z(G)$ of the $3$-transposition group 
$G$ (from which the Fischer space $\Gm$ was obtained) over its center. 

Another way to describe the action of the Miyamoto involution $\tau_c$ is in 
terms of the Fischer space $\Gm$. Namely, $\tau_c$ fixes $c$ and all points 
non-collinear with $c$ and it switches the two points other than $c$ on any 
line through $c$. This provides a very efficient way of computing the image 
under $\tau_c$.

\section{Double axes}

\label{double}

Suppose $M=M_\eta(\Gm)$ is the Matsuo algebra for the Fischer space $\Gm$ 
corresponding to a $3$-transposition group $(G,D)$. 

\begin{definition}
Axes $a,b\in D$ are called \emph{orthogonal} if $ab=0$. In this 
case we call the idempotent $x:=a+b$ a \emph{double axis}.
\end{definition}

Indeed, $(a+b)^2=a^2+ab+ba+b^2=a+0+0+b=a+b$, and so $a+b$ is an idempotent. 

To distinguish double axes from the Matsuo axes in $D$, we will sometimes 
call the latter \emph{single axes}. Thus, single axes satisfy the fusion law 
$\cJ(\eta)$ and each double axis is made out of two orthogonal single axes.

The goal of this section is to show that double axes satisfy a specific nice 
fusion law and so they are indeed nice axes in the sense of axial algebras. 
Throughout this section, $a,b\in D$ are orthogonal single axes and $x=a+b$.

Note that a linear transformation $\phi$ is semisimple if and only if its 
minimal polynomial, say $f$, has no multiple roots. If $U$ is a subspace 
invariant under $\phi$ then, clearly, $f(\phi)=0$ on $U$. In other 
words, the minimal polynomial of $\phi$ in its action on $U$ is a divisor 
of $f$, and in particular, it cannot have multiple roots. Hence $\phi$ is 
also semisimple on each invariant subspace $U$. We apply this to 
$\phi=\ad_b$ acting on the eigenspaces of $\ad_a$. First of all, note that 
the eigenspaces $M_\lm(a)$ are indeed invariant under $\ad_b$. This is 
because $b\in M_0(a)$ and, according to the fusion law $\cJ(\eta)$, we have
$0\ast\lm\subseteq\{\lm\}$ for all $\lm\in X=\{1,0,\eta\}$. Hence 
$\ad_b(M_\lm(a))=bM_\lm(a)\subseteq M_\lm(a)$. 

We use the following notation: for $\lm,\mu\in\{1,0,\eta\}$, 
$M_{\lm,\mu}:=M_\lm(a)\cap M_\mu(b)$.

\begin{lemma}
For $\lm\in\{1,0,\eta\}$, $M_\lm(a)=\oplus_{\mu\in\{1,0,\eta\}}M_{\lm,\mu}$.
\end{lemma}

\begin{proof}
We have already seen that $M_\lm(a)$ is invariant under $\ad_b$. Since $\ad_b$ 
is semisimple on $A$, it is also semisimple on $M_\lm(a)$, and so $M_\lm(a)$ 
decomposes as the direct sum of eigenspaces with respect to $\ad_b$. Manifestly, 
$M_{\lm,\mu}=M_\lm(a)\cap M_\mu(b)$ are the only possible nontrivial eigenspaces 
of $\ad_b$ in $M_\lm(a)$. This yields the desired decomposition.
\end{proof}

Since $M=\oplus_{\lm\in\{1,0,\eta\}}M_\lm(a)$, we also obtain the following. 

\begin{corollary} \label{fine}
$M=\oplus_{\lm,\mu\in\{1,0,\eta\}}M_{\lm,\mu}$.
\end{corollary}

We next observe that $\ad_x$ acts on each $M_{\lm,\mu}$ as a scalar.

\begin{lemma} \label{value}
For $u\in M_{\lm,\mu}$, we have $xu=(\lm+\mu)u$.
\end{lemma}

\begin{proof}
Indeed, $xu=(a+b)u=au+bu=\lm u+\mu u=(\lm+\mu)u$, since $u$ lies in 
$M_\lm(a)$ and in $M_\mu(b)$.
\end{proof}

At this point, we already see that $\ad_x$ is semisimple on $M$. Next, 
we investigate which eigenvalues arise and identify the eigenspaces. 

In total there are $3^2=9$ pieces $M_{\lm,\mu}$. We first note that three 
of them are trivial and further two are just the $1$-dimensional spaces 
$\langle a\rangle$ and $\langle b\rangle$.

\begin{lemma}
\begin{itemize}
\item[(a)] $M_{1,1}=M_{1,\eta}=M_{\eta,1}=0$; and 
\item[(b)] $M_{1,0}=M_1(a)=\langle a\rangle$ and $M_{0,1}=M_1(b)=\langle 
b\rangle$.
\end{itemize}
\end{lemma}

\begin{proof}
Since $M_1(a)=\langle a\rangle$ is $1$-dimensional and since $a\in M_0(b)$, 
we conclude that $M_{1,1}=0$, $M_{1,0}=\langle a\rangle$, and $M_{1,\eta}=0$. 
Similarly, since $M_1(b)=\langle b\rangle$ and since $b\in M_0(a)$, we also 
have that $M_{0,1}=\langle b\rangle$, and $M_{\eta,1}=0$.
\end{proof}

So in total we have at most six nontrivial pieces $M_{\lm,\mu}$. By Lemma 
\ref{value}, $M_{1,0}\oplus M_{0,1}\subseteq M_1(x)$, $M_{0,0}\subseteq 
M_0(x)$, $M_{\eta,\eta}\subseteq M_{2\eta}(x)$, and $M_{0,\eta}\oplus 
M_{\eta,0}\subseteq M_\eta(x)$. Combining this with Corollary \ref{fine}, 
yields that all these inclusions are equalities. Hence we have the 
following.

\begin{proposition} \label{double eigenspaces}
$M=M_1(x)\oplus M_0(x)\oplus M_{2\eta}(x)\oplus M_\eta(x)$, where
\begin{itemize} 
\item[(a)] $M_1(x)=M_{1,0}\oplus M_{0,1}=\langle a,b\rangle$;
\item[(b)] $M_0(x)=M_{0,0}$;
\item[(c)] $M_{2\eta}(x)=M_{\eta,\eta}$; and
\item[(d)] $M_\eta(x)=M_{\eta,0}\oplus M_{0,\eta}$.
\end{itemize}
\end{proposition}

From this point and until the end of the paper we assume that 
$\eta\neq\frac{1}{2}$. This is in addition to our earlier assumption that 
$\eta\neq 1,0$. This guarantees that $2\eta\not\in\{1,0,\eta\}$.

We have established that the eigenvalues of $\ad_x$ are contained in 
$\{1,0,2\eta,\eta\}$. It remains to see which fusion law is satisfied. 
\begin{table}[h]
\begin{center}
\begin{tabular}{|c||c|c|c|c|}
\hline
$\diamond$&$1$&$0$&$2\eta$&$\eta$\\
\hline\hline
$1$&$1$&&$2\eta$&$\eta$\\
\hline
$0$&&$0$&$2\eta$&$\eta$\\
\hline
$2\eta$&$2\eta$&$2\eta$&$1+0$&$\eta$\\
\hline
$\eta$&$\eta$&$\eta$&$\eta$&$1+0+2\eta$\\
\hline
\end{tabular}
\end{center}
\caption{Fusion law for the double axis $x$}\label{double law}
\end{table} 

\begin{proposition} \label{double fusion rules}
The double axis $x$ satisfies the fusion law $\cM(\al,\bt)$ with 
$\al=2\eta$ and $\bt=\eta$ (see Table \ref{double law}).
\end{proposition}

Note that we consider the fusion law for $x$ in parallel with 
the Jordan type fusion law for $a$ and $b$. Hence we use a different 
symbol, $\diamond$, in the table for $x$. 

We start with two lemmas.

\begin{lemma} \label{observation 1}
$M_{\lm,\mu}M_{\gm,\dl}\subseteq M_{\lm\ast\gm}(a)\cap M_{\mu\ast\dl}(b)$.
\end{lemma}

\begin{proof}
Indeed, 
\begin{align*}
M_{\lm,\mu}M_{\gm,\dl}&=(M_\lm(a)\cap M_\mu(b))(M_\gm(a)\cap M_\dl(b))\\
&\subseteq(M_\lm(a)M_\gm(a))\cap(M_\mu(b)M_\dl(b))\\
&\subseteq M_{\lm\ast\gm}(a)\cap M_{\mu\ast\dl}(b).
\end{align*}
\end{proof}

The second lemma is an extension of Corollary \ref{fine}.

\begin{lemma} \label{observation 2}
Suppose $S,T\subseteq\{1,0,\eta\}$. Then $M_S(a)\cap M_T(b)=\oplus_{\lm\in S,\mu\in T}M_{\lm,\mu}$.
\end{lemma}

\begin{proof}
Clearly, $M_S(a)\cap M_T(b)\supseteq\oplus_{\lm\in S,\mu\in T}M_{\lm,\mu}$. 
Hence we just need to show the reverse inclusion.

Define $g\in\F[z]$ to be the product $\prod_{\lm\in S}(z-\lm)$. Then 
$u\in M$ lies in $M_S(a)$ if and only if $g(\ad_a)u=0$. Let $u\in M_S(a)\cap 
M_T(b)$ and write it as $u=\sum_{\mu\in T}u_\mu$, where each $u_\mu$ lies in 
$M_\mu(b)$.

Note that $0=g(\ad_a)u=\sum_{\mu\in T}g(\ad_a)u_\mu$ and that each component 
$g(\ad_a)u_\mu$ is contained in the corresponding $M_\mu(b)$, since the 
latter is invariant under $\ad_a$. It follows that $g(\ad_a)u_\mu=0$ for 
each $\mu\in T$, that is, $u_\mu\in M_S(a)\cap M_\mu(b)$. 

Now we repeat this argument switching the r\^oles of $a$ and $b$. Namely, 
we decompose $u_\mu$ as $u_\mu=\sum_{\lm\in S} u_{\lm,\mu}$, where 
$u_{\lm,\mu}\in M_\lm(a)$. Taking $h\in\F[z]$ equal to $z-\mu$, we see that 
$0=h(\ad_b)u_\mu=\sum_{\lm\in S}h(\ad_b)u_{\lm,\mu}$ and, as above, we 
conclude that $h(\ad_b)u_{\lm,\mu}=0$ for all $\lm\in S$, that is, 
$u_{\lm,\mu}\in M_\lm(a)\cap M_\mu(b)=M_{\lm,\mu}$. Thus, $u=\sum_{\mu\in 
T}u_\mu=\sum_{\lm\in S,\mu\in T}u_{\lm,\mu}\in\oplus_{\lm\in S,\mu\in 
T}M_{\lm,\mu}$. We have shown that the reverse inclusion does indeed hold.
\end{proof}

We can now proceed with the proof of Proposition \ref{double fusion rules}.

\begin{proof}
We have already established that $\ad_x$ has eigenvalues in 
$\{1,0,2\eta,\eta\}$. Since $M$ is commutative, we only need to check the 
upper triangular part of the fusion law.

We start with the first row of the table:
\begin{itemize}
\item $M_1(x)M_1(x)=\langle a,b\rangle\langle a,b\rangle\subseteq\langle 
aa,ab,bb\rangle=\langle a,b\rangle=M_1(x)$, since $aa=a$, $ab=0$, and 
$bb=b$. Thus, $1\diamond 1=\{1\}$.
\item $M_1(x)M_0(x)=\langle a,b\rangle M_{0,0}\subseteq aM_{0,0}+ 
bM_{0,0}\subseteq aM_0(a)+bM_0(b)=0+0=0$. Thus, $1\diamond 0=\emptyset$. 
\item $M_1(x)M_{2\eta}(x)=\langle a,b\rangle M_{\eta,\eta}\subseteq 
aM_{\eta,\eta}+bM_{\eta,\eta}=M_{\eta,\eta}+M_{\eta,\eta}=M_{\eta,\eta}= 
M_{2\eta}(x)$, since $M_{\eta,\eta}=M_\eta(a)\cap M_\eta(b)$ and so 
$aM_{\eta,\eta}=M_{\eta,\eta}=bM_{\eta,\eta}$. Thus, $1\diamond 2\eta= 
\{2\eta\}$.
\item $M_1(x)M_\eta(x)=\langle a,b\rangle(M_{\eta,0}+M_{0,\eta})\subseteq 
aM_{\eta,0}+aM_{0,\eta}+bM_{\eta,0}+bM_{0,\eta}=M_{\eta,0}+0+0+M_{0,\eta}= 
M_{\eta,0}+M_{0,\eta}=M_\eta(x)$. Thus, $1\diamond\eta=\{\eta\}$.
\end{itemize}

We now turn to the second row and here we will start using Lemma 
\ref{observation 1}:
\begin{itemize}
\item $M_0(x)M_0(x)=M_{0,0}M_{0,0}\subseteq M_{0\ast 0}(a)\cap M_{0\ast 
0}(b)=M_0(a)\cap M_0(b)=M_{0,0}=M_0(x)$. Thus, $0\diamond 0=\{0\}$.
\item $M_0(x)M_{2\eta}(x)=M_{0,0}M_{\eta,\eta}\subseteq 
M_{0\ast\eta}(a)\cap M_{0\ast\eta}(b)=M_\eta(a)\cap M_\eta(b)= 
M_{\eta,\eta}=M_{2\eta}(x)$. Hence, $0\diamond 2\eta=\{2\eta\}$.
\item $M_0(x)M_\eta(x)=M_{0,0}(M_{\eta,0}+M_{0,\eta})\subseteq 
M_{0,0}M_{\eta,0}+M_{0,0}M_{0,\eta}\subseteq (M_{0\ast\eta}(a)\cap 
M_{0\ast 0}(b))+(M_{0\ast 0}(a)\cap M_{0\ast\eta}(b))=(M_\eta(a)\cap 
M_0(b))+(M_0(a)\cap M_\eta(b))=M_{\eta,0}+M_{0,\eta}=M_\eta(x)$. We have 
shown that $0\diamond\eta=\{\eta\}$.
\end{itemize}

Next, the third row:
\begin{itemize}
\item $M_{2\eta}(x)M_{2\eta}(x)=M_{\eta,\eta}M_{\eta,\eta}\subseteq 
M_{\eta\ast\eta}(a)\cap M_{\eta\ast\eta}(b)=(M_1(a)+M_0(a))\cap 
(M_1(b)+M_0(b))=M_{1,1}+M_{1,0}+M_{0,1}+M_{0,0}=M_1(x)+M_0(x)$. (Here we 
used Lemma \ref{observation 2} for the first time.) Therefore, 
$2\eta\diamond 2\eta=\{1,0\}$.
\item $M_{2\eta}(x)M_\eta(x)=M_{\eta,\eta}(M_{\eta,0}+M_{0,\eta})\subseteq 
M_{\eta,\eta}M_{\eta,0}+M_{\eta,\eta}M_{0,\eta}\subseteq(M_{\eta\ast\eta}(a)
\cap M_{\eta\ast 0}(b))+(M_{\eta\ast 0}(a)\cap M_{\eta\ast\eta}(b))=
((M_1(a)+M_0(a))\cap M_\eta(b))+(M_\eta(a)\cap(M_1(b)+M_0(b)))=
M_{1,\eta}+M_{0,\eta}+M_{\eta,1}+M_{\eta,0}=M_{0,\eta}+M_{\eta,0}= 
M_\eta(x)$. Thus, $2\eta\diamond\eta=\{\eta\}$.
\end{itemize}

Finally, there is only one entry in the fourth row for us to check:
\begin{itemize}
\item $M_\eta(x)M_\eta(x)=(M_{\eta,0}+M_{0,\eta})(M_{\eta,0}+ 
M_{0,\eta})\subseteq M_{\eta,0}M_{\eta,0}+M_{\eta,0}M_{0,\eta}+ 
M_{0,\eta}M_{\eta,0}+M_{0,\eta}M_{0,\eta}\subseteq (M_{\eta\ast\eta}(a)\cap 
M_{0\ast 0}(b))+(M_{\eta\ast 0}(a)\cap M_{0\ast\eta}(b))+ 
(M_{0\ast\eta}(a)\cap M_{\eta\ast 0}(b))+(M_{0\ast 0}(a)\cap 
M_{\eta\ast\eta}(b))=((M_1(a)+M_0(a))\cap M_0(b))+(M_\eta(a)\cap 
M_\eta(b))+(M_\eta(a)\cap M_\eta(b))+(M_0(a)\cap (M_1(b)+M_0(b)))= 
M_{1,0}+M_{0,0}+M_{\eta,\eta}+M_{\eta,\eta}+M_{0,1}+ M_{0,0}= 
M_1(x)+M_0(x)+M_{2\eta}(x)$. Therefore, $\eta\diamond\eta=\{1,0,2\eta\}$.
\end{itemize}

This completes the proof of the proposition.
\end{proof}

The fusion law $\cM(2\eta,\eta)$, like any law $\cM(\al,\bt)$, is graded by 
the group $C_2$. Hence for each double axis $x\in M$ we have the corresponding 
Miyamoto involution $\tau_x$.

\begin{proposition} \label{Miyamoto involution}
For a double axis $x=a+b$, we have $\tau_x=\tau_a\tau_b$.
\end{proposition}

\begin{proof}
The involution $\tau_a$ acts as identity on $M_1(a)+M_0(a)$ and as minus 
identity on $M_\eta(a)$. Similarly, $\tau_b$ acts as identity on 
$M_1(b)+M_0(b)$ and as minus identity on $M_\eta(b)$. Also, $\tau_x$ acts as 
identity on $M_1(x)+M_0(x)+M_{2\eta}(x)$ and as minus identity on $M_\eta(x)$. 

Let us compare the action of $\tau_x$ and $\tau_a\tau_b$. Clearly, $\tau_x$, 
$\tau_a$ and $\tau_b$ act as identity on $M_1(x)+M_0(x)=
M_{1,0}+M_{0,1}+M_{0,0}$ and so the actions agree here. On 
$M_{2\eta}(x)=M_{\eta,\eta}$, $\tau_x$ acts as identity, while both $\tau_a$ 
and $\tau_b$ act as minus identity. However, this means that $\tau_a\tau_b$ 
acts as identity, and so $\tau_x$ and $\tau_a\tau_b$ agree on $M_{2\eta}(x)$. 
Finally, $\tau_x$ acts as minus identity on $M_\eta(x)=
M_{\eta,0}+M_{0,\eta}$. On the first summand, $\tau_a$ acts as minus identity 
and $\tau_b$ as identity. Hence, $\tau_a\tau_b$ acts as minus identity, 
agreeing with $\tau_x$. Similarly, on the second summand, $\tau_a$ acts as 
identity and $\tau_b$ as minus identity. Hence again $\tau_a\tau_b$ agrees 
with $\tau_x$.     
\end{proof} 

Note that double axes are not primitive, as $M_1(x)=\langle a,b\rangle$ is 
$2$-dimensional. However, $x$ may be primitive within a proper subalgebra 
of the algebra $M$. We will explore this idea in the remaining sections.

\section{2-Generated subalgebras}

\label{2-generated}

Manifestly, the Jordan type fusion law $\cJ(\eta)$ is a minor (sublaw) of 
the Monster type fusion law $\cM(2\eta,\eta)$. This means that single axes
satisfy the same fusion law $\cM(2\eta,\eta)$ as double axes, except that 
the $2\eta$-eigenspace for them is trivial. Therefore, any subalgebra
$A$ of $M$ that is generated by single and double axes is an algebra of 
Monster type $(2\eta,\eta)$ and it is primitive as long as the double axes 
we use are primitive in $A$. The principal goal of this section is to find 
all primitive subalgebras in Matsuo algebras $M$, generated by two axes, 
at least one of which is a double axis.

As in Section \ref{double}, $M=M_\eta(\Gm)$ is a Matsuo algebra coming from a 
Fischer space $\Gm$, and $(G,D)$ is the $3$-transposition group leading to 
$\Gm$. Recall from Section \ref{Matsuo} that $M$ admits a Frobenius form 
$(\cdot,\cdot)$ and, clearly, this form is inherited by any subalgebra  
$A\subseteq M$.

\subsection{Flip}

We start by showing that primitivity of the subalgebra $A$ forces a certain 
symmetry in the underlying set of single axes. 

\begin{definition}
The support $\supp(x)$ of a single axis $x=a$ is the set $\{a\}$ and the 
support $\supp(y)$ of a double axis $y=b+c$ is the set $\{b,c\}$. More 
generally, the \emph{support} $\supp(Y)$ of a set $Y$ of single and double 
axes is the union of supports of all axes from $Y$.
\end{definition}

Manifestly, $\supp(Y)$ is a subset of the basis $D$ of $M$. By the 
\emph{diagram} on the set $Z\subseteq D$ we mean the graph (denoted by 
the same symbol $Z$) whose vertex set is $Z$ and where two vertices are 
adjacent if and only if they are collinear as points of the Fischer space 
$\Gm$.

In what follows, we use double angular brackets $\dla~\dra$ to indicate 
subalgebra generation, leaving single brackets for the linear span. 

\begin{proposition} \label{same value}
Suppose that $x=a+b$ is a double axis and $y$ is any axis. If $x$ is primitive 
in $A=\dla x,y\dra$ then $(a,y)=(b,y)$.
\end{proposition}

\begin{proof}
Decompose $y$ with respect to $\ad_x$ as $y=y_1+y_0+y_{2\eta}+y_\eta$, where 
each $y_\lm$ is contained in $M_\lm(x)$. In particular, $y_1\in M_1(x)=
\la a,b\ra$. Hence $y_1=\al a+\bt b$. Since $a\in M_1(x)$ and the eigenspaces 
$M_\lm(x)$ are pairwise orthogonal, we deduce that $(a,y)=(a,y_1)=
(a,\al a+\bt b)=\al(a,a)+\bt(a,b)=\al$, since $(a,a)=1$ and $(a,b)=0$ (because 
$b\in M_0(a)$). Similarly, $(b,y)=\bt$. 

Now recall that every component $y_\lm$ is contained in every subalgebra 
containing $x$ and $y$. In particular, $y_1\in A_1(x)$. If $x$ is primitive 
then $A_1(x)=\la x\ra=\la a+b\ra$. This means that $y_1=\al a+\bt b$ is a 
multiple of $a+b$, and so $\al=\bt$, implying that $(a,y)=(b,y)$.
\end{proof}

\begin{corollary} \label{disjoint}
Under the assumptions of Proposition \ref{same value}, if $x\neq y$ and 
$x$ is primitive in $A=\dla x,y\dra$ then $\supp(x)\cap\supp(y)=\emptyset$ 
and, furthermore, $a$ and $b$ are adjacent in the diagram on $Z=\supp(\{x,y\})$ 
to the same number of vertices from $\supp(y)$. 
\end{corollary}

\begin{proof}
Indeed, if $u\in\{a,b\}$ is contained in $\supp(y)$ then $(u,y)=1$; if $u$ 
is adjacent to one single axis in $\supp(y)$ then $(u,y)=\frac{\eta}{2}$; 
and if $u$ is adjacent to two single axes in $\supp(y)$ (when $y$ is a double 
axis) then $(u,y)=\eta$. This follows from the description of the Frobenius 
form on $M$ given in Section \ref{Matsuo}.

Finally, we remark that if both $a$ and $b$ are contained in $\supp(y)$ then 
$x=y$, a contradiction.
\end{proof}

We can now show the symmetry forced by primitivity.

\begin{definition} \label{flip} 
Let $Y$ be a set of single and double axes such that $\dla Y\dra$ is 
primitive. The \emph{flip} is the permutation of $Z=\supp(Y)$ that fixes every 
single axis from $Y$ and switches the two single axes in the support of each 
double axis from $Y$.
\end{definition}

Note that by Corollary \ref{disjoint} different axes from $Y$ have disjoint 
supports. So the flip is well-defined.

\begin{theorem} \label{symmetry}
If the subalgebra $A=\dla Y\dra$ is primitive then the flip is an 
automorphism of the diagram on $Z=\supp(Y)$.
\end{theorem}

\begin{proof}
Consider a double axis $x=a+b\in Y$ and let $y$ be any other axis from $Y$. 
According to Corollary \ref{disjoint}, $a$ and $b$ are adjacent in the diagram 
to the same number $s\in\{0,1,2\}$ of single axes from $\supp(y)$. If $s=0$ or 
$2$ then clearly the flip preserves edges between $\supp(x)$ and $\supp(y)$. 
If $s=1$ then the only non-symmetric situation is where $y$ is a double axis 
and $a$ and $b$ are adjacent to the same single axis in $\supp(y)$. However, 
switching the r\^oles of $x$ and $y$, we see that this situation is impossible.
\end{proof}

Clearly, this theorem restricts the diagrams that can lead to primitive 
subalgebras. We now turn to $2$-generated primitive subalgebras.

\subsection{Setup}

Let $A\subseteq M$ be generated by two axes $x$ and $y$. When both $x$ and 
$y$ are single, we are in the Matsuo algebra situation and so we know what 
$A=\dla x,y\dra$ is; namely, $A\cong 2B$ or $3C(\eta)$. 

Thus, we focus on the situation, where at least one axis is a double axis. 
There are two cases to consider: (A) $x=a$ is a single axis and $y=b+c$ is a 
double axis; and (B) both $x=a+b$ and $y=c+d$ are double axes. 

Let $(H,C)$ be the $3$-transposition group, where $H=\la\supp(\{x,y\})\ra$ and 
so $H=\la a,b,c\ra$ in case (A) and $H=\la a,b,c,d\ra$ in case (B). The normal 
set of involutions $C$ is given by $C=\{a^H\}\cup\{b^H\}\cup\{c^H\}$ and
$C=\{a^H\}\cup\{b^H\}\cup\{c^H\}\cup\{d^H\}$ in the respective cases.
Let $\Sg$ be the Fischer space of $(H,C)$. It is clear that the subalgebra
$A=\dla x,y\dra$ is fully contained in the Matsuo subalgebra $M_\eta(\Sg)$ of 
$M$. In view of this, we can assume that $\Gm=\Sg$ and so $(G,D)=(H,C)$.

Then the group $G$, being generated by at most four $3$-transpositions, belongs
to a short list of small groups. They have been classified by Fischer \cite{f} 
for three generators and by Zara \cite{z} for four generators. (Also available 
is the list of $5$-generated $3$-transposition groups due to Hall and Soicher 
\cite{hs}.) We will look at all possible configurations, organizing them 
according to the diagram on the set $Z=\supp(\{x,y\})$. In view of Theorem 
\ref{symmetry}, we only need to consider diagrams that admit the flip symmetry.

The diagram represents a Coxeter group $\hat G$ given by the presentation 
encoded in the diagram. Clearly, $G$ is a factor group of $\hat G$ and, in 
many cases, $G=\hat G$. These will be the easier cases, and in the remaining 
harder cases we will have to find additional relations identifying $G$ as a 
factor group of the Coxeter group $\hat G$.

Lastly, where possible, we will embed the group $G$ into the symmetric
group $S_n$, so that the involutions $a,b,\ldots$ fall into the class
of transpositions. This embeds the Fischer space $\Gm$ into the Fischer 
space of $S_n$ and it will significantly simplify our notation. 

\subsection{Type A}

\begin{figure}[ht]
\begin{center}
\includegraphics[scale=0.5]{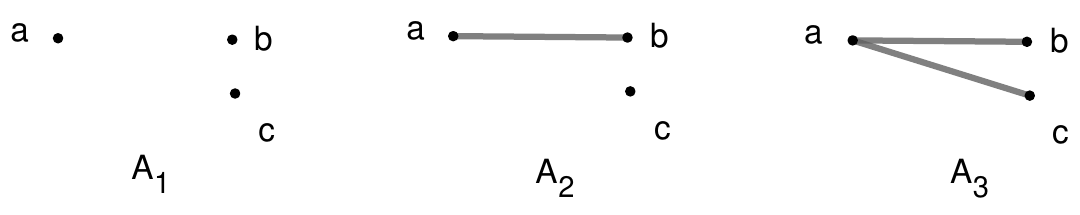}
\caption{A single axis and a double axis}\label{Type A}
\end{center}
\end{figure}

Here we deal with case (A). Recall that $x=a$, $y=b+c$, $Z=\{a,b,c\}$,
$H=\la Z\ra=\la a,b,c\ra$. Since $y=b+c$ is a double axis, the axes $b$ and 
$c$ satisfy $bc=0$, that is, they are not collinear in the Fischer space 
$\Gm$. Thus, we only have diagrams on $Z$ shown in Figure \ref{Type A}.

The second diagram, $A_2$, is not symmetric with respect to the flip, hence we 
discard it. Let us look at each of the remaining two diagrams in turn.

\bigskip\noindent
{\bf Diagram $\mathbf{A_1}$:} Here the Coxeter group $\hat G$ is an
elementary abelian $2$-group of rank $3$ and $G=\hat G$ or $G$ is the
unique factor of $\hat G$ of order $2^2$, in which the three generators
of $\hat G$ remain distinct. In either case, $xy=a(b+c)=ab+ac=0+0=0$,
so $A=2B\cong\F^2$ and $y$ is primitive in $A$.

\bigskip\noindent
{\bf Diagram $\mathbf{A_3}$:} Here $\hat G\cong S_4$ and $G=\hat G$, as $G$ 
cannot be isomorphic to the factor groups $S_3$ and $C_2$ of $S_4$. We 
identify $G$ with $S_4$ by equating $a$ with $(1,2)$, $b$ with $(1,3)$ and $c$ 
with $(2,4)$.

The algebra $A$ is invariant under $\tau_x=\tau_a$, acting on $D$ as 
conjugation by $(1,2)$ (c.f. the discussion in the penultimate paragraph of 
Section \ref{Matsuo}), and under $\tau_y=\tau_b\tau_c$, acting as conjugation 
by $(1,3)(2,4)$ (see Proposition \ref{Miyamoto involution}). Hence 
$z:=x^{\tau_y}=(1,2)^{(1,3)(2,4)}=(3,4)$ is in $A$ and also 
$t:=y^{\tau_x}=((1,3)+(2,4))^{(1,2)}=(2,3)+(1,4)=(1,4)+(2,3)\in A$. We 
claim that $\{x,z,y,t\}$ is a basis of $A$. Manifestly, these four vectors are 
linearly independent in $M$ as their supports are disjoint. Hence we just need 
to check that their span is closed for the algebra product. Throughout the 
paper and especially in the multiplication tables, it will be convenient to 
use the notation $s_i$ for single axes and $d_i$ for double axes. Here we set 
$s_1:=x=(1,2)$, $s_2:=z=(3,4)$, $d_1:=y=(1,3)+(2,4)$ and $d_2:=t=(1,4)+(2,3)$. 

Clearly, all four axes are idempotents and $s_1s_2=0$. Also,
$s_1d_1=(1,2)((1,3)+(2,4))=\frac{\eta}{2}((1,2)+(1,3)-(2,3))+
\frac{\eta}{2}((1,2)+(2,4)-(1,4))=\eta s_1+\frac{\eta}{2}d_1- 
\frac{\eta}{2}d_2=\frac{\eta}{2}(2s_1+d_1-d_2)$. Applying the automorphisms 
$\tau_x$ (switching $d_1$ and $d_2$) and $\tau_y$ (switching $s_1$ and $s_2$), 
we immediately obtain the values of $s_1d_2$, $s_2d_1$, and $s_2d_2$. Finally, 
$d_1d_2=((1,3)+(2,4))((1,4)+(2,3))=\frac{\eta}{2}((1,3)+(1,4)-(3,4))+
\frac{\eta}{2}((1,3)+(2,3)-(1,2))+\frac{\eta}{2}((2,4)+(1,4)-(1,2))+
\frac{\eta}{2}((2,4)+(2,3)-(3,4))=\eta(-s_1-s_2+d_1+d_2)$. Thus, indeed $A$ is 
$4$-dimensional with basis $\{s_1,s_2,d_1,d_2\}$ and the multiplication table 
is as in Table \ref{4-dim}.
\begin{table}[ht]
\begin{center}
\begingroup
\setlength{\tabcolsep}{6pt}
\scalebox{.85}{
$\begin{tabu}[ht]{|c||c|c|c|c|}
\hline
&s_1&s_2&d_1&d_2\\
\hline\hline
s_1&s_1&0&\frac{\eta}{2}(2s_1+d_1-d_2)&\frac{\eta}{2}(2s_1+d_2-d_1)\\
\hline
s_2&0&s_2&\frac{\eta}{2}(2s_2+d_1-d_2)&\frac{\eta}{2}(2s_2+d_2-d_1)\\
\hline
d_1&\frac{\eta}{2}(2s_1+d_1-d_2)&\frac{\eta}{2}(2s_2+d_1-d_2)&d_1&
\eta(-s_1-s_2+d_1+d_2)\\
\hline
d_2&\frac{\eta}{2}(2s_1+d_2-d_1)&\frac{\eta}{2}(2s_2+d_2-d_1)&
\eta(-s_1-s_2+d_1+d_2)&d_2\\
\hline
\end{tabu}$}
\caption{The $4$-dimensional algebra $Q_2(\eta)$}\label{4-dim}
\endgroup
\end{center}
\end{table}
Clearly, $b,c\not\in A$, so $d_1=y=b+c$ is primitive in $A$. We will denote 
this algebra $Q_2(\eta)$. It is part of an infinite series developed in 
Section \ref{symmetric group}.

Let us briefly discuss the properties of this new algebra. We will see 
in Section \ref{symmetric group} that the algebras $Q_k(\eta)$ are simple 
except when $\eta$ belongs to a short list of exceptional values. In this 
case, the only such value is $-\frac{1}{2}$. When $\eta=-\frac{1}{2}$, the 
algebra $Q_2(\eta)$ has a $1$-dimensional radical and the $3$-dimensional
factor algebra over the radical is simple. Note also that, when $F$ is of 
characteristic $3$, we have that $-\frac{1}{2}=1$, and so this situation 
cannot occur.

To summarize, in case (A) we found two primitive algebras: diagram
$A_1$ leads to the familiar $2$-dimensional algebra $2B$ and diagram $A_3$ 
leads to the new $4$-dimensional algebra $Q_2(\eta)$.

\subsection{Type B}

\begin{figure}[ht]
\begin{center}
\includegraphics[scale=0.5]{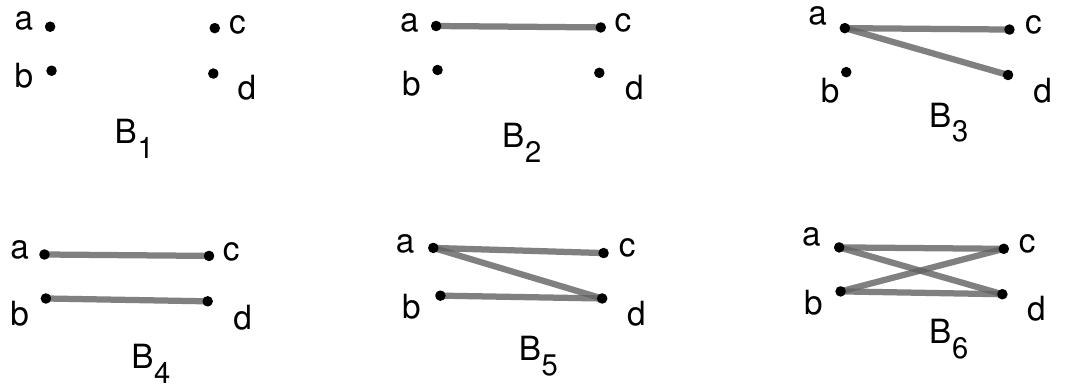}
\caption{Two double axes}\label{Type B}
\end{center}
\end{figure}

Now we deal with case (B), where $x=a+b$, $y=c+d$, and $Z=\{a,b,c,d\}$.
Here $a$ and $b$ are non-collinear in $\Gm$ and also $c$ and $d$ are
non-collinear. Hence the diagram on $Z$ belongs to the list shown in
Figure \ref{Type B}.

Diagrams $B_2$, $B_3$, and $B_5$ do not admit the flip. Hence we discard them.
Again, we do the remaining diagrams one at a time.

\bigskip\noindent
{\bf Diagram $\mathbf{B_1}$:} In this first case, $\hat G$ is an elementary
abelian $2$-group of rank $4$ and $G=\hat G$ or a factor group of $\hat G$
of rank $3$, in which the four generators of $\hat G$ remain distinct.
Furthermore, in all these cases, $xy=(a+b)(c+d)=ac+ad+bc+bd=0+0+0+0=0$, and
so $A$ is the familiar algebra $2B$.

\bigskip\noindent
{\bf Diagram $\mathbf{B_4}$:} Here $\Gm$ has two components, lines $\{a,c,e\}$ 
and $\{b,d,f\}$. Note that $z:=x^{\tau_y}=(a+b)^{\tau_c\tau_d}=e+f=
(c+d)^{\tau_a\tau_b}=y^{\tau_x}$ must be in $A$. Furthermore, $xy=(a+b)(c+d)= 
ac+ad+bc+bd=ac+0+0+bd=\frac{\eta}{2}(a+c-e)+\frac{\eta}{2}(b+d-f)=
\frac{\eta}{2}(x+y-z)$. Applying $\tau_x$ and $\tau_y$, we also obtain 
equalities $xz=\frac{\eta}{2}(x+z-y)$ and $yz=\frac{\eta}{2}(y+z-x)$, showing 
that $A$ is isomorphic to the $3$-dimensional Matsuo algebra $3C(\eta)$. 
Clearly, $x$ and $y$ are both primitive in $A$.

\bigskip\noindent
{\bf Diagram $\mathbf{B_6}$:} Here the Coxeter group $\hat G$ defined by the 
diagram is infinite and so we need an extra relation to identify $G$. This 
extra relation comes from the fact that the order $p$ of $a^db^c$ must also be 
in $\{1,2,3\}$. Let $\hat G(p)$ be the factor group of $\hat G$ over the 
normal subgroup generated by $(a^db^c)^p$. All three groups $\hat G(p)$, 
$p\in\{1,2,3\}$, are finite and their structure is well-known. In fact, they 
are all $3$-transposition groups. We note that $\hat G(1)\cong S_4$ is the 
factor group of both $\hat G(2)$ and $\hat G(3)$. So we just need to describe 
the latter two groups. Let $p\in\{2,3\}$. Consider the $4$-dimensional 
permutational module $V$ of $S_4$ over $\F_p$. Let the vectors $e_i$, 
$i\in\{1,2,3,4\}$, form the basis of $V$ permuted by $S_4$ and $U\subseteq 
V$ be the $3$-dimensional ``sum-zero'' submodule of $V$. Then $\hat G(p)$ is 
isomorphic to the semi-direct product $U:S_4$. Note that, for $p=2$, $U$ 
contains a $1$-dimensional ``all-one'' submodule $D$, which is the center of 
$\hat G(2)$. When $p=3$, $U$ is irreducible. In both cases, $U$ is the unique 
minimal non-central normal subgroup of $\hat G(p)$ and $\hat G(p)/U\cong\hat 
G(1)\cong S_4$. Since $S_4$ does not have proper factor groups containing 
commuting involutions, we conclude that, up to the center (which does not 
influence the Fischer space of $G$), we have that $G=\hat G(p)$ for $p=1$, 
$2$, or $3$.

We consider these possibilities in turn. 

\bigskip
For $p=1$, there is a unique, up to conjugation, identification of the 
generators $a$, $b$, $c$, and $d$ of $G$ with involutions in $S_4$. We set 
$a=(1,2)$, $b=(3,4)$, $c=(1,3)$, and $d=(2,4)$. Then $x=(1,2)+(3,4)$, 
$y=(1,3)+(2,4)$. We also set  $z:=(1,4)+(2,3)$. Then $x$, $y$, and $z$ span 
$A$, which is isomorphic to the 3-dimensional Matsuo algebra $3C(2\eta)$. 
Indeed, $x\cdot y=((1,2)+(3,4))\cdot((1,3)+(2,4))=(1,2)\cdot(1,3)+ (1,2)\cdot(2,4)+(3,4)\cdot(1,3)+(3,4)\cdot(2,4)= 
\frac{\eta}{2}((1,2)+(1,3)-(2,3))+\frac{\eta}{2}((1,2)+(2,4)-(1,4))+ 
\frac{\eta}{2}((3,4)+(1,3)-(1,4))+\frac{\eta}{2}((3,4)+(2,4)-(2,3))= 
\eta(x+y-z)$, and similarly for the other products. Clearly, $x$ and $y$ are 
primitive in $A=3C(2\eta)$.

\bigskip
For $p=2$, the Fischer space $\Gm(2)$ of $\hat G(2)=U:S_4\leq E:S_4$ consists of 
$2\cdot 6=12$ points: $b_{i,j}=(i,j)$ and $c_{i,j}=(e_i+e_j)(i,j)$, for $1\leq 
i<j\leq 4$; and $4\cdot 4=16$ lines, where each ``b'' line 
$\{b_{i,j},b_{i,k},b_{j,k}\}$, $1\leq i<j<k\leq 4$, is complemented by three 
``bc'' lines $\{b_{i,j},c_{i,k},c_{j,k}\}$, $\{b_{i,k},c_{i,j},c_{j,k}\}$, 
and $\{b_{j,k},c_{i,j},c_{i,k}\}$. We can identify our generators $a$, $b$, 
$c$, and $d$ with involutions in $\hat G(2)$ as follows: $a=c_{1,2}=
(e_1+e_2)(1,2)$, $b=b_{3,4}=(3,4)$, $c=b_{1,3}=(1,3)$ and $d=b_{2,4}=(2,4)$. 
Then $a^db^c=e_1+e_4$ is of order $2$. Therefore, this indeed provides the 
required isomorphism from $G$ onto $\hat{G}(2)$.

Recall from the last paragraph of Section \ref{Matsuo} how the Miyamoto 
involutions act on the Fischer space. Since $A$ is invariant under $\tau_x$ 
and $\tau_y$, the following elements lie in $A$:
\begin{equation*}
\begin{aligned}
d_1:=&\,x=c_{1,2}+b_{3,4},\\
d_2:=&\,x^{\tau_y}=((c_{1,2}+b_{3,4})^{\tau_c})^{\tau_d}=
(c_{2,3}+b_{1,4})^{\tau_d}=c_{3,4}+b_{1,2}=b_{1,2}+c_{3,4},\\
d_3:=&\,y=b_{1,3}+b_{2,4},\\
d_4:=&\,y^{\tau_x}=((b_{1,3}+b_{2,4})^{\tau_a})^{\tau_b}=(c_{2,3}+
c_{1,4})^{\tau_b}=c_{2,4}+c_{1,3}=c_{1,3}+c_{2,4}.
\end{aligned}
\end{equation*}
Additionally, $A$ is closed for the algebra product and so it also contains:
\begin{equation*}
\begin{aligned}
w:=&\,-\frac{2}{\eta}xy+2(x+y)\\
=&\,-\frac{2}{\eta}(c_{1,2}+b_{3,4})(b_{1,3}+b_{2,4})+2(c_{1,2}+b_{3,4}+
b_{1,3}+b_{2,4})\\
=&\,-((c_{1,2}+b_{1,3}-c_{2,3})+(c_{1,2}+b_{2,4}-c_{1,4})+(b_{3,4}+b_{1,3}-
b_{1,4})\\
&\,+(b_{3,4}+b_{2,4}-b_{2,3}))+2(c_{1,2}+b_{3,4}+b_{1,3}+b_{2,4})\\
=&\,c_{2,3}+c_{1,4}+b_{1,4}+b_{2,3}=b_{1,4}+c_{1,4}+b_{2,3}+c_{2,3}.
\end{aligned}
\end{equation*}
The five elements above have disjoint support and so they are linearly 
independent. We skip the straightforward calculation of the products and 
simply present them in Table \ref{t:B6(2)}. In particular, $A$ is 
$5$-dimensional with basis $\{d_1,d_2,d_3,d_4,w\}$. It is a new algebra and both 
$x$ and $y$ are primitive in it, since $A$ contains none of the single axes $a$, 
$b$, $c$, and $d$.
\begin{table}[ht]
\begin{center}
\begingroup
\setlength{\tabcolsep}{6pt} 
\scalebox{.60}{
$\begin{tabu}[h!]{|c||c|c|c|c|c|}
\hline
&d_1&d_2&d_3&d_4&w\\
\hline\hline
d_1&d_1&0&\frac{\eta}{2}(2d_1+2d_3-w)&\frac{\eta}{2}(2d_1+2d_4-w)& 
\eta(2d_1-d_3-d_4+w)\\
\hline
d_2&0&d_2&\frac{\eta}{2}(2d_2+2d_3-w)&\frac{\eta}{2}(2d_2+2d_4-w)&
\eta(2d_2-d_3-d_4+w)\\
\hline
d_3&\frac{\eta}{2}(2d_1+2d_3-w)&\frac{\eta}{2}(2d_2+2d_3-w)&d_3&0& 
\eta(2d_3-d_1-d_2+w)\\
\hline
d_4&\frac{\eta}{2}(2d_1+2d_4-w)&\frac{\eta}{2}(2d_2+2d_4-w)&0&d_4&
\eta(2d_4-d_1-d_2+w)\\
\hline                                                         w&\eta(2d_1-d_3-d_4+w)&\eta(2d_2-d_3-d_4+w)&\eta(2d_3-d_1-d_2+w)&
\eta(2d_4-d_1-d_2+w)&w\\
\hline
\end{tabu}$}
\caption{The $5$-dimensional algebra}\label{t:B6(2)}
\endgroup
\end{center}
\end{table}

We check this algebra for simplicity using GAP \cite{GAP} and the structure 
theory from Subsection \ref{ideals}. We will provide details for this algebra 
and skip them for the later examples. The Gram matrix for the Frobenius form 
with respect to the basis $\{d_1,d_2,d_3,d_4,w\}$ can easily be found and it is 
as follows:
$$
\left(
\begin{array}{rrrrr}
2&0&2\eta&2\eta&4\eta\\
0&2&2\eta&2\eta&4\eta\\
2\eta&2\eta&2&0&4\eta\\
2\eta&2\eta&0&2&4\eta\\
4\eta&4\eta&4\eta&4\eta&4
\end{array}
\right).
$$

Manifestly, the projection graph on $\{x,y\}=\{d_1,d_3\}$ is connected. This 
means that the algebra has no proper ideals containing axes. The determinant of 
the above Gram matrix is $1024\eta^3-768\eta^2+64$ and it has roots $\frac{1}{2}$ 
(twice) and $-\frac{1}{4}$. As $\eta\neq\frac{1}{2}$ by assumption, the only 
special value is $\eta=-\frac{1}{4}$, for which the radical of the Frobenius form 
is $1$-dimensional. The factor over the radical is a simple $4$-dimensional 
algebra. 

Note that if the characteristic of $\F$ is $3$ then $-\frac{1}{4}=\frac{1}{2}$ 
and if the characteristic is $5$ then $-\frac{1}{4}=1$. So for these 
characteristics the $5$-dimensional algebra is always simple.

\bigskip
Let, finally, $p=3$. We start again by describing the Fischer space $\Gm(3)$ of 
$\hat G(3)=U:S_4\leq E:S_4$, where $E$ is now the permutational module over $\F_3$. 
The Fischer space has $3\cdot 6=18$ points. Namely, for each point $b_{i,j}=(i,j)= 
b_{j,i}$, contained in the complement $S_4$, there are two further points $c_{i,j}= 
(e_i-e_j)(i,j)$ and $c_{j,i}=(e_j-e_i)(i,j)$. The lines in $\Gm(3)$ are of several 
types. First, for each $1\leq i<j\leq 4$, the triple (1) $\{b_{i,j},c_{i,j},c_{j,i}\}$ 
is a line. Secondly, for distinct $i$, $j$, and $k$ in $\{1,2,3,4\}$, the triples 
(2) $\{b_{i,j},b_{i,k},b_{j,k}\}$, (3) $\{b_{i,j},c_{i,k},c_{j,k}\}$, (4) 
$\{b_{j,k},c_{i,j},c_{i,k}\}$, and (5) $\{c_{i,j},c_{j,k},c_{k,i}\}$ are lines. This 
gives the total of $42$ lines of $\Gm(3)$, including six lines of type (1), four lines 
of type (2), twelve lines of type (3), twelve lines of type (4), and eight lines of 
type (5). This information is used to multiply in $M$ and to act by the Miyamoto 
involutions, as described in Section \ref{Matsuo}.   

We identify $G$ with $\hat G(3)$ by taking $a=c_{1,2}$, $b=b_{3,4}$, $c=b_{1,3}$, and 
$d=c_{2,4}$. Then $a^db^c=c_{1,2}^{\tau_d}b_{3,4}^{\tau_c}=c_{4,1}b_{1,4}=
(e_4-e_1)(1,4)(1,4)=e_4-e_1$ is of order $3$, so the required relation hold. Since $A$ 
is invariant under $\tau_x$ and $\tau_y$, the following elements lie in $A$:
\begin{equation*}     
\begin{aligned}
d_1:=&\,x=c_{1,2}+b_{3,4},\\      
d_2:=&\,x^{\tau_y}=((c_{1,2}+b_{3,4})^{\tau_c})^{\tau_d}=(c_{3,2}+b_{1,4})^{\tau_d}= 
c_{4,3}+c_{2,1}=c_{2,1}+c_{4,3},\\ 
d_3:=&\,d_2^{\tau_x}=((c_{2,1}+c_{4,3})^{\tau_a})^{\tau_b}=(b_{1,2}+c_{4,3})^{\tau_b}= 
b_{1,2}+c_{3,4},\\ 
d_4:=&\,y=b_{1,3}+c_{2,4},\\ 
d_5:=&\,y^{\tau_x}=((b_{1,3}+c_{2,4})^{\tau_a})^{\tau_b}=(c_{3,2}+c_{4,1})^{\tau_b}=
c_{4,2}+c_{3,1}=c_{3,1}+c_{4,2},\\ 
d_6:=&\,d_5^{\tau_y}=((c_{3,1}+c_{4,2})^{\tau_c})^{\tau_d}=(c_{1,3}+c_{4,2})^{\tau_d}=
c_{1,3}+b_{2,4},\\ 
u:=&\,-\frac{1}{\eta}xd_6+x+d_6\\
=&\,-\frac{1}{\eta}(c_{1,2}+b_{3,4})(c_{1,3}+b_{2,4})+c_{1,2}+b_{3,4}+c_{1,3}+b_{2,4}\\
=&\,-\frac{1}{2}((c_{1,2}+c_{1,3}-b_{2,3})+(c_{1,2}+b_{2,4}-c_{1,4})+
(b_{3,4}+c_{1,3}-c_{1,4})\\
&\,+(b_{3,4}+b_{2,4}-b_{2,3}))+c_{1,2}+b_{3,4}+c_{1,3}+b_{2,4}\\
=&\,c_{1,4}+b_{2,3},\\ 
w:=&\,-\frac{2}{\eta}xy+2(x+y)\\
=&\,-\frac{2}{\eta}(c_{1,2}+b_{3,4})(b_{1,3}+c_{2,4})+
2(c_{1,2}+b_{3,4}+b_{1,3}+c_{2,4})\\
=&\,-((c_{1,2}+b_{1,3}-c_{3,2})+(c_{1,2}+c_{2,4}-c_{4,1})+(b_{3,4}+b_{1,3}-b_{1,4})\\
&\,+(b_{3,4}+c_{2,4}-c_{2,3}))+2(c_{1,2}+b_{3,4}+b_{1,3}+c_{2,4})\\
=&\,c_{3,2}+c_{4,1}+b_{1,4}+c_{2,3}=b_{1,4}+c_{4,1}+c_{2,3}+c_{3,2}. 
\end{aligned}
\end{equation*}
These eight elements have disjoint support hence they are linearly independent.  
Calculations show that the pairwise products of these elements are as shown in Table 
\ref{t:B6(3)}. So this is a basis of the (new) subalgebra $A$. Since $a$, $b$, $c$, and 
$d$ are not in $A$, both $x$ and $y$ are primitive in $A$.
\begin{table}[ht]
\begin{center}
\begingroup
\setlength{\tabcolsep}{20pt}
\scalebox{.5}{
$\begin{tabu}[ht]{|c||c|c|c|c|c|c|c|c|}
\hline
&d_1&d_2&d_3&d_4&d_5&d_6&u&w\\
\hline
\hline
d_1&d_1&\frac{\eta}{2}(d_1+d_2-d_3)&\frac{\eta}{2}(d_1+d_3-d_2)
&\frac{\eta}{2}(2d_1+2d_4-w)&\frac{\eta}{2}(2d_1+2d_5-w)&\eta(d_1+d_6-u)
&\eta(d_1-d_6+u)&\begin{tabu}{@{}c@{}}\eta(2d_1+w)\\-\eta(d_4+d_5)\end{tabu}\\
\hline
d_2&\frac{\eta}{2}(d_1+d_2-d_3)&d_2&\frac{\eta}{2}(d_2+d_3-d_1)
&\frac{\eta}{2}(2d_2+2d_4-w)&\eta(d_2+d_5-u)&\frac{\eta}{2}(2d_2+2d_6-w)
&\eta(d_2-d_5+u)&\begin{tabu}{@{}c@{}}\eta(2d_2+w)\\-\eta(d_4+d_6)\end{tabu}\\
\hline
d_3&\frac{\eta}{2}(d_1+d_3-d_2)&\frac{\eta}{2}(d_2+d_3-d_1)&d_3&\eta(d_3+d_4-u)
&\frac{\eta}{2}(2d_3+2d_5-w)&\frac{\eta}{2}(2d_3+2d_6-w)&\eta(d_3-d_4+u)
&\begin{tabu}{@{}c@{}}\eta(2d_3+w)\\-\eta(d_6+d_5)\end{tabu}\\
\hline
d_4&\frac{\eta}{2}(2d_1+2d_4-w)&\frac{\eta}{2}(2d_2+2d_4-w)&\eta(d_4+d_3-u)&d_4
&\frac{\eta}{2}(d_4+d_5-d_6)&\frac{\eta}{2}(d_4+d_6-d_5)&\eta(d_4-d_3+u)
&\begin{tabu}{@{}c@{}}\eta(2d_4+w)\\-\eta(d_1+d_2)\end{tabu}\\
\hline
d_5&\frac{\eta}{2}(2d_1+2d_5-w)&\eta(d_2+d_5-u)&\frac{\eta}{2}(2d_3+2d_5-w)
&\frac{\eta}{2}(d_4+d_5-d_6)&d_5&\frac{\eta}{2}(d_5+d_6-d_4)&\eta(d_5-d_2+u)
&\begin{tabu}{@{}c@{}}\eta(2d_5+w)\\-\eta(d_1+d_3)\end{tabu}\\
\hline
d_6&\eta(d_1+d_6-u)&\frac{\eta}{2}(2d_2+2d_6-w)&\frac{\eta}{2}(2d_6+2d_3-w)
&\frac{\eta}{2}(d_4+d_6-d_5)&\frac{\eta}{2}(d_5+d_6-d_4)&d_6
&\eta(d_6-d_1+u)&\begin{tabu}{@{}c@{}}\eta(2d_6+w)\\-\eta(d_2+d_3)\end{tabu}\\
\hline
u&\eta(d_1-d_6+u)&\eta(d_2-d_5+u)&\eta(d_3-d_4+u)&\eta(d_4-d_3+u)&\eta(d_5-d_2+u)
&\eta(d_6-d_1+u)&u&\eta u\\
\hline
w&\begin{tabu}{@{}c@{}}\eta(2d_1+w)\\-\eta(d_4+d_5)\end{tabu}
&\begin{tabu}{@{}c@{}}\eta(2d_2+w)\\-\eta(d_4+d_6)\end{tabu}
&\begin{tabu}{@{}c@{}}\eta(2d_3+w)\\-\eta(d_6+d_5)\end{tabu}
&\begin{tabu}{@{}c@{}}\eta(2d_4+w)\\-\eta(d_1+d_2)\end{tabu}
&\begin{tabu}{@{}c@{}}\eta(2d_5+w)\\-\eta(d_1+d_3)\end{tabu}
&\begin{tabu}{@{}c@{}}\eta(2d_6+w)\\-\eta(d_2+d_3)\end{tabu}
&\eta u&(\eta+1)w-\eta u\\
\hline
\end{tabu}$}
\caption{The $8$-dimensional algebra}\label{t:B6(3)}
\endgroup
\end{center}
\end{table}

We build in GAP the Gram matrix of the Frobenius form and find its determinant, 
which turns out to be $-448\eta^8+4864\eta^7-21808\eta^6+51664\eta^5-
68320\eta^4+48256\eta^3-14848\eta^2+256\eta+512$. This has roots $-\frac{1}{7}$, 
$\frac{1}{2}$ of multiplicity $2$, and $2$ of multiplicity $5$. Consequently, 
the algebra is simple unless $\eta=-\frac{1}{7}$ or $\eta=2$. In the former 
case, the radical is $1$-dimensional and the latter case it is $5$-dimensional.

If the characteristic of $\F$ is $3$ then both $-\frac{1}{7}$ and $2$ are equal 
to $\frac{1}{2}$, so the algebra is always simple, as $\eta\neq\frac{1}{2}$. 
Another special characteristic for this algebra is $5$, because then 
$-\frac{1}{7}=2$.

\subsection{Discussion}

The list of $2$-generated algebras of Monster type $(2\eta,\eta)$ 
found in this section consists of three versions of the algebra $2B=\F^2$ 
(generated by two single axes, one single and one double axes, and two double 
axes), two versions of $3C(\eta)$ (generated by two single axes and by two 
double axes), algebra $3C(2\eta)$ (two double axes only), the new 
$4$-dimensional algebra $Q_2(\eta)$ (single and double axes), the new 
$5$-dimensional algebra (two double axes), and the new $8$-dimensional algebra 
(two double axes). We put forward the following question.

\begin{question}
Is it true that the above is the complete list of $2$-generated primitive 
algebras of Monster type $(2\eta,\eta)$?
\end{question}

It is interesting to compare this with the list of eight Norton-Sakuma 
algebras, which is the complete list of $2$-generated primitive algebras of 
Monster type $(\qu,\thi)$. The most glaring difference is the absence of an 
equivalent of the $6$-dimensional algebra $5A$, where $|\tau_x\tau_y|=5$. 
Hence we also put forward the following partial case of the above question.

\begin{question}
Is it true that there is no $2$-generated primitive algebra $\dla x,y\dra$ of 
Monster type $(2\eta,\eta)$ satisfying $|\tau_x\tau_y|=5$?
\end{question}

In our list, we operate in terms of single and double axis, but these are only 
fully defined in the context of an enveloping Matsuo algebra $M$.

\begin{question}
In an arbitrary primitive algebra of Monster type $(2\eta,\eta)$, can the axes 
be always split into single and double?
\end{question}

In the presence of a Frobenius form, a possible approach to this question is 
via the length of the axes. Namely, it might be possible to scale the form so 
that all ``single'' axes $x$ satisfy $(x,x)=1$ and all ``double'' axes $x$ satisfy 
$(x,x)=2$. 

Another idea is to base the distinction between ``single'' and ``double'' on 
the absence or presence of the $2\eta$-eigenspace.

\section{Some 3-generated subalgebras}

\label{3-generated}

The natural next step is to try and classify $3$-generated primitive subalgebras. 
Again, we can assume that at least one of the three generators is a 
double axis. This gives us three further cases: (C) two single axes and one 
double axis; (D) one single axis and two double axes; and (E) three double axes.

It is easy to list the relevant diagrams on the support $Z$ of the set 
$Y=\{x,y,z\}$ of three generating axes. Note that here we only need to list the 
diagrams admitting the flip as only they can lead to primitive algebras.

\begin{figure}[ht]
\begin{center}
\includegraphics[scale=0.5]{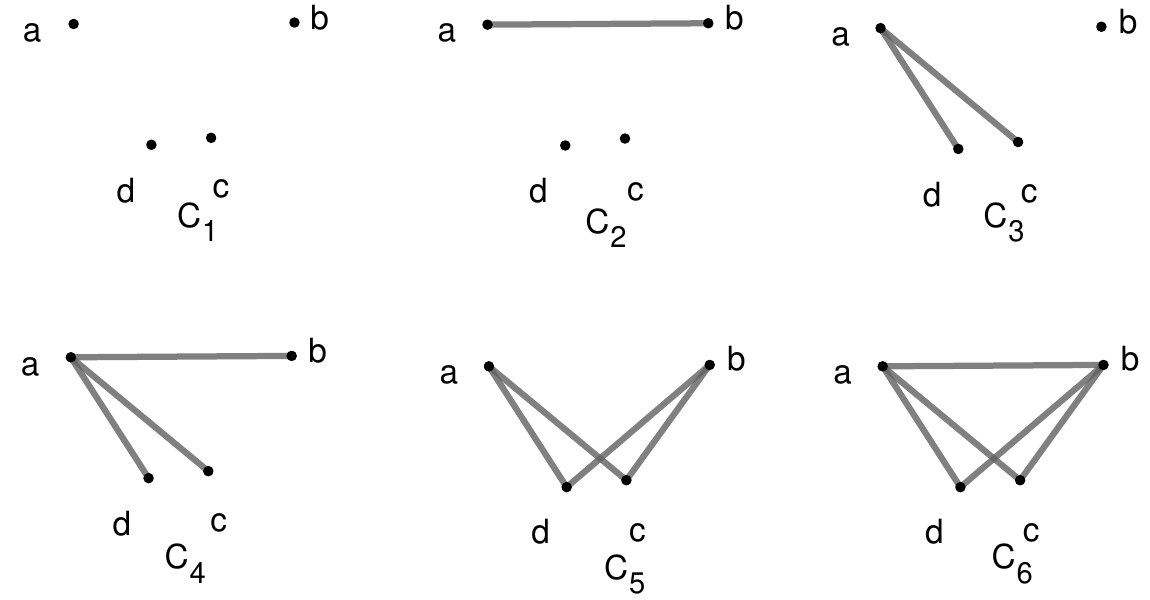}
\caption{Two single axes and one double axis}\label{Type C}
\end{center}
\end{figure}

All possible symmetric diagrams for case (C) are shown in Figure \ref{Type C}.
In this section we will completely enumerate the primitive algebras of 
Monster type $(2\eta,\eta)$ arising in this case. We also discuss the 
properties of these algebras.

\begin{figure}[ht]
\begin{center}
\includegraphics[scale=0.4]{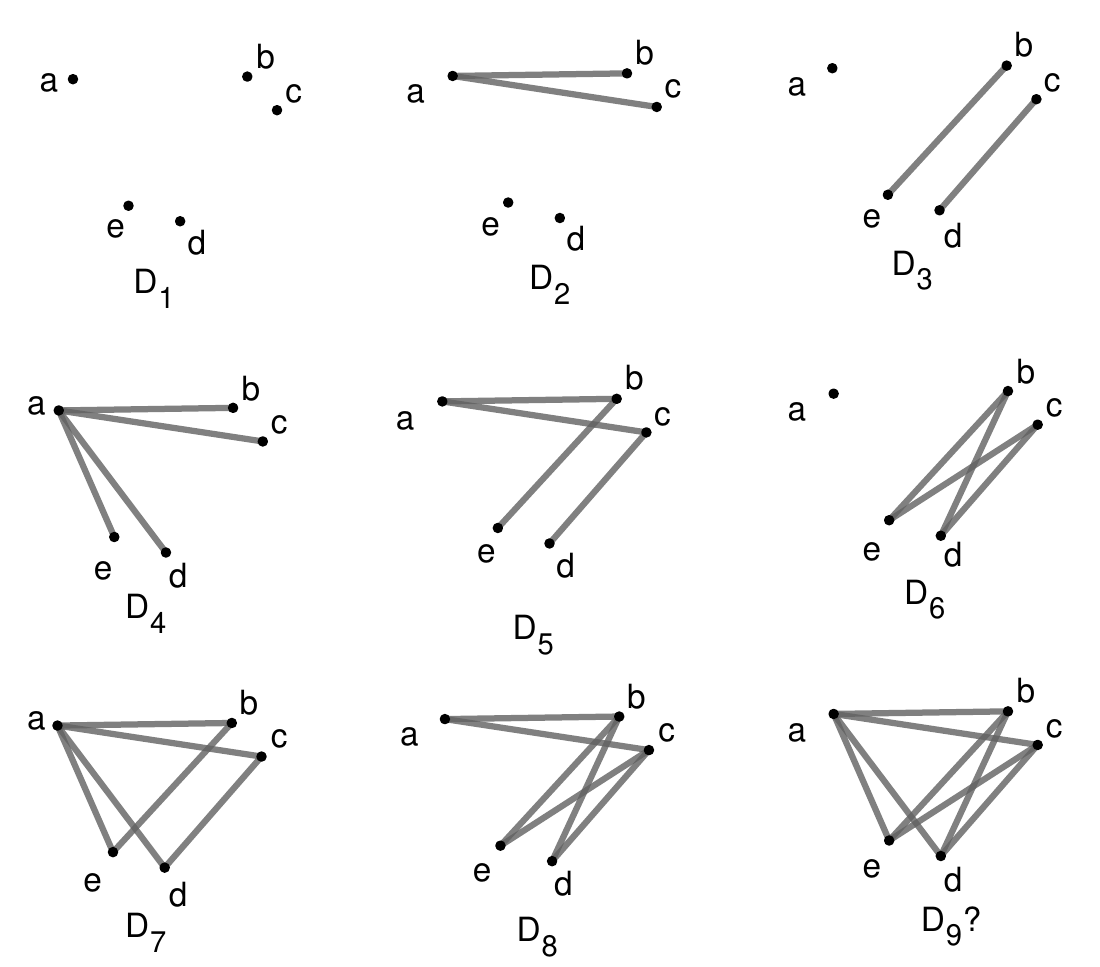}
\caption{One single axis and two double axes}\label{Type D}
\end{center}
\end{figure}

\begin{figure}[ht]
\begin{center}
\includegraphics[scale=0.4]{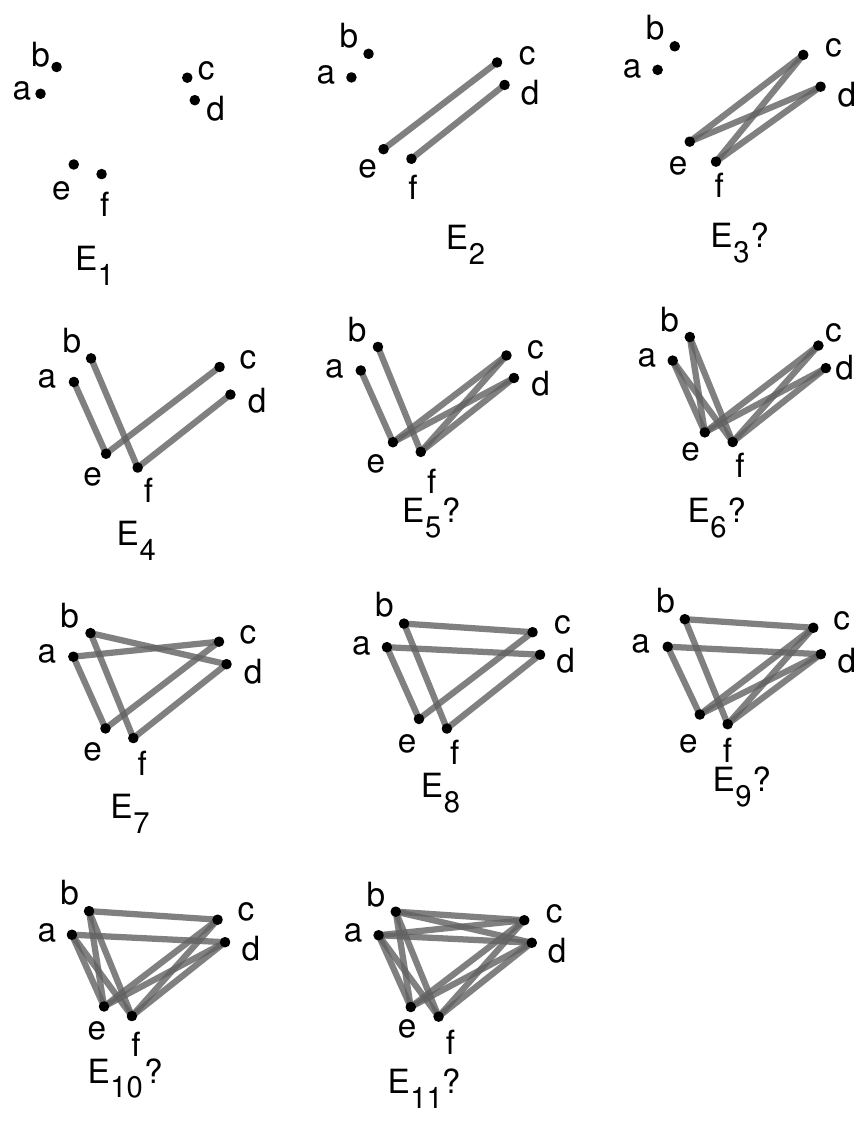}
\caption{Three double axes}\label{Type E}
\end{center}
\end{figure}

The full list of symmetric diagrams arising in cases (D) and (E) are 
shown in Figures \ref{Type D} and \ref{Type E}, respectively. 
These cases constitute an ongoing project and we are hoping to report on the 
results in a separate paper.

\subsection{Type C}

In this case,  $x=a$ and $y=b$ are single axes while $z=c+d$ is a double axis. 
Therefore, $Z=\supp(\{x,y,z\})=\{a,b,c,d\}$.

\bigskip\noindent
{\bf Diagram $\mathbf{C_1}$:} In this first case, $\hat G$ is an elementary 
abelian group of order $2^4$ and $G=\hat G$ or a factor group of $\hat G$
of order $2^3$, in which the images of the four generators of $\hat G$ remain 
distinct. Futhermore, $xy=0$, $xz=0$, $yz=0$, which means that $A$ is 
isomorphic to $\F\oplus\F\oplus\F=\F^3$.

\bigskip\noindent
{\bf Diagram $\mathbf{C_2}$:} Here $\hat G\cong\la a,b\ra\times\la c\ra\times 
\la d\ra\cong S_3\times C_2\times C_2$ and it is easy to see that $G=\hat G$, 
as it cannot be any of the proper factor groups. Since $xz=0$ and $yz=0$, the 
algebra $A=\dla x,y,z\dra$ decomposes as the direct sum $\dla x,y\dra\oplus 
\dla z\dra\cong 3C(\eta)\oplus\F$. 

\bigskip\noindent
{\bf Diagram $\mathbf{C_3}$:} In this case, $\hat G=\la a,c,d\ra\times\la 
b\ra\cong S_4\times C_2$. Again, $G$ can only be the full $\hat G$. The 
algebra $A$ decomposes, in turn, as $\dla x,z\dra\oplus\dla y\dra\cong 
Q_2(\eta)\oplus\F$, where $Q_2(\eta)$ is the algebra from case $A_3$.

\bigskip\noindent
{\bf Diagram $\mathbf{C_4}$:} Here the diagram on $Z$ coincides with 
the Coxeter diagram $D_4$ and hence $\hat G$ is isomorphic to the Weyl group 
$W(D_4)\cong 2^3:S_4$ of order 192.  We claim that $G=\hat G$ or 
$\hat G/Z(\hat G)$. Indeed, $G_b:=\la a,c,d\ra\leq G$ is a factor group of 
$W(A_3)\cong S_4$ containing distinct commuting transpositions. Hence  
$G_b\cong S_4$ and, similarly, $G_c:=\la a,b,d\ra\cong S_4$ and 
$G_d:=\la a,b,c\ra\cong S_4$. 

Now we claim that $G_b\cap G_c=\la a,d\ra$. If not, then $G_b=G_c$, as $\la 
a,d\ra\cong S_3$ is maximal in both groups. Hence $b\in G_b$. We can identify 
$a$, $c$, and $d$ with transpositions $(1,2)$, $(1,3)$, and $(2,4)$, 
respectively. Then $b$ must be a transposition commuting with $c$ and $d$, 
so $b$ must be $c$ or $d$; a contradiction. 

Thus, $|G|\geq 24+24-6=42$ and hence $G$ is a factor group of $\hat G$ over a 
normal subgroup of order at most $\frac{192}{42}<5$. However, the only normal 
subgroup of $\hat G$ of order at most $4$ is the center $Z(\hat G)$, and so, 
indeed, $G=\hat G$ or $\hat G/Z(\hat G)$. Note that these two groups have the 
same Fischer space, so we can select whichever one we prefer. Thus, we may 
assume that $G=\hat G$.

Note that the present $G=\hat G$ appeared in Section \ref{2-generated}, in the 
case of diagram $B_6$, as the group $\hat G(2)$. In particular, $\hat G$ 
is isomorphic to the index $2$ subgroup $U:S_4$ in the group $E:S_4$. Here 
$E$ is the permutational module for $S_4$, with the basis 
$\{e_1,e_2,e_3,e_4\}$ permuted by the complement $S_4$, and $U$ is the 
``sum-zero'' submodule of $E$. 

Recall the Fischer space of $\hat G(2)$ consisting of $2\cdot 
6=12$ points: $b_{i,j}=(i,j)$ and $c_{i,j}=(e_i+e_j)(i,j)$, for $1\leq i<j\leq 
4$; and $4\cdot 4=16$ lines $\{b_{i,j},b_{i,k},b_{j,k}\}$, 
$\{b_{i,j},c_{i,k},c_{j,k}\}$, $\{b_{i,k},c_{i,j},c_{j,k}\}$, 
and $\{b_{j,k},c_{i,j},c_{i,k}\}$, for $1\leq i<j<k\leq 4$.

We have already identified $a$, $c$, and $d$ with $b_{1,2}=(1,2)$, 
$b_{1,3}=(1,3)$, and $b_{2,4}=(2,4)$, respectively. We can, for example, 
identify $b$ with $c_{1,3}=(e_1+e_3)(1,3)$, and then all relations are 
satisfied. 

Recall that $A=\dla x,y,z\dra$ is invariant under $\tau_x=\tau_a$, 
$\tau_y=\tau_b$, and $\tau_z=\tau_c\tau_d$. Recall also that the Miyamoto 
involution $\tau_f$, for a point $f$, fixes $f$ and every point non-collinear 
with $f$ and it switches the two points other than $f$ on each line through 
$f$. 

Applying this, the algebra $A=\dla x,y,z\dra$ contains the following axes:
\begin{center}
$\begin{aligned}
s_1:=&\,x=a=b_{1,2},\\
s_2:=&\,y=b=c_{1,3},\\
s_3:=&\,x^{\tau_y}=y^{\tau_x}=c_{2,3},\\
s_4:=&\,s_3^{\tau_z}=(c_{2,3}^{\tau_c})^{\tau_d}=c_{1,2}^{\tau_d}=c_{1,4},\\
s_5:=&\,s_4^{\tau_x}=c_{1,4}^{\tau_x}=c_{2,4},\\
s_6:=&\,x^{\tau_z}=(x^{\tau_c})^{\tau_d}=b_{2,3}^{\tau_d}=b_{3,4},\\
d_2:=&\,z=b_{1,3}+b_{2,4},\\
d_3:=&\,z^{\tau_x}=b_{2,3}+b_{1,4}=b_{1,4}+b_{2,3},\\ 
d_1:=&\,d_3^{\tau_y}=c_{3,4}+c_{1,2}=c_{1,2}+c_{3,4}.\\
\end{aligned}$
\end{center}
Here, as in Section~\ref{2-generated}, $s_1,\ldots,s_6$ are single axes and 
$d_1,d_2,d_3$ are double axes. (The order of the double axes is altered so 
that the multiplication table looks more symmetric.) We claim that 
the nine axes form a basis of $A$. The supports of the axes are disjoint and 
so they are certainly linearly independent. It remains to see that the 
subspace they span is closed for the algebra product. It is a straightforward 
calculation using the description of lines above and utilizing the substantial 
symmetry we have. Note that the six single axes span a subalgebra isomorphic 
to the Matsuo algebra of $S_4$, and so these products are immediate. Here is 
a sample calculation involving double axes: $d_1d_3=
(c_{1,2}+c_{3,4})(b_{1,4}+b_{2,3})=
c_{1,2}b_{1,4}+c_{1,2}b_{2,3}+c_{3,4}b_{1,4}+c_{3,4}b_{2,3}=
\frac{\eta}{2}(c_{1,2}+b_{1,4}-c_{2,4})+\frac{\eta}{2}(c_{1,2}+b_{2,3}-
c_{1,3})+\frac{\eta}{2}(c_{3,4}+b_{1,4}-c_{1,3})+\frac{\eta}{2}(c_{3,4}+
b_{2,3}-c_{2,4})=\eta(d_1+d_3-s_2-s_5)$. By the way, this calculation indicates, 
using the symmetry, that the three double axes span a subalgebra $3C(2\eta)$. 

The complete table of products is in Table~\ref{t:algebra C4}. 
\begin{table}[ht]
\begin{center}
\begingroup
\setlength{\tabcolsep}{20pt}
\scalebox{.42}{
$\begin{tabu}[h!]{|c||c|c|c|c|c|c|c|c|c|}
\hline
&s_1&s_2&s_3&s_4&s_5&s_6&d_1&d_2&d_3\\
\hline\hline
s_1&s_1&\frac{\eta}{2}(s_1+s_2-s_3)&\frac{\eta}{2}(s_1+s_3-s_2)& 
\frac{\eta}{2}(s_1+s_4-s_5)&\frac{\eta}{2}(s_1+s_5-s_4)&0&
0&\frac{\eta}{2}(2s_1+d_2-d_3)&\frac{\eta}{2}(2s_1+d_3-d_2)\\
\hline 
s_2&\frac{\eta}{2}(s_1+s_2-s_3)&s_2&\frac{\eta}{2}(s_2+s_3-s_1)&
\frac{\eta}{2}(s_2+s_4-s_6)&0&\frac{\eta}{2}(s_2+s_6-s_4)&
\frac{\eta}{2}(2s_2+d_1-d_3)&0&\frac{\eta}{2}(2s_2+d_3-d_1)\\
\hline
s_3&\frac{\eta}{2}(s_1+s_3-s_2)&\frac{\eta}{2}(s_2+s_3-s_1)&s_3&0&
\frac{\eta}{2}(s_3+s_5-s_6)&\frac{\eta}{2}(s_3+s_6-s_5)&
\frac{\eta}{2}(2s_3+d_1-d_2)&\frac{\eta}{2}(2s_3+d_2-d_1)&0\\
\hline 
s_4&\frac{\eta}{2}(s_1+s_4-s_5)&\frac{\eta}{2}(s_2+s_4-s_6)&0&s_4&
\frac{\eta}{2}(s_4+s_5-s_1)&\frac{\eta}{2}(s_4+s_6-s_2)&
\frac{\eta}{2}(2s_4+d_1-d_2)&\frac{\eta}{2}(2s_4+d_2-d_1)&0\\
\hline 
s_5&\frac{\eta}{2}(s_1+s_5-s_4)&0&\frac{\eta}{2}(s_3+s_5-s_6)&
\frac{\eta}{2}(s_4+s_5-s_1)&s_5&\frac{\eta}{2}(s_5+s_6-s_3)&
\frac{\eta}{2}(2s_5+d_1-d_3)&0&\frac{\eta}{2}(2s_5+d_3-d_1)\\ 
\hline 
s_6&0&\frac{\eta}{2}(s_2+s_6-s_4)&\frac{\eta}{2}(s_3+s_6-s_5)&
\frac{\eta}{2}(s_4+s_6-s_2)&\frac{\eta}{2}(s_5+s_6-s_3)&s_6& 
0&\frac{\eta}{2}(2s_6+d_2-d_3)&\frac{\eta}{2}(2s_6+d_3-d_2)\\
\hline 
d_1&0&\frac{\eta}{2}(2s_2+d_1-d_3)&\frac{\eta}{2}(2s_3+d_1-d_2)&
\frac{\eta}{2}(2s_4+d_1-d_2)&\frac{\eta}{2}(2s_5+d_1-d_3)&0&
d_1&\eta(d_1+d_2-s_3-s_4)&\eta(d_1+d_3-s_2-s_5)\\
\hline
d_2&\frac{\eta}{2}(2s_1+d_2-d_3)&0&\frac{\eta}{2}(2s_3+d_2-d_1)&
\frac{\eta}{2}(2s_4+d_2-d_1)&0&\frac{\eta}{2}(2s_6+d_2-d_3)&
\eta(d_1+d_2-s_3-s_4)&d_2&\eta(d_2+d_3-s_1-s_6)\\
\hline 
d_3&\frac{\eta}{2}(2s_1+d_3-d_2)&\frac{\eta}{2}(2s_2+d_3-d_1)&0&
0&\frac{\eta}{2}(2s_5+d_3-d_1)&\frac{\eta}{2}(2s_6+d_3-d_2)&
\eta(d_1+d_3-s_2-s_5)&\eta(d_2+d_3-s_1-s_6)&d_3\\
\hline 
\end{tabu}$}
\caption{The 9-dimensional algebra $Q^3(\eta)$}\label{t:algebra C4}
\endgroup
\end{center}
\end{table}
As a consequence, we infer that this $9$-dimensional subalgebra\footnote{
This algebra has been since included by Alsaeedi \cite{al} in an infinite 
series $Q^k(\eta)$ that is dual, in a sense, to our series $Q_k(\eta)$ 
from Section \ref{symmetric group}. Also in this new series is the algebra
$Q_2(\eta)=Q^2(\eta)$ from the case of diagram $A_3$.}
is primitive, since otherwise $c$ and $d$ would belong to the subalgebra.

The Gram matrix of the Frobenius form was computed in GAP and its determinant 
is $128\eta^3-96\eta^2+8$, having roots $\frac{1}{2}$ (with multiplicity two) 
and $-\frac{1}{4}$. Hence the algebra is not simple only when 
$\eta=-\frac{1}{4}$. In characteristics $3$ and $5$, $-\frac{1}{4}=\frac{1}{2}$ 
and $1$, respectively, so for such fields $\F$, the algebra is always simple.

\bigskip\noindent
{\bf Diagram $\mathbf{C_5}$:} Here the diagram on $Z$ is the same as in case 
$B_6$ in Section \ref{2-generated}. Correspondingly, we again need to consider 
the $3$-transposition factor groups $\hat G(p)$, $p=1$, $2$, or $3$. 

\bigskip
First let $p=1$. As $G=\hat G(1)\cong S_4$, we can identify $a$, $c$, and $d$ 
with $(1,2)$, $(1,3)$, and $(2,4)$, respectively. Since $b$ is a transposition 
commuting with $a$, we must have that $b=(3,4)$. Now, the subalgebra $\dla 
x,z\dra=\dla a,c+d\dra$ falls into the case $A_3$ of Section 
\ref{2-generated}. Hence this subalgebra is none other than the algebra 
$Q_2(\eta)$ and it also contains $y=b=x^{\tau_z}$. Thus, in this 
case, $\dla x,y,z\dra\cong Q_2(\eta)$.  

\bigskip
Suppose now that $p=2$. The group $G=\hat G=U:S_4\leq E:S_4$ is the same 
group that we met for the diagrams $B_6$ and $C_4$. In particular, we can 
reuse the description of the Fischer space of $G$. We  can identify $a$, $b$, 
$c$, and $d$ with $b_{1,2}$, $c_{3,4}$, $b_{1,3}$, and $b_{2,4}$, 
respectively. Applying Miyamoto involutions, we obtain the following axes in 
$A=\dla x,y,z\dra$:
\begin{center}
$\begin{aligned}
s_1:=&\,x=a=b_{1,2},\\
s_2:=&\,y^{\tau_z}=(c_{3,4}^{\tau_c})^{\tau_d}=c_{1,4}^{\tau_d}=c_{1,2},\\
s_3:=&\,x^{\tau_z}=(b_{1,2}^{\tau_c})^{\tau_d}=b_{2,3}^{\tau_d}=b_{3,4},\\
s_4:=&\,y=b=c_{3,4},\\
d_1:=&\,z=c+d=b_{1,3}+b_{2,4},\\
d_3:=&\,z^{\tau_x}=(b_{1,3}+b_{2,4})^{\tau_a}=b_{2,3}+b_{1,4}=b_{1,4}+b_{2,3},\\
d_4:=&\,z^{\tau_y}=(b_{1,3}+b_{2,4})^{\tau_b}=c_{1,4}+c_{2,3},\\
d_2:=&\,d_4^{\tau_x}=(c_{1,4}+c_{2,3})^{\tau_a}=c_{2,4}+c_{1,3}=c_{1,3}+c_{2,4}.\\
\end{aligned}$
\end{center}
(We selected the order of the double axes that exhibits the symmetry of the 
multiplication table.)

Again, it is immediate to see that the supports of these axes partition the 
Fischer space. In particular, the axes are linearly independent. To see that 
the $8$-dimensional subspace they span is the subalgebra $A=\dla x,y,z\dra$, 
we compute the products. The results of the computation in 
Table~\ref{t:algebra C5} show that $A$ is indeed $8$-dimensional spanned by 
the above eight axes. 
\begin{table}
\begin{center}
\begingroup
\setlength{\tabcolsep}{20pt}
\scalebox{.44}{
$\begin{tabu}[h!]{|c||c|c|c|c|c|c|c|c|}
\hline
&s_1&s_2&s_3&s_4&d_1&d_2&d_3&d_4\\
\hline\hline
s_1&s_1&0&0&0&\frac{\eta}{2}(2s_1+d_1-d_3)&\frac{\eta}{2}(2s_1+d_2-d_4)& 
\frac{\eta}{2}(2s_1+d_3-d_1)&\frac{\eta}{2}(2s_1+d_4-d_2)\\
\hline 
s_2&0&s_2&0&0&\frac{\eta}{2}(2s_2+d_1-d_4)&\frac{\eta}{2}(2s_2+d_2-d_3)& 
\frac{\eta}{2}(2s_2+d_3-d_2)&\frac{\eta}{2}(2s_2+d_4-d_1)\\
\hline
s_3&0&0&s_3&0&\frac{\eta}{2}(2s_3+d_1-d_3)&\frac{\eta}{2}(2s_3+d_2-d_4)& 
\frac{\eta}{2}(2s_3+d_3-d_1)&\frac{\eta}{2}(2s_3+d_4-d_2)\\
\hline
s_4&0&0&0&s_4&\frac{\eta}{2}(2s_4+d_1-d_4)&\frac{\eta}{2}(2s_4+d_2-d_3)&
\frac{\eta}{2}(2s_4+d_3-d_2)&\frac{\eta}{2}(2s_4+d_4-d_1)\\
\hline
d_1&\frac{\eta}{2}(2s_1+d_1-d_3)&\frac{\eta}{2}(2s_2+d_1-d_4)& 
\frac{\eta}{2}(2s_3+d_1-d_3)&\frac{\eta}{2}(2s_4+d_1-d_4)&d_1&0& 
\eta(-s_1-s_3+d_1+d_3)&\eta(-s_2-s_4+d_1+d_4)\\
\hline
d_2&\frac{\eta}{2}(2s_1+d_2-d_4)&\frac{\eta}{2}(2s_2+d_2-d_3)& 
\frac{\eta}{2}(2s_3+d_2-d_4)&\frac{\eta}{2}(2s_4+d_2-d_3)&0&d_2& 
\eta(-s_2-s_4+d_2+d_3)&\eta(-s_1-s_3+d_2+d_4)\\
\hline
d_3&\frac{\eta}{2}(2s_1+d_3-d_1)&\frac{\eta}{2}(2s_2+d_3-d_2)& 
\frac{\eta}{2}(2s_3+d_3-d_1)&\frac{\eta}{2}(2s_4+d_3-d_2)& 
\eta(-s_1-s_3+d_1+d_3)&\eta(-s_2-s_4+d_2+d_3)&d_3&0\\
\hline
d_4&\frac{\eta}{2}(2s_1+d_4-d_2)&\frac{\eta}{2}(2s_2+d_4-d_1)& 
\frac{\eta}{2}(2s_3+d_4-d_2)&\frac{\eta}{2}(2s_4+d_4-d_1)& 
\eta(-s_2-s_4+d_1+d_4)&\eta(-s_1-s_3+d_2+d_4)&0&d_4\\
\hline
\end{tabu}$}
\caption{The $8$-dimensional algebra $2Q_2(\eta)$}\label{t:algebra C5}
\endgroup
\end{center}
\end{table}
In particular, we note that $c$ and $d$ are not in $A$, and so $z$ is 
primitive in $A$.

The determinant of the Gram matrix of the Frobenius form on $A$ was 
computed in GAP, and it is $256\eta^3-192\eta^2+16$. The roots of 
this polynomial $\frac{1}{2}$ (of multiplicity $2$) and $-\frac{1}{4}$. 
So the only value of $\eta$ for which the algebra is not simple is, 
as in the case of diagram $C_4$, the value $\eta=-\frac{1}{4}$. 
When $\F$ is of characteristic $3$ or $5$, this cannot happen and so 
the algebra is always simple.

\bigskip
Finally, we need to consider the case $p=3$. Again, we can identify $a$, $b$, 
$c$, and $d$, respectively, with $b_{1,2}$, $c_{3,4}$, $b_{1,3}$, and $b_{2,4}$, 
using the notation $b_{i,j}=(i,j)=b_{j,i}$, $c_{i,j}=(e_i-e_j)(i,j)\neq 
c_{j,i}=(e_j-e_i)(i,j)$ established in Section~\ref{2-generated}, where we 
discussed the group $\hat G(3)$ and its Fischer space. 

Again we act by $\tau_x=\tau_a$, $\tau_y=\tau_b$, and $\tau_z=\tau_c\tau_d$ to 
find all axes contained in $A$. They are: 
\begin{center}
$\begin{aligned}
s_1:=&\,a=b_{1,2},\\
s_2:=&\,b^{\tau_z}=(c_{3,4}^{\tau_c})^{\tau_d}=c_{1,4}^{\tau_d}=c_{1,2},\\
s_3:=&\,s_2^{\tau_a}=c_{1,2}^{\tau_a}=c_{2,1},\\
s_4:=&\,a^{\tau_z}=(b_{1,2}^{\tau_c})^{\tau_d}=b_{2,3}^{\tau_d}=b_{3,4},\\
s_5:=&\,b=c_{3,4},\\
s_6:=&\,s_4^{\tau_b}=b_{3,4}^{\tau_b}=c_{4,3},\\
d_1:=&\,c+d=b_{1,3}+b_{2,4},\\
d_2:=&\,((c+d)^{\tau_a})^{\tau_b}=((b_{1,3}+b_{2,4})^{\tau_a})^{\tau_b}=
(b_{2,3}+b_{1,4})^{\tau_b}\\
&=c_{2,4}+c_{3,1}=c_{3,1}+c_{2,4},\\
d_3:=&\,d_2^{\tau_z}=((c_{3,1}+c_{2,4})^{\tau_c})^{\tau_d}=
(c_{1,3}+c_{2,4})^{\tau_d}=c_{1,3}+c_{4,2},\\
d_4:=&\,d_1^{\tau_a}=(b_{1,3}+b_{2,4})^{\tau_a}=b_{2,3}+b_{1,4}=
b_{1,4}+b_{2,3},\\
d_5:=&\,d_2^{\tau_a}=(c_{3,1}+c_{2,4})^{\tau_a}=c_{3,2}+c_{1,4}=
c_{1,4}+c_{3,2},\\
d_6:=&\,d_3^{\tau_a}=(c_{1,3}+c_{4,2})^{\tau_a}=c_{2,3}+c_{4,1}=
c_{4,1}+c_{2,3}.
\end{aligned}$
\end{center}

As in all previous cases, these axes have disjoint support and so they are 
linearly independent. The products are shown in Table~\ref{t:algebra C5-2}. 
It demonstrates that $A=\dla x,y,z\dra$ is 12-dimensional spanned by the above 
twelve axes. The axis $z$ is primitive in $A$ since $A$ does not contain $c$ 
and $d$.
\begin{sidewaystable}  
\begin{center}
\begingroup
\setlength{\tabcolsep}{10pt} %
\renewcommand{\arraystretch}{4}
\scalebox{.45}{
$\begin{tabu}[h!]{|c||c|c|c|c|c|c|c|c|c|c|c|c|}
\hline
& s_1 & s_2 & s_3 & s_4 & s_5 & s_6 & d_1 & d_2 & d_3 & d_4 & d_5 & d_6 \\ \hline\hline 
s_1 & s_1 &\frac{\eta}{2}(s_1+s_2-s_3) & \frac{\eta}{2}(s_1+s_3-s_2) & 0 & 0 & 0 & \frac{\eta}{2}(2s_1+d_1-d_4) & \frac{\eta}{2}(2s_1+d_2-d_5) & \frac{\eta}{2}(2s_1+d_3-d_6) & \frac{\eta}{2}(2s_1+d_4-d_1) & \frac{\eta}{2}(2s_1+d_5-d_2) & \frac{\eta}{2}(2s_1+d_6-d_3) \\ \hline
s_2 & \frac{\eta}{2}(s_1+s_2-s_3) & s_2  & \frac{\eta}{2}(s_2+s_3-s_1) & 0 & 0 & 0 & \frac{\eta}{2}(2s_2+d_1-d_5) & \frac{\eta}{2}(2s_2+d_2-d_6) & \frac{\eta}{2}(2s_2+d_3-d_4) & \frac{\eta}{2}(2s_2+d_4-d_3) & \frac{\eta}{2}(2s_2+d_5-d_1) & \frac{\eta}{2}(2s_2+d_6-d_2) \\ \hline
s_3 & \frac{\eta}{2}(s_1+s_3-s_2) & \frac{\eta}{2}(s_2+s_3-s_1) & s_3 & 0 & 0 & 0 & \frac{\eta}{2}(2s_3+d_1-d_6) & \frac{\eta}{2}(2s_3+d_2-d_4) & \frac{\eta}{2}(2s_3+d_3-d_5) & \frac{\eta}{2}(2s_3+d_4-d_2) & \frac{\eta}{2}(2s_3+d_5-d_3) & \frac{\eta}{2}(2s_3+d_6-d_1) \\ \hline
s_4 & 0 & 0 & 0 & s_4 & \frac{\eta}{2}(s_4+s_5-s_6) & \frac{\eta}{2}(s_4+s_6-s_5) & \frac{\eta}{2}(2s_4+d_1-d_4) & \frac{\eta}{2}(2s_4+d_2-d_6) & \frac{\eta}{2}(2s_4+d_3-d_5) & \frac{\eta}{2}(2s_4+d_4-d_1) & \frac{\eta}{2}(2s_4+d_5-d_3) & \frac{\eta}{2}(2s_4+d_6-d_2) \\ \hline
s_5 & 0 & 0 & 0 & \frac{\eta}{2}(s_4+s_5-s_6) & s_5 & \frac{\eta}{2}(s_5+s_6-s_4) & \frac{\eta}{2}(2s_5+d_1-d_5) & \frac{\eta}{2}(2s_5+d_2-d_4) & \frac{\eta}{2}(2s_5+d_3-d_6) & \frac{\eta}{2}(2s_5+d_4-d_2) & \frac{\eta}{2}(2s_5+d_5-d_1) & \frac{\eta}{2}(2s_5+d_6-d_3) \\ \hline
s_6 & 0 & 0 & 0 & \frac{\eta}{2}(s_4+s_6-s_5) & \frac{\eta}{2}(s_5+s_6-s_4) & s_6 & \frac{\eta}{2}(2s_6+d_1-d_6) & \frac{\eta}{2}(2s_6+d_2-d_5) & \frac{\eta}{2}(2s_6+d_3-d_4) & \frac{\eta}{2}(2s_6+d_4-d_3) & \frac{\eta}{2}(2s_6+d_5-d_2) & \frac{\eta}{2}(2s_6+d_6-d_1) \\ \hline
d_1 & \frac{\eta}{2}(2s_1+d_1-d_4) & \frac{\eta}{2}(2s_2+d_1-d_5) & \frac{\eta}{2}(2s_3+d_1-d_6) & \frac{\eta}{2}(2s_4+d_1-d_4) & \frac{\eta}{2}(2s_5+d_1-d_5) &\frac{\eta}{2}(2s_6+d_1-d_6) & d_1 & \frac{\eta}{2}(d_1+d_2-d_3) & \frac{\eta}{2}(d_1+d_3-d_2) & \eta(d_1+d_4-s_1-s_4) & \eta(d_1+d_5-s_2-s_5) & \eta(d_1+d_6-s_3-s_6) \\ \hline
d_2 & \frac{\eta}{2}(2s_1+d_2-d_5) & \frac{\eta}{2}(2s_2+d_2-d_6) & \frac{\eta}{2}(2s_3+d_2-d_4) & \frac{\eta}{2}(2s_4+d_2-d_6) & \frac{\eta}{2}(2s_5+d_2-d_4) &\frac{\eta}{2}(2s_6+d_2-d_5) & \frac{\eta}{2}(d_1+d_2-d_3) & d_2 & \frac{\eta}{2}(d_2+d_3-d_1) & \eta(d_2+d_4-s_3-s_5) & \eta(d_2+d_5-s_1-s_6) & \eta(d_2+d_6-s_2-s_4)  \\ \hline
d_3 & \frac{\eta}{2}(2s_1+d_3-d_6) & \frac{\eta}{2}(2s_2+d_3-d_4) & \frac{\eta}{2}(2s_3+d_3-d_5) & \frac{\eta}{2}(2s_4+d_3-d_5) & \frac{\eta}{2}(2s_5+d_3-d_6) &\frac{\eta}{2}(2s_6+d_3-d_4) & \frac{\eta}{2}(d_1+d_3-d_2) & \frac{\eta}{2}(d_2+d_3-d_1) & d_3 & \eta(d_3+d_4-s_2-s_6) & \eta(d_3+d_5-s_3-s_4) & \eta(d_3+d_6-s_1-s_5) \\ \hline
d_4 & \frac{\eta}{2}(2s_1+d_4-d_1) & \frac{\eta}{2}(2s_2+d_4-d_3) & \frac{\eta}{2}(2s_3+d_4-d_2) & \frac{\eta}{2}(2s_4+d_4-d_1) & \frac{\eta}{2}(2s_5+d_4-d_2) & \frac{\eta}{2}(2s_6+d_4-d_3) & \eta(d_1+d_4-s_1-s_4) & \eta(d_2+d_4-s_3-s_5) & \eta(d_3+d_4-s_2-s_6) & d_4  & \frac{\eta}{2}(d_4+d_5-d_6) & \frac{\eta}{2}(d_4+d_6-d_5) \\ \hline
d_5 & \frac{\eta}{2}(2s_1+d_5-d_2) & \frac{\eta}{2}(2s_2+d_5-d_1) & \frac{\eta}{2}(2s_3+d_5-d_3) & \frac{\eta}{2}(2s_4+d_5-d_3) & \frac{\eta}{2}(2s_5+d_5-d_1) & \frac{\eta}{2}(2s_6+d_5-d_2) & \eta(d_1+d_5-s_2-s_5) & \eta(d_2+d_5-s_1-s_6) & \eta(d_3+d_5-s_3-s_4) & \frac{\eta}{2}(d_4+d_5-d_6) & d_5  & \frac{\eta}{2}(d_5+d_6-d_4)  \\ \hline
d_6 & \frac{\eta}{2}(2s_1+d_6-d_3) & \frac{\eta}{2}(2s_2+d_6-d_2) & \frac{\eta}{2}(2s_3+d_6-d_1) & \frac{\eta}{2}(2s_4+d_6-d_2) & \frac{\eta}{2}(2s_5+d_6-d_3) & \frac{\eta}{2}(2s_6+d_6-d_1) & \eta(d_1+d_6-s_3-s_6) & \eta(d_2+d_6-s_2-s_4) & \eta(d_3+d_6-s_1-s_5) & \frac{\eta}{2}(d_4+d_6-d_5) & \frac{\eta}{2}(d_5+d_6-d_4) & d_6 \\ \hline
\end{tabu}$}
\caption{The 12-dimensional algebra $3Q_2(\eta)$}\label{t:algebra C5-2}
\endgroup
\end{center}
\end{sidewaystable}
%
The determinant of the Gram matrix of the Frobenius form on this algebra, computed in GAP, is 
$7\eta^{12}-111\eta^{11}+\frac{3051}{4}\eta^{10}-2939\eta^9+\frac{27153}{4}\eta^8-8964\eta^7+
4452\eta^6+5040\eta^5-10152\eta^4+6976\eta^3-1920\eta^2+64$. Its roots are $2$ (with 
multiplicity $8$), $\frac{1}{2}$ (with multiplicity $2$), $-\frac{1}{7}$, and $-1$. So the 
algebra is not simple only if $\eta=2$, $-\frac{1}{7}$, or $-1$. If $\F$ is of characteristic 
$3$ then $2=-\frac{1}{7}=-1=\frac{1}{2}$, so the algebra is always simple. If the 
characteristic of $\F$ is $5$ then $2=-\frac{1}{7}$, so for such fields, there are only two 
values of $\eta$, for which the the algebra is not simple.

\bigskip\noindent
{\bf Diagram $\mathbf{C_6}$:} The Coxeter group defined by the diagram $C_6$ is infinite. So we 
need to involve additional relations coming from the $3$-transposition property of $G=\la 
a,b,c,d\ra$. 

Consider $e:=a^{\tau_z}$, which is a single axis in $A$. Furthermore, $e$ is contained in 
$\dla a,c+d\dra$. According to our analysis of the diagram $A_3$, the latter algebra is the 
$4$-dimensional algebra $Q_2(\eta)$, and we see from Table \ref{4-dim} that the two single axes 
in this algebra, $a$ and $e$, are not collinear. In particular, $e\neq b$, since $b$ is 
collinear with $a$. Clearly, $A=\dla a,b,c+d\dra=\dla e,b,c+d\dra$. Indeed, we saw that 
$e=a^{\tau_z}\in A$ and, similarly, $a=e^{\tau_z}\in\dla e,b,c+d\dra$. If $e$ is non-collinear 
with $b$ then the diagram on the support set $\{e,b,c,d\}$ is $C_5$ and $A$ must be one of the 
algebras we obtained in that case. Thus, we may assume without loss of generality that $e$ is 
collinear with $b$. As an element of $G$, the axis $e$ is equal to $a^{cd}\in D$. So we get the 
additional relation $(a^{cd}b)^3=1$. 

Next consider the single axis $f:=a^{\tau_b}$, the third point on the line through $a$ and $b$. 
Clearly, $f\in A=\dla x,y,z\dra=\dla a,b,c+d\dra$. Note that $f\notin\{c,d\}$. Indeed, if, say, 
$f=c$ then $A$ contains $c$ and $z=c+d$, and hence all of $M_1(z)=\la c,d\ra$ lies in $A$, 
which contradicts primitivity of $z$. Thus, $f\notin\{c,d\}$. Since $A$ is primitive, so is 
also the subalgebra $\dla f,c+d\dra$. By Theorem \ref{symmetry}, the diagram on the support set 
$\{f,c,d\}$ has to be flip-symmetric. That is, either $f$ is collinear to both $c$ and $d$, or 
it is non-collinear to both $c$ and $d$. Note that $A=\dla a,b,c+d\dra=\dla f,b,c,d\dra$, since 
$a=f^{\tau_b}\in\dla f,b,c+d\dra$. If $f$ is non-collinear to $c$ and $d$ then the diagram 
on the support set $\{f,b,c,d\}$ is our $C_4$, and so it follows that $A=Q^3(\eta)$. Hence, we 
may assume without loss of generality that $f$ is collinear to $c$ and $d$. In terms of $G$, 
the axis $f$ is the involution $a^b\in D$, and so we get two additional relations: 
$(a^bc)^3=1=(a^bd)^3$.

Combining the Coxeter presentation coming with the diagram $C_6$ and the three additional 
relations $(a^{cd}b)^3=(a^bc)^3=(a^bd)^3=1$, we obtain via a coset enumeration in GAP that 
$G$ is a factor group of a group $\hat G\cong 2^{1+6}:SU_3(2)'$ in the notation of \cite{hs}, 
the (unique) irreducible 3-transposition group generated by four involutions and having the 
class of 3-transpositions of cardinality $36$. (This also follows from the presentation (A.8) 
in the appendix of \cite{hs}.) The subgroup $2^{1+6}$ is the unique minimal normal 
non-central subgroup. It is easy to check that $G$ cannot be the factor group $SU_3(2)'$, 
we can assume that $G=\hat G$.

We do not currently have a good description of the Fischer space of this group and it is a 
bit too big to do computations by hand in any case. Using GAP, we found that $A=\dla 
x,y,z\dra=\dla a,b,c+d\dra$ is a primitive 24-dimensional subalgebra of the $36$-dimensional 
Matsuo algebra of $G$. It is spanned by its $12$ single axes and $12$ double axes, which are 
as follows: 
\begin{center}
$\begin{aligned}
s_1:=&\,a, s_2:=a^{\tau_y},\\
s_3:=&\,b, s_4:=a^{\tau_z},\\
s_5:=&\,s_2^{\tau_z}, s_6:=s_4^{\tau_y},\\
s_7:=&\,s_3^{\tau_z}, s_8:=s_5^{\tau_x},\\
s_{9}:=&\,s_5^{\tau_y}, s_{10}:=s_8^{\tau_z},\\
s_{11}:=&\,s_7^{\tau_x}, s_{12}:=s_{11}^{\tau_y},\\
d_1:=&\,z^{\tau_x}, d_3:=z=c+d,\\
d_4:=&\,d_1^{\tau_y}, d_2:=d_4^{\tau_x},\\
d_5:=&\,d_3^{\tau_y}, d_6:=d_5^{\tau_x},\\
d_7:=&\,d_6^{\tau_z}, d_8:=d_4^{\tau_z},\\
d_9:=&\,d_7^{\tau_x}, d_{10}:=d_8^{\tau_x},\\
d_{11}:=&\,d_8^{\tau_y}, d_{12}:=d_{11}^{\tau_x},\\
\end{aligned}$
\end{center}
They are, clearly, forming a basis of $A$. The multiplication table of $A$, also found in GAP, 
is in Table~\ref{t:algebra C6-24}.

We also computed in GAP the determinant of the Gram matrix, which is $16777216\eta^9-
66060288\eta^8+113246208\eta^7-110100480\eta^6+66060288\eta^5-24772608\eta^4+5505024\eta^3-
589824\eta^2+4096$. The roots of this polynomial are $\frac{1}{2}$ (of multiplicity $8$) and 
$-\frac{1}{16}$. Hence the algebra is simple unless $\eta=-\frac{1}{16}$. In characteristic 
$3$, $-\frac{1}{16}=\frac{1}{2}$ and, in characteristic $17$, $-\frac{1}{16}=1$. Hence, for 
fields $\F$ in these two characteristics, the algebra is always simple.

\begin{center}
\begin{sidewaystable}  
\begingroup
\setlength{\tabcolsep}{20pt} %
\renewcommand{\arraystretch}{1.1}
\scalebox{.3}{
$\begin{tabu}[h!]{|c||c|c|c|c|c|c|c|c|c|c|c|c|c|c|c|c|c|c|c|c|c|c|c|c|}
\hline
&&&&&&&&&&&&&&&&&&&&&&&&\\
 & s_1 & s_4 & s_9 & s_{12} & s_2 & s_5 & s_{11} & s_6 & s_3 & s_7 & s_{10} & s_8 & d_1 & d_{11} & d_2 & d_{10} & d_5 & d_9 & d_3 & d_{12} & d_6 & d_7 & d_4 & d_8 \\ 
&&&&&&&&&&&&&&&&&&&&&&&&\\\hline \hline

&&&&&&&&&&&&&&&&&&&&&&&&\\
s_1 & s_1 & 0 & 0 & 0 & s_1+s_2-s_3 & s_1+s_5-s_8 & s_1+s_{11}-s_7 & s_1+s_6-s_{10} & s_1-s_2+s_3 & s_1-s_{11}+s_7 & 
s_1-s_6+s_{10} & s_1-s_5+s_8 & 2s_1+d_1-d_3 & 2s_1+d_{11}-d_{12} & 2s_1+d_2-d_4 & 2s_1+d_{10}-d_8 & 2s_1+d_5-d_6 & 2s_1+d_9-d_7 & 2s_1-d_1+d_3 & 2s_1-d_{11}+d_{12} & 
2s_1-d_5+d_6 & 2s_1-d_9+d_7 & 2s_1-d_2+d_4 & 2s_1-d_{10}+d_8 \\ 
&&&&&&&&&&&&&&&&&&&&&&&&\\\hline 

&&&&&&&&&&&&&&&&&&&&&&&&\\
s_4 & 0 & s_4 & 0 & 0 & s_4+s_2-s_{10} & s_4+s_5-s_7 & s_4+s_{11}-s_8 & s_4+s_6-s_3 & s_4-s_6+s_3 & s_4-s_5+s_7 & s_4-s_2+s_{10} & s_4-s_{11}+s_8 & 2s_4+d_1-d_3 & 2s_4+d_{11}-d_{12} & 
2s_4+d_2-d_8 & 2s_4+d_{10}-d_4 & 2s_4+d_5-d_7 & 2s_4+d_9-d_6 & 2s_4-d_1+d_3 & 
2s_4-d_{11}+d_{12} & 2s_4-d_9+d_6 & 2s_4-d_5+d_7 & 2s_4-d_{10}+d_4 & 2s_4-d_2+d_8 \\ 
&&&&&&&&&&&&&&&&&&&&&&&&\\\hline 

&&&&&&&&&&&&&&&&&&&&&&&&\\
s_9 & 0 & 0 & s_9 & 0 & s_9+s_2-s_8 & s_9+s_5-s_3 & s_9+s_{11}-s_{10} & s_9+s_6-s_7 & s_9-s_5+s_3 & s_9-s_6+s_7 & s_9-s_{11}+s_{10} & s_9-s_2+s_8 & 2s_9+d_1-d_{12} & 2s_9+d_{11}-d_3
 & 2s_9+d_2-d_8 & 2s_9+d_{10}-d_4 & 2s_9+d_5-d_6 & 2s_9+d_9-d_7 & 2s_9-d_{11}+d_3 & 2s_9-d_1+d_{12} & 2s_9-d_5+d_6 & 2s_9-d_9+d_7 & 2s_9-d_{10}+d_4 & 2s_9-d_2+d_8 \\ 
&&&&&&&&&&&&&&&&&&&&&&&&\\\hline 

&&&&&&&&&&&&&&&&&&&&&&&&\\
s_{12} & 0 & 0 & 0 & s_{12} & s_{12}+s_2-s_7 & s_{12}+s_5-s_{10} & s_{12}+s_{11}-s_3 & s_{12}+s_6-s_8 & s_{12}-s_{11}+s_3 & s_{12}-s_2+s_7 & s_{12}-s_5+s_{10} & s_{12}-s_6+s_8 & 2s_{12}+d_1-d_{12} & 2s_{12}+d_{11}-d_3 &
2s_{12}+d_2-d_4 & 2s_{12}+d_{10}-d_8 & 2s_{12}+d_5-d_7 & 2s_{12}+d_9-d_6 & 2s_{12}-d_{11}+d_3 & 2s_{12}-d_1+d_{12} & 2s_{12}-d_9+d_6 & 2s_{12}-d_5+d_7 & 2s_{12}-d_2+d_4 & 2s_{12}-d_{10}+d_8 \\ 
&&&&&&&&&&&&&&&&&&&&&&&&\\\hline 

&&&&&&&&&&&&&&&&&&&&&&&&\\
s_2 & s_1+s_2-s_3 & s_4+s_2-s_{10} & s_9+s_2-s_8 & s_{12}+s_2-s_7 & s_2 & 0 & 0 & 0 & -s_1+s_2+s_3 & -s_{12}+s_2+s_7 & -s_4+s_2+s_{10} & -s_9+s_2+s_8 & 2s_2+d_1-d_6 & 2s_2+d_{11}-d_7 &
2s_2+d_2-d_3 & 2s_2+d_{10}-d_{12} & 2s_2+d_5-d_4 & 2s_2+d_9-d_8 & 2s_2-d_2+d_3 & 2s_2-d_{10}+d_{12} & 2s_2-d_1+d_6 & 2s_2-d_{11}+d_7 & 2s_2-d_5+d_4 & 2s_2-d_9+d_8 \\ 
&&&&&&&&&&&&&&&&&&&&&&&&\\\hline 

&&&&&&&&&&&&&&&&&&&&&&&&\\
s_5 & s_1+s_5-s_8 & s_4+s_5-s_7 & s_9+s_5-s_3 & s_{12}+s_5-s_{10} & 0 & s_5 & 0 & 0 & -s_9+s_5+s_3 & -s_4+s_5+s_7 & -s_{12}+s_5+s_{10} & -s_1+s_5+s_8  & 2s_5+d_1-d_7 & 2s_5+d_{11}-d_6 & 2s_5+d_2-d_3 & 2s_5+d_{10}-d_{12} & 2s_5+d_5-d_8 & 2s_5+d_9-d_4 & 2s_5-d_2+d_3 & 2s_5-d_{10}+d_{12} & 2s_5-d_{11}+d_6 & 2s_5-d_1+d_7  & 2s_5-d_9+d_4 & 2s_5-d_5+d_8 \\ 
&&&&&&&&&&&&&&&&&&&&&&&&\\\hline 

&&&&&&&&&&&&&&&&&&&&&&&&\\
s_{11} & s_1+s_{11}-s_7 & s_4+s_{11}-s_8 & s_9+s_{11}-s_{10} & s_{12}+s_{11}-s_3 & 0 & 0 & s_{11} & 0 & -s_{12}+s_{11}+s_3 & -s_1+s_{11}+s_7 & -s_9+s_{11}+s_{10} & -s_4+s_{11}+s_8 & 2s_{11}+d_1-d_6 & 2s_{11}+d_{11}-d_7 & 2s_{11}+d_2-d_{12} & 2s_{11}+d_{10}-d_3 & 2s_{11}+d_5-d_8 & 2s_{11}+d_9-d_4 & 2s_{11}-d_{10}+d_3 & 2s_{11}-d_2+d_{12} & 2s_{11}-d_1+d_6 & 2s_{11}-d_{11}+d_7 & 2s_{11}-d_9+d_4 & 2s_{11}-d_5+d_8 \\ 
&&&&&&&&&&&&&&&&&&&&&&&&\\\hline 

&&&&&&&&&&&&&&&&&&&&&&&&\\
s_6 & s_1+s_6-s_{10} & s_4+s_6-s_3 & s_9+s_6-s_7 & s_{12}+s_6-s_8 & 0 & 0 & 0 & s_6 & -s_4+s_6+s_3 & -s_9+s_6+s_7 & -s_1+s_6+s_{10} & -s_{12}+s_6+s_8 & 2s_6+d_1-d_7 & 2s_6+d_{11}-d_6 & 2s_6+d_2-d_{12} & 2s_6+d_{10}-d_3 & 2s_6+d_5-d_4 & 2s_6+d_9-d_8 & 2s_6-d_{10}+d_3 & 2s_6-d_2+d_{12} & 2s_6-d_{11}+d_6 & 2s_6-d_1+d_7 & 2s_6-d_5+d_4 & 2s_6-d_9+d_8 \\ 
&&&&&&&&&&&&&&&&&&&&&&&&\\\hline 

&&&&&&&&&&&&&&&&&&&&&&&&\\
s_3 & s_1-s_2+s_3 & s_4-s_6+s_3 & s_9-s_5+s_3 & s_{12}-s_{11}+s_3 & -s_1+s_2+s_3 & -s_9+s_5+s_3 & -s_{12}+s_{11}+s_3 & -s_4+s_6+s_3 & s_3 & 0 & 0 & 0 & 2s_3+d_1-d_4 & 2s_3+d_{11}-d_8 & 2s_3+d_2-d_6 & 2s_3+d_{10}-d_7 & 2s_3+d_5-d_3 & 2s_3+d_9-d_{12} & 2s_3-d_5+d_3 & 2s_3-d_9+d_{12} & 2s_3-d_2+d_6 & 2s_3-d_{10}+d_7 & 2s_3-d_1+d_4 & 2s_3-d_{11}+d_8 \\ 
&&&&&&&&&&&&&&&&&&&&&&&&\\\hline 

&&&&&&&&&&&&&&&&&&&&&&&&\\
s_7 & s_1-s_{11}+s_7 & s_4-s_5+s_7 & s_9-s_6+s_7 & s_{12}-s_2+s_7 & -s_{12}+s_2+s_7 & -s_4+s_5+s_7 & -s_1+s_{11}+s_7 & -s_9+s_6+s_7 & 0 & s_7 & 0 & 0 & 2s_7+d_1-d_8 & 
2s_7+d_{11}-d_4 & 2s_7+d_2-d_7 & 2s_7+d_{10}-d_6 & 2s_7+d_5-d_3 & 2s_7+d_9-d_{12} & 2s_7-d_5+d_3 & 2s_7-d_9+d_{12} & 2s_7-d_{10}+d_6 & 2s_7-d_2+d_7 & 2s_7-d_{11}+d_4 & 2s_7-d_1+d_8 \\ 
&&&&&&&&&&&&&&&&&&&&&&&&\\\hline 

&&&&&&&&&&&&&&&&&&&&&&&&\\
s_{10} & s_1-s_6+s_{10} & s_4-s_2+s_{10} & s_9-s_{11}+s_{10} & s_{12}-s_5+s_{10} & -s_4+s_2+s_{10} & -s_{12}+s_5+s_{10} & -s_9+s_{11}+s_{10} & -s_1+s_6+s_{10} & 0 & 0 & s_{10} & 0 & 2s_{10}+d_1-d_8 & 2s_{10}+d_{11}-d_4 & 2s_{10}+d_2-d_6 & 2s_{10}+d_{10}-d_7 & 2s_{10}+d_5-d_{12} & 2s_{10}+d_9-d_3 & 2s_{10}-d_9+d_3 & 2s_{10}-d_5+d_{12} & 2s_{10}-d_2+d_6 & 2s_{10}-d_{10}+d_7 & 2s_{10}-d_{11}+d_4 & 2s_{10}-d_1+d_8 \\ 
&&&&&&&&&&&&&&&&&&&&&&&&\\\hline 

&&&&&&&&&&&&&&&&&&&&&&&&\\
s_8 & s_1-s_5+s_8 & s_4-s_{11}+s_8 & s_9-s_2+s_8 & s_{12}-s_6+s_8 & -s_9+s_2+s_8 & -s_1+s_5+s_8 & -s_4+s_{11}+s_8 & -s_{12}+s_6+s_8 & 0 & 0 & 0 & s_8 & 2s_8+d_1-d_4 & 2s_8+d_{11}-d_8 & 2s_8+d_2-d_7 & 2s_8+d_{10}-d_6 & 2s_8+d_5-d_{12} & 2s_8+d_9-d_3 & 2s_8-d_9+d_3 & 2s_8-d_5+d_{12} & 2s_8-d_{10}+d_6 & 2s_8-d_2+d_7 & 2s_8-d_1+d_4 & 2s_8-d_{11}+d_8 \\ 
&&&&&&&&&&&&&&&&&&&&&&&&\\\hline 

&&&&&&&&&&&&&&&&&&&&&&&&\\
d_1 & 2s_1+d_1-d_3 & 2s_4+d_1-d_3 & 2s_9+d_1-d_{12} & 2s_{12}+d_1-d_{12} & 2s_2+d_1-d_6 & 2s_5+d_1-d_7 & 2s_{11}+d_1-d_6 & 2s_6+d_1-d_7 & 2s_3+d_1-d_4 & 2s_7+d_1-d_8 & 2s_{10}+d_1-d_8 & 2s_8+d_1-d_4 & 
d_1 & 0 & 
\begin{tabu}{@{}c@{}} 2(d_1+d_2) \\ -(d_5+d_9) \end{tabu} & 
\begin{tabu}{@{}c@{}} 2(d_1+d_{10}) \\ -(d_5+d_9) \end{tabu} & 
\begin{tabu}{@{}c@{}} 2(d_1+d_5) \\ -(d_2+d_{10}) \end{tabu} & 
\begin{tabu}{@{}c@{}} 2(d_1+d_9) \\ -(d_2+d_{10}) \end{tabu} & 
\begin{tabu}{@{}c@{}} 2(d_1+d_3) \\ -2(s_1+s_4) \end{tabu} & 
\begin{tabu}{@{}c@{}} 2(d_1+d_{12}) \\ -2(s_9+s_{12}) \end{tabu} & 
\begin{tabu}{@{}c@{}} 2(d_1+d_6) \\ -2(s_2+s_{11}) \end{tabu} & 
\begin{tabu}{@{}c@{}} 2(d_1+d_7) \\ -2(s_5+s_6) \end{tabu} & 
\begin{tabu}{@{}c@{}} 2(d_1+d_4) \\ -2(s_3+s_8) \end{tabu} &  
\begin{tabu}{@{}c@{}} 2(d_1+d_8) \\ -2(s_7+s_{10}) \end{tabu} \\ \hline 

d_{11} & 2s_1+d_{11}-d_{12} & 2s_4+d_{11}-d_{12} & 2s_9+d_{11}-d_3 & 2s_{12}+d_{11}-d_3 & 2s_2+d_{11}-d_7 & 2s_5+d_{11}-d_6 & 2s_{11}+d_{11}-d_7 & 2s_6+d_{11}-d_6 & 2s_3+d_{11}-d_8 & 2s_7+d_{11}-d_4 & 2s_{10}+d_{11}-d_4 & 2s_8+d_{11}-d_8 & 
0 & d_{11} & 
\begin{tabu}{@{}c@{}} 2(d_{11}+d_2) \\ -(d_5+d_9) \end{tabu} & 
\begin{tabu}{@{}c@{}} 2(d_{11}+d_{10}) \\ -(d_5+d_9) \end{tabu} & 
\begin{tabu}{@{}c@{}} 2(d_{11}+d_5) \\ -(d_2+d_{10}) \end{tabu} & 
\begin{tabu}{@{}c@{}} 2(d_{11}+d_9) \\ -(d_2+d_{10}) \end{tabu} & 
\begin{tabu}{@{}c@{}} 2(d_{11}+d_3) \\ -2(s_9+s_{12}) \end{tabu} & 
\begin{tabu}{@{}c@{}} 2(d_{11}+d_{12}) \\ -2(s_1+s_4) \end{tabu} & 
\begin{tabu}{@{}c@{}} 2(d_{11}+d_6) \\ -2(s_5+s_6) \end{tabu} & 
\begin{tabu}{@{}c@{}} 2(d_{11}+d_7) \\ -2(s_2+s_{11}) \end{tabu} & 
\begin{tabu}{@{}c@{}} 2(d_{11}+d_4) \\ -2(s_7+s_{10}) \end{tabu} & 
\begin{tabu}{@{}c@{}} 2(d_{11}+d_8) \\ -2(s_3+s_8) \end{tabu}  \\ \hline 

d_2 & 2s_1+d_2-d_4 & 2s_4+d_2-d_8 & 2s_9+d_2-d_8 & 2s_{12}+d_2-d_4 & 2s_2+d_2-d_3 & 2s_5+d_2-d_3 & 2s_{11}+d_2-d_{12} & 2s_6+d_2-d_{12} & 2s_3+d_2-d_6 & 2s_7+d_2-d_7 & 2s_{10}+d_2-d_6 & 2s_8+d_2-d_7 & 
\begin{tabu}{@{}c@{}} 2(d_1+d_2) \\ -(d_5+d_9) \end{tabu} & 
\begin{tabu}{@{}c@{}} 2(d_{11}+d_2) \\ -(d_5+d_9) \end{tabu} & 
d_2 & 0 & 
\begin{tabu}{@{}c@{}} 2(d_2+d_5) \\ -(d_1+d_{11}) \end{tabu} & 
\begin{tabu}{@{}c@{}} 2(d_2+d_9) \\ -(d_1+d_{11}) \end{tabu} & 
\begin{tabu}{@{}c@{}} 2(d_2+d_3) \\ -2(s_2+s_5) \end{tabu} & 
\begin{tabu}{@{}c@{}} 2(d_2+d_{12}) \\ -2(s_{11}+s_6) \end{tabu} & 
\begin{tabu}{@{}c@{}} 2(d_2+d_6) \\ -2(s_3+s_{10}) \end{tabu} & 
\begin{tabu}{@{}c@{}} 2(d_2+d_7) \\ -2(s_7+s_8) \end{tabu} & 
\begin{tabu}{@{}c@{}} 2(d_2+d_4) \\ -2(s_1+s_{12}) \end{tabu} & 
\begin{tabu}{@{}c@{}} 2(d_2+d_8) \\ -2(s_4+s_9) \end{tabu} \\ \hline 

d_{10} & 2s_1+d_{10}-d_8 & 2s_4+d_{10}-d_4 & 2s_9+d_{10}-d_4 & 2s_{12}+d_{10}-d_8 & 2s_2+d_{10}-d_{12} & 2s_5+d_{10}-d_{12} & 2s_{11}+d_{10}-d_3 & 2s_6+d_{10}-d_3 & 2s_3+d_{10}-d_7 & 2s_7+d_{10}-d_6 & 2s_{10}+d_{10}-d_7 & 2s_8+d_{10}-d_6 & 
\begin{tabu}{@{}c@{}} 2(d_1+d_{10}) \\ -(d_5+d_9) \end{tabu} & 
\begin{tabu}{@{}c@{}} 2(d_{11}+d_{10}) \\ -(d_5+d_9) \end{tabu} & 
0 & d_{10} & 
\begin{tabu}{@{}c@{}} 2(d_{10}+d_5) \\ -(d_1+d_{11}) \end{tabu}& 
\begin{tabu}{@{}c@{}} 2(d_{10}+d_9) \\ -(d_1+d_{11}) \end{tabu} & 
\begin{tabu}{@{}c@{}} 2(d_{10}+d_3) \\ -2(s_{11}+s_6) \end{tabu} & 
\begin{tabu}{@{}c@{}} 2(d_{10}+d_{12}) \\ -2(s_2+s_5) \end{tabu} & 
\begin{tabu}{@{}c@{}} 2(d_{10}+d_6) \\ -2(s_7+s_8) \end{tabu} & 
\begin{tabu}{@{}c@{}} 2(d_{10}+d_7)  \\ -2(s_3+s_{10}) \end{tabu} & 
\begin{tabu}{@{}c@{}} 2(d_{10}+d_4) \\ -2(s_4+s_9) \end{tabu} & 
\begin{tabu}{@{}c@{}} 2(d_{10}+d_8) \\ -2(s_1+s_{12}) \end{tabu} \\ \hline 

d_5 & 2s_1+d_5-d_6 & 2s_4+d_5-d_7 & 2s_9+d_5-d_6 & 2s_{12}+d_5-d_7 & 2s_2+d_5-d_4 & 2s_5+d_5-d_8 & 2s_{11}+d_5-d_8 & 2s_6+d_5-d_4 & 2s_3+d_5-d_3 & 2s_7+d_5-d_3 & 2s_{10}+d_5-d_{12} & 2s_8+d_5-d_{12} & 
\begin{tabu}{@{}c@{}} 2(d_1+d_5) \\ -(d_2+d_{10}) \end{tabu} & 
\begin{tabu}{@{}c@{}} 2(d_{11}+d_5) \\ -(d_2+d_{10}) \end{tabu} & 
\begin{tabu}{@{}c@{}} 2(d_2+d_5) \\ -(d_1+d_{11}) \end{tabu} & 
\begin{tabu}{@{}c@{}} 2(d_{10}+d_5) \\ -(d_1+d_{11}) \end{tabu} & 
d_5 & 0 & 
\begin{tabu}{@{}c@{}} 2(d_5+d_3) \\ -2(s_3+s_7) \end{tabu} & 
\begin{tabu}{@{}c@{}} 2(d_5+d_{12}) \\ -2(s_{10}+s_8) \end{tabu} & 
\begin{tabu}{@{}c@{}} 2(d_5+d_6) \\ -2(s_1+s_9) \end{tabu} & 
\begin{tabu}{@{}c@{}} 2(d_5+d_7) \\ -2(s_4+s_{12}) \end{tabu} & 
\begin{tabu}{@{}c@{}} 2(d_5+d_4) \\ -2(s_2+s_6) \end{tabu} & 
\begin{tabu}{@{}c@{}} 2(d_5+d_8) \\ -2(s_5+s_{11}) \end{tabu} \\ \hline 

d_9 & 2s_1+d_9-d_7 & 2s_4+d_9-d_6 & 2s_9+d_9-d_7 & 2s_{12}+d_9-d_6 & 2s_2+d_9-d_8 & 2s_5+d_9-d_4 & 2s_{11}+d_9-d_4 & 2s_6+d_9-d_8 & 2s_3+d_9-d_{12} & 2s_7+d_9-d_{12} & 2s_{10}+d_9-d_3 & 2s_8+d_9-d_3 & 
\begin{tabu}{@{}c@{}} 2(d_1+d_9) \\ -(d_2+d_{10}) \end{tabu} & 
\begin{tabu}{@{}c@{}} 2(d_{11}+d_9) \\ -(d_2+d_{10}) \end{tabu} & 
\begin{tabu}{@{}c@{}} 2(d_2+d_9) \\ -(d_1+d_{11}) \end{tabu} & 
\begin{tabu}{@{}c@{}} 2(d_{10}+d_9) \\ -(d_1+d_{11}) \end{tabu} & 
0 & d_9 & 
\begin{tabu}{@{}c@{}} 2(d_9+d_3) \\ -2(s_{10}+s_8) \end{tabu} & 
\begin{tabu}{@{}c@{}} 2(d_9+d_{12}) \\ -2(s_3+s_7) \end{tabu} & 
\begin{tabu}{@{}c@{}} 2(d_9+d_6) \\ -2(s_4+s_{12}) \end{tabu} & 
\begin{tabu}{@{}c@{}} 2(d_9+d_7) \\ -2(s_1+s_9) \end{tabu} & 
\begin{tabu}{@{}c@{}} 2(d_9+d_4) \\ -2(s_5+s_{11}) \end{tabu} & 
\begin{tabu}{@{}c@{}} 2(d_9+d_8) \\ -2(s_2+s_6) \end{tabu} \\ \hline 

d_3 & 2s_1-d_1+d_3 & 2s_4-d_1+d_3 & 2s_9-d_{11}+d_3 & 2s_{12}-d_{11}+d_3 & 2s_2-d_2+d_3 & 2s_5-d_2+d_3 & 2s_{11}-d_{10}+d_3 & 2s_6-d_{10}+d_3 & 2s_3-d_5+d_3 & 2s_7-d_5+d_3 & 2s_{10}-d_9+d_3 & 2s_8-d_9+d_3 & 
\begin{tabu}{@{}c@{}} 2(d_1+d_3) \\ -2(s_1+s_4) \end{tabu} & 
\begin{tabu}{@{}c@{}} 2(d_{11}+d_3) \\ -2(s_9+s_{12}) \end{tabu} & 
\begin{tabu}{@{}c@{}} 2(d_2+d_3) \\ -2(s_2+s_5) \end{tabu} & 
\begin{tabu}{@{}c@{}} 2(d_{10}+d_3) \\ -2(s_{11}+s_6) \end{tabu} & 
\begin{tabu}{@{}c@{}} 2(d_5+d_3) \\ -2(s_3+s_7) \end{tabu} & 
\begin{tabu}{@{}c@{}} 2(d_9+d_3) \\ -2(s_{10}+s_8) \end{tabu} & 
d_3 & 0 & 
\begin{tabu}{@{}c@{}} 2(d_3+d_6) \\ -(d_4+d_8) \end{tabu} & 
\begin{tabu}{@{}c@{}} 2(d_3+d_7) \\ -(d_4+d_8) \end{tabu} & 
\begin{tabu}{@{}c@{}} 2(d_3+d_4) \\ -(d_6+d_7) \end{tabu} & 
\begin{tabu}{@{}c@{}} 2(d_3+d_8) \\ -(d_6+d_7) \end{tabu} \\ \hline 

d_{12} & 2s_1-d_{11}+d_{12} & 2s_4-d_{11}+d_{12} & 2s_9-d_1+d_{12} & 2s_{12}-d_1+d_{12} & 2s_2-d_{10}+d_{12} & 2s_5-d_{10}+d_{12} & 2s_{11}-d_2+d_{12} & 2s_6-d_2+d_{12} & 2s_3-d_9+d_{12} & 2s_7-d_9+d_{12} & 2s_{10}-d_5+d_{12} & 2s_8-d_5+d_{12} & 
\begin{tabu}{@{}c@{}} 2(d_1+d_{12}) \\ -2(s_9+s_{12}) \end{tabu} & 
\begin{tabu}{@{}c@{}} 2(d_{11}+d_{12}) \\ -2(s_1+s_4) \end{tabu} & 
\begin{tabu}{@{}c@{}} 2(d_2+d_{12}) \\ -2(s_{11}+s_6) \end{tabu} & 
\begin{tabu}{@{}c@{}} 2(d_{10}+d_{12}) \\ -2(s_2+s_5) \end{tabu} & 
\begin{tabu}{@{}c@{}} 2(d_5+d_{12}) \\ -2(s_{10}+s_8) \end{tabu} & 
\begin{tabu}{@{}c@{}} 2(d_9+d_{12}) \\ -2(s_3+s_7) \end{tabu} & 
0 & d_{12} & 
\begin{tabu}{@{}c@{}} 2(d_{12}+d_6) \\ -(d_4+d_8) \end{tabu} & 
\begin{tabu}{@{}c@{}} 2(d_{12}+d_7) \\ -(d_4+d_8) \end{tabu} & 
\begin{tabu}{@{}c@{}} 2(d_{12}+d_4) \\ -(d_6+d_7) \end{tabu} & 
\begin{tabu}{@{}c@{}} 2(d_{12}+d_8) \\ -(d_6+d_7) \end{tabu} \\ \hline 

d_6 & 2s_1-d_5+d_6 & 2s_4-d_9+d_6 & 2s_9-d_5+d_6 & 2s_{12}-d_9+d_6 & 2s_2-d_1+d_6 & 2s_5-d_{11}+d_6 & 2s_{11}-d_1+d_6 & 2s_6-d_{11}+d_6 & 2s_3-d_2+d_6 & 2s_7-d_{10}+d_6 & 2s_{10}-d_2+d_6 & 2s_8-d_{10}+d_6 & 
\begin{tabu}{@{}c@{}} 2(d_1+d_6) \\ -2(s_2+s_{11}) \end{tabu} & 
\begin{tabu}{@{}c@{}} 2(d_{11}+d_6) \\ -2(s_5+s_6) \end{tabu} & 
\begin{tabu}{@{}c@{}} 2(d_2+d_6) \\ -2(s_3+s_{10}) \end{tabu} & 
\begin{tabu}{@{}c@{}} 2(d_{10}+d_6) \\ -2(s_7+s_8) \end{tabu} & 
\begin{tabu}{@{}c@{}} 2(d_5+d_6) \\ -2(s_1+s_9) \end{tabu} & 
\begin{tabu}{@{}c@{}} 2(d_9+d_6) \\ -2(s_4+s_{12}) \end{tabu} & 
\begin{tabu}{@{}c@{}} 2(d_3+d_6) \\ -(d_4+d_8) \end{tabu} & 
\begin{tabu}{@{}c@{}} 2(d_{12}+d_6) \\ -(d_4+d_8) \end{tabu} & 
d_6 & 0 & 
\begin{tabu}{@{}c@{}} 2(d_6+d_4) \\ -(d_3+d_{12}) \end{tabu} & 
\begin{tabu}{@{}c@{}} 2(d_6+d_8) \\ -(d_3+d_{12}) \end{tabu} \\ \hline 

d_7 & 2s_1-d_9+d_7 & 2s_4-d_5+d_7 & 2s_9-d_9+d_7 & 2s_{12}-d_5+d_7 & 2s_2-d_{11}+d_7 & 2s_5-d_1+d_7 & 2s_{11}-d_{11}+d_7 & 2s_6-d_1+d_7 & 2s_3-d_{10}+d_7 & 2s_7-d_2+d_7 & 2s_{10}-d_{10}+d_7 & 2s_8-d_2+d_7 & 
\begin{tabu}{@{}c@{}} 2(d_1+d_7) \\ -2(s_5+s_6) \end{tabu} & 
\begin{tabu}{@{}c@{}} 2(d_{11}+d_7) \\ -2(s_2+s_{11}) \end{tabu} & 
\begin{tabu}{@{}c@{}} 2(d_2+d_7) \\ -2(s_7+s_8) \end{tabu} & 
\begin{tabu}{@{}c@{}} 2(d_{10}+d_7) \\ -2(s_3+s_{10}) \end{tabu} & 
\begin{tabu}{@{}c@{}} 2(d_5+d_7) \\ -2(s_4+s_{12}) \end{tabu} & 
\begin{tabu}{@{}c@{}} 2(d_9+d_7) \\ -2(s_1+s_9) \end{tabu} & 
\begin{tabu}{@{}c@{}} 2(d_3+d_7) \\ -(d_4+d_8) \end{tabu} & 
\begin{tabu}{@{}c@{}} 2(d_{12}+d_7) \\ -(d_4+d_8) \end{tabu} & 
0 & d_7 & 
\begin{tabu}{@{}c@{}} 2(d_7+d_4) \\ -(d_3+d_{12}) \end{tabu} & 
\begin{tabu}{@{}c@{}} 2(d_7+d_8) \\ -(d_3+d_{12}) \end{tabu} \\ \hline 

d_4 & 2s_1-d_2+d_4 & 2s_4-d_{10}+d_4 & 2s_9-d_{10}+d_4 & 2s_{12}-d_2+d_4 & 2s_2-d_5+d_4 & 2s_5-d_9+d_4 & 2s_{11}-d_9+d_4 & 2s_6-d_5+d_4 & 2s_3-d_1+d_4 & 2s_7-d_{11}+d_4 & 2s_{10}-d_{11}+d_4 & 2s_8-d_1+d_4 & 
\begin{tabu}{@{}c@{}} 2(d_1+d_4) \\ -2(s_3+s_8) \end{tabu} & 
\begin{tabu}{@{}c@{}} 2(d_{11}+d_4) \\ -2(s_7+s_{10}) \end{tabu} & 
\begin{tabu}{@{}c@{}} 2(d_2+d_4) \\ -2(s_1+s_{12}) \end{tabu} & 
\begin{tabu}{@{}c@{}} 2(d_{10}+d_4) \\ -2(s_4+s_9) \end{tabu} & 
\begin{tabu}{@{}c@{}} 2(d_5+d_4) \\ -2(s_2+s_6) \end{tabu} & 
\begin{tabu}{@{}c@{}} 2(d_9+d_4) \\ -2(s_5+s_{11}) \end{tabu}  & 
\begin{tabu}{@{}c@{}} 2(d_3+d_4) \\ -(d_6+d_7) \end{tabu} & 
\begin{tabu}{@{}c@{}} 2(d_{12}+d_4) \\ -(d_6+d_7) \end{tabu} & 
\begin{tabu}{@{}c@{}} 2(d_6+d_4) \\ -(d_3+d_{12}) \end{tabu} & 
\begin{tabu}{@{}c@{}} 2(d_7+d_4) \\ -(d_3+d_{12}) \end{tabu} & 
d_4 & 0 \\ \hline 

d_8 & 2s_1-d_{10}+d_8 & 2s_4-d_2+d_8 & 2s_9-d_2+d_8 & 2s_{12}-d_{10}+d_8 & 2s_2-d_9+d_8 & 2s_5-d_5+d_8 & 2s_{11}-d_5+d_8 & 2s_6-d_9+d_8 & 2s_3-d_{11}+d_8 & 2s_7-d_1+d_8 & 2s_{10}-d_1+d_8 & 
2s_8-d_{11}+d_8 & 
\begin{tabu}{@{}c@{}} 2(d_1+d_8) \\ -2(s_7+s_{10}) \end{tabu} & 
\begin{tabu}{@{}c@{}} 2(d_{11}+d_8) \\ -2(s_3+s_8) \end{tabu} & 
\begin{tabu}{@{}c@{}} 2(d_2+d_8) \\ -2(s_4+s_9) \end{tabu} & 
\begin{tabu}{@{}c@{}} 2(d_{10}+d_8) \\ -2(s_1+s_{12}) \end{tabu} & 
\begin{tabu}{@{}c@{}} 2(d_5+d_8) \\ -2(s_5+s_{11}) \end{tabu} & 
\begin{tabu}{@{}c@{}} 2(d_9+d_8) \\ -2(s_2+s_6) \end{tabu} & 
\begin{tabu}{@{}c@{}} 2(d_3+d_8) \\ -(d_6-d_7) \end{tabu} & 
\begin{tabu}{@{}c@{}} 2(d_{12}+d_8) \\ -(d_6-d_7) \end{tabu} & 
\begin{tabu}{@{}c@{}} 2(d_6+d_8) \\ -(d_3-d_{12}) \end{tabu} & 
\begin{tabu}{@{}c@{}} 2(d_7+d_8) \\ -(d_3-d_{12}) \end{tabu} & 
0 & d_8 \\ \hline 
\end{tabu}$}
\caption{The 24-dimensional algebra of type $C_6$ (all non-diagonal entries 
have an additional factor $\frac{\eta}{2}$ omitted here)}\label{t:algebra C6-24}
\endgroup
\end{sidewaystable}  
\end{center}


\section{Flip subalgebras}

\label{flip subalgebras}

Let us return to Theorem \ref{symmetry}, stating that a subalgebra in a Matsuo 
algebra $M$, generated by a set of single and double axes $Y$, can only be 
primitive when the diagram on the support of $Y$ admits a specific 
automorphism of order two, the flip. In this section, we shift focus from the 
set of generators (and in particular, from the number of generators) onto the 
flip itself. In doing so, we uncover an ``industrial'' method of building primitive 
algebras of Monster type $(2\eta,\eta)$.

As before, let $M=M_\eta(\Gm)$ be a Matsuo algebra over a field $\F$, 
corresponding to the Fischer space $\Gm$ coming from a $3$-transposition 
group $(G,D)$. The group $\Aut\Gm$ acts on $M$ and it is the full group of 
automorphisms of $M$ preserving the basis $D$ consisting of all single axes. 

\subsection{Fixed subalgebra}

Suppose $H\leq\Aut\Gm$. 

\begin{definition}
The \emph{fixed subalgebra} of $H$ in $M$ is defined as
$$M_H:=\{u\in M\mid u^h=u\mbox{ for all }h\in H\}.$$
\end{definition}

Clearly, $M_H$ is a subspace and it is closed for the algebra product: if 
$u,v\in M_H$ then $(uv)^h=u^hv^h=uv$ for each $h\in H$. Thus, $M_H$ is indeed 
a subalgebra of $M$. 

\begin{definition}
For an $H$-orbit $O\subseteq D$, define the \emph{orbit vector} as 
$$u_O:=\sum_{c\in O}c.$$ 
\end{definition}

Manifestly, $u_O$ is fixed by $H$ and so $u_O\in M_H$. In fact, the orbit 
vectors provide us with a nice basis of $M_H$.

\begin{proposition} \label{basis}
Let $O_1,\ldots,O_k$ be all the $H$-orbits on $D$. Then the vectors 
$u_i:=u_{O_i}$ form a basis of the fixed subalgebra $M_H$.
\end{proposition}

\begin{proof}
First of all, the vectors $u_1,\ldots,u_k$ are linearly independent as their 
supports in the basis $D$ of $M$ ({\it i.e.}, the orbits $O_i$) partition $D$ 
and hence are disjoint. 

Consider an arbitrary $u\in M_H$, say, $u=\sum_{c\in D}\al_cc$ for some 
scalars $\al_c\in\F$. If $d,e\in D$ are in the same $H$-orbit, then $e=d^h$ 
for some $h\in H$, and so $u=u^h=\sum_{c\in D}\al_cc^h$ and this shows that 
$\al_d=\al_e$. That is, the coefficients $\al_c$ stay constant on the 
$H$-orbits $O_i$. Clearly, this means that, for each orbit $O_i$, there is a 
single value $\al_i$ such that $\al_c=\al_i$ for all $c\in O_i$. 
Consequently, $u$ can be written as $u=\sum_{i=1}^k\al_i u_i$, and so 
$u_1,\ldots,u_k$ also spans $M_H$. 
\end{proof}

\subsection{Flip subalgebra}

Select an element $\tau\in\Aut\Gm$ with $|\tau|=2$. Let $H=\la\tau\ra$. Then 
the $H$-orbits $O_1,\ldots,O_k$ have each length $1$ or $2$. Let us classify 
the orbit vectors $u_i$ into three groups:
\begin{itemize}
\item $u_i=a\in D$, corresponding to orbits $O_i=\{a\}$ of length $1$;
\item $u_i=a+b$, corresponding to orbits $O_i=\{a,b\}$ with $ab=0$ (orthogonal 
orbits);
\item $u_i=a+b$, with $O_i=\{a,b\}$ satisfying $ab\neq 0$ (non-orthogonal 
orbits). 
\end{itemize}

We call the $u_i$ of the first kind \emph{singles}, of the second kind 
\emph{doubles}, and of the third kind \emph{extras}. 

It is easy to see that singles are simply all single axes from $M$ that are 
contained in $M_H$. 

\begin{proposition} \label{primitive}
Each double $u_i=a+b$ is a double axis and it is primitive in $M_H$. In fact, 
the set of doubles consists of all double axes that are contained in $M_H$ and 
are primitive in it.
\end{proposition}

\begin{proof}
Manifestly, $u_i$ is a double axis, so we just need to see that $u_i$ is 
primitive in $B:=M_H$. Note that $B_1(u_i)=M_1(u_i)\cap B=\la a,b\ra\cap M_H$. 
We know that $u_i=a+b\in M_H$. On the other hand, $a\not\in M_H$, because 
$a^\tau=b\neq a$. Thus, $B_1(u_i)=\la u_i\ra$ is $1$-dimensional and so $u_i$ 
is indeed primitive in $B=M_H$.

For the second claim, if a double axis $x=c+d$ is contained in $M_H$ then $H$ 
acts on the support $\{c,d\}$ of $x$. If the action is trivial then both $c$ 
and $d$ are in $M_H$, and so $x$ is not primitive in $M_H$. On the other hand, 
if $c^\tau=d$ then $\{c,d\}$ is one of the orbits $O_i$ and $x=u_i$. Hence the 
claim follows.
\end{proof}

We are ready for the key definition.

\begin{definition}
The \emph{flip subalgebra} $A=A(\tau)$ corresponding to an involution 
$\tau\in\Aut\Gm$ is the subalgebra of $M_H$ (where $H=\la\tau\ra$) generated 
by all singles and all doubles.
\end{definition} 

Singles are Matsuo axes and so they are primitive in all of $M$. According to 
Proposition \ref{primitive}, doubles are primitive in $M_H$ and so they are 
also primitive in $A(\tau)$. Thus, each $A=A(\tau)$ is a primitive algebra of 
Monster type $(2\eta,\eta)$. 

Let $Y$ be the set of generators of $A$, that is, $Y$ consists of all singles 
and all doubles. We note that $\tau$ is the flip of the diagram on 
$Z=\supp(Y)$, and so we will call $\tau$ the flip, even though $\tau$ is an 
arbitrary involution in $\Aut\Gm$.

If $\tau$ and $\tau'$ are conjugate in $\Aut\Gm$, say $\tau'=\tau^\sg$ for 
some $\sg\in\Aut\Gm$, then $A(\tau')=A(\tau)^\sg$ and so these flip 
subalgebras are isomorphic. Hence, we just need to deal with possible flips up 
to conjugation in $\Aut\Gm$.  

The extras $x=a+b$, where $ab\neq 0$, are not idempotents, and in particular, 
they are not axes. This is why we discard them here. However, it is not 
impossible that the extras may lead to interesting idempotents/axes in a 
different, more intricate way. Hence we record the following.

\begin{question}
Is there a way to produce interesting idempotents/axes from the extras? What 
kind of fusion law do such new axes satisfy?
\end{question}

\section{Symmetric group}

\label{symmetric group}

Until now, there have been only a handful of known Monster type axial 
algebras, all of type $(\qu,\thi)$, that were outside of the smaller class of 
algebras of Jordan type. The results of Section \ref{flip subalgebras} allow 
us to populate the class of algebras of Monster type with many new interesting 
axial algebras and to show that this class is truly much larger than the class 
of algebras of Jordan type. 

There are many families of $3$-transposition groups and for each there are 
many conjugacy classes of flips $\tau$. In this section we study the simplest 
family of $3$-transposition groups, the symmetric groups. For them, we 
identify the flip subalgebra $A(\tau)$ for each flip $\tau$. Hence, in this 
section, $G=S_n$, the symmetric group on $n$ symbols, and $D=(1,2)^G$ is the 
conjugacy class of transpositions ($2$-cycles). We note that $D$, the point 
set of $\Gm$ and the basis of $M=M_\eta(\Gm)$, has cardinality ${n\choose 
2}=\frac{n(n-1)}{2}$. We also note that, in this case, $\Aut\Gm$ coincides with 
$G=S_n$. The classes of involutions $\tau$ in $S_n$ are identified by a single 
parameter, the number $k$, $1\leq k\leq\frac{n}{2}$, of $2$-cycles in the 
decomposition of $\tau$ as a product of independent cycles. 

\subsection{Infinite series}

We first deal with a special situation, where $n=2k$. Without loss of 
generality we may assume that 
$$\tau=(1,2)(3,4)\ldots(2k-1,2k).$$

As above, let $H=\la\tau\ra$ be the cyclic subgroup generated by $\tau$. 
Associated with $\tau$ and $H$ is a partition $\cP$ of $\{1,2,\ldots,n\}$ into 
$k$ parts of size two, each part $P_i=\{2i-1,2i\}$ being an orbit of $H$. 
Each $s\in\{1,2,\ldots,n\}$ lies in some part $P_i$. Let $\bar s$ denote the 
other number in this $P_i$. For example, $1\in P_1=\{1,2\}$ and so $\bar 1=2$.
Note that $\tau$ sends every $s$ to $\bar s$.

We need to classify singles, doubles and extras, and we do it in terms of the 
above partition $\cP$.

\begin{proposition} \label{Qk}
When $n=2k$, the fixed subalgebra $M_H$ has dimension $k^2$. 
Among the orbit vectors, there are (a) $k$ singles $(s,\bar s)$; (b) 
$k^2-k$ doubles $(s,t)+(\bar s,\bar t)$, with $s$ and $t$ contained in 
different parts of $\cP$; and (c) no extras. 
\end{proposition}

\begin{proof}
Let $O$ be the orbit $a^H$, where $a=(s,t)\in D$. If $s$ and $t$ are in 
the same part of the partition $\cP$ then $t=\bar s$ and so 
$(s,t)^\tau=(t,s)=(s,t)$. Hence in this case we obtain $k$ orbits of 
length $1$. These are the $k$ singles appearing in (a) above.

Let us now consider the complementary case, where $s$ and $t$ belong to 
different parts of $\cP$. In this case, $(s,t)^\tau=(\bar s,\bar t)$ 
contains no $s$ and no $t$. So $O=\{(s,t),(\bar s,\bar t)\}$ is of 
length two and of orthogonal type, since the two transpositions in $O$ 
are independent and hence commute, which corresponds to product zero in 
$M$. Thus, all such orbit vectors are doubles. Let us count them. Suppose 
that $s\in P_i$ and $t\in P_j$. This pair of parts gives us two 
doubles: $(s,t)+(\bar s,\bar t)$ and $(s,\bar t)+(\bar s,t)$. There are 
${k\choose 2}=\frac{k(k-1)}{2}$ pairs $\{i,j\}$, each leading to two 
doubles. Hence the total number of doubles is $2\frac{k(k-1)}{2}=k^2-k$, 
as claimed in (b).

We have already covered all possibilities for $O=(s,t)^H$, hence there 
is no room for extras. 

For the claim on the dimension of $M_H$, recall from Proposition 
\ref{basis} that the orbit vectors form a basis of $M_H$. Adding $k$ to 
$k^2-k$ gives the dimension equal to $k^2$. 
\end{proof}

\begin{corollary}
When $n=2k$, the flip algebra $A(\tau)$ coincides with the 
fixed subalgebra $M_H$. In particular, its dimension is $k^2$ and it has 
a basis consisting of singles and doubles. 
\end{corollary}

\begin{proof}
Indeed, as there are no extras, $A(\tau)$ contains the whole basis of 
$M_H$ and thus coincides with it.
\end{proof}

We call this algebra $A(\tau)$ the $k^2$-algebra and denote it 
$Q_k(\eta)$. The $4$-dimensional algebra $Q_2(\eta)$, we found in Section 
\ref{2-generated} for the diagram $A_3$, belongs in this series. (The 
algebra $Q_1(\eta)$ is $1$-dimensional generated by the single axis 
$(1,2)$.)

As we have here a whole infinite series of algebras, one for each value 
of $k$, we cannot show the multiplication table as we did earlier. 
Instead, we will now describe the $2$-generated subalgebras $\dla 
x,y\dra$ involved in $A(\tau)$. Clearly, all of them must be on the list 
we obtained in Section \ref{2-generated}. There are several cases for 
the pair $Y=\{x,y\}$. We do not need to provide a formal proof in each 
case, because all one needs to do is to identify the diagram arising on 
$Z=\supp(Y)$.
\begin{itemize}
\item If both $x$ and $y$ are singles then $xy=0$ and $\dla x,y\dra\cong 
2B$. 
\item Suppose that $x=(s,\bar s)$ is a single, with $\{s,\bar s\}=P_m$, 
and $y$ is a double corresponding to a pair $\{i,j\}$.
\begin{itemize}
\item If $m\not\in\{i,j\}$ then the diagram on $Z$ has no edges (diagram 
$A_1$) and so $\dla x,y\dra\cong 2B$.
\item If $m\in\{i,j\}$ then the diagram on $Z$ is $A_3$, and so $\dla 
x,y\dra$ is the $2^2$-algebra $Q_2(\eta)$ arising for that diagram. 
\end{itemize}
\item Finally suppose that both $x$ and $y$ are doubles, corresponding 
to pairs $\{i,j\}$ and $\{i',j'\}$ respectively. 
\begin{itemize}
\item If $\{i,j\}\cap\{i',j'\}=\emptyset$ then $xy=0$ and again $\dla 
x,y\dra\cong 2B$. (This is diagram $B_1$.)
\item If $i=i'$, but $j\neq j'$, then the diagram on $Z$ is $B_4$ and 
hence $\dla x,y\dra\cong 3C(\eta)$.
\item Lastly, if $i=i'$ and $j=j'$ then $x=(s,t)+(\bar s,\bar t)$ and 
$y=(s,\bar t)+(\bar s,t)$ for some $s\in P_i$ and $t\in P_j$. Here the 
diagram is $B_6$ and, according to the discussion of this case in 
Section \ref{2-generated}, we need to determine the order $p$ of 
$(s,t)^{(\bar s,t)}(\bar s,\bar t)^{(s,\bar t)}$. Since 
$(s,t)^{(\bar s,t)}=(s,\bar s)$ and $(\bar s,\bar t)^{(s,\bar t)}=
(\bar s,s)=(s,\bar s)$, the product is equal to the identity, that is, $p=1$. 
Therefore, in this last case, $\dla x,y\dra\cong 3C(2\eta)$.
\end{itemize}
\end{itemize}
Now the formula for the product $xy$, when it is non-zero, can be looked 
up in the multiplication table of the corresponding algebra $\dla 
x,y\dra$.  

\subsection{General case}

In this subsection we deal with the general case: arbitrary $n$ and $k$. 
Let $m:=n-2k$.  

\begin{theorem} \label{general}
For arbitrary $n$ and $k$, the flip algebra $A(\tau)$ is a direct sum of 
the $k^2$-algebra $Q_k(\eta)$ and the Matsuo algebra $M_\eta(\Gm')$ 
where $\Gm'$ is the Fischer space of $G'\cong S_m$ acting on the 
$m$-element set $\{2k+1,\ldots,n\}$.
\end{theorem}

\begin{proof}
We split the set $\{1,2,\ldots,n\}$ as a disjoint union $R\cup S$, where 
$R=\{1,2,\ldots,2k\}$ and $S=\{2k+1,\ldots,n\}$. Clearly, $H=\la\tau\ra$ 
acts on $R$ with $r$ orbits of length $2$ and it acts on $S$ trivially. 
Consider an arbitrary single axis $a=(s,t)\in D$. If $s\in R$ and $t\in 
S$ then $a^\tau=(\bar s,t)$ (we adapt the bar notation from the 
preceding subsection). Clearly, all such orbits $a^H=\{a,a^\tau\}$ are 
non-orthogonal and so the corresponding orbit vectors are extras. Hence 
each axis of $A(\tau)$, whether a single or a double, comes fully from 
$R$ or fully from $S$. The former kind clearly span the $k^2$-algebra 
$Q_k(\eta)$ and the latter kind are all singles, because $H$ acts 
trivially on $S$, and they span the Matsuo algebra $M_\eta(\Gm')$ of the 
symmetric group of $S$. 

Finally, since $R$ and $S$ are disjoint, any transposition involved in 
an axis $x$ of the first kind is independent from any axis $y$ of the 
second kind. Therefore, $xy=0$, proving that $A(\tau)$ is a direct sum, 
as claimed in the theorem.
\end{proof}

\subsection{Simplicity}

Here we tackle the following question: for which values of $\eta$ is the 
$k^2$-algebra $A=Q_k(\eta)$ not simple and what is then the dimension of the 
ideal? 

Recall from Subsection \ref{ideals} that ideals in an axial algebra are 
classified into two types: the ones containing a generating axis and the ones 
contained in the radical. For the first kind, we need to know the projection 
graph $\Dl$ (see Definition \ref{projection graph}). The vertices of $\Dl$ 
are the generating axes and the edges, in the presence of a Frobenius form, 
are pairs of generating axes, $\{a,b\}$, $a\neq b$, with $(a,b)\neq 0$. So let 
us describe the Gram matrix of the Frobenius form on the set of axes.

Let us give all axes standard names in terms of the partition $\cP$. First of 
all, we have the single $a_i=(2i-1,2i)$ corresponding to the part $P_i$, 
$i=1,\ldots,k$. Also, for $1\leq i<j\leq k$, the two doubles corresponding to 
$\{i,j\}$ are $b_{i,j}:=(2i-1,2j-1)+(2i,2j)$ and 
$c_{i,j}:=(2i-1,2j)+(2i,2j-1)$. Now, let us compute the entries of the Gram 
matrix. First of all, $(a_i,a_i)=1$ and $(a_i,a_j)=0$ for all $i\neq j$. 
Next, $(b_{i,j},b_{i,j})=2=(c_{i,j},c_{i,j})$ and $(b_{i,j},c_{i,j})=
4\frac{\eta}{2}=2\eta$. Taking a single $x=a_i$ and a double $y=b_{s,t}$ or 
$c_{s,t}$, we obtain that $(x,y)=0$ if $i\not\in\{s,t\}$ and $(x,y)=
2\frac{\eta}{2}=\eta$ if $i\in\{s,t\}$. Finally, taking two doubles 
$x=b_{i,j}$ or $c_{i,j}$ and, similarly, $y=b_{s,t}$ or $c_{s,t}$, we 
calculate that $(x,y)=0$ if $\{i,j\}$ and $\{s,t\}$ are disjoint and 
$(x,y)=2\frac{\eta}{2}=\eta$ if $\{i,j\}$ and $\{s,t\}$ meet in one element.  

\begin{proposition} \label{connected}
The projection graph $\Dl$ of $A=Q_k(\eta)$ is connected. In particular, 
$A$ has no proper ideals containing any of the generating axes.
\end{proposition}

\begin{proof}
From the above description of the Gram matrix, it is clear that $a_1$ is 
adjacent to all $b_{1,j}$ and $c_{1,j}$, $2\leq j\leq k$, and these are 
adjacent to all the remaining vertices of $\Dl$. According to the discussion 
around Theorem \ref{unoriented}, connectedness of $\Dl$ means that $A$ has 
no proper ideals containing axes. 
\end{proof}

Therefore, every proper ideal of $A$ is contained in the radical $R(A)$. 
Furthermore, by Theorem \ref{perp}, since none of the axes is singular, $R(A)$ 
coincides with $A^\perp$, the radical of the Frobenius form. Therefore, $A$ is 
not simple exactly when the determinant of the Gram matrix is zero. 

We next compute this determinant. Our plan is to find a basis with respect to 
which the Gram matrix splits into several blocks. Let us define 
$d_{i,j}:=\frac{b_{i,j}-c_{i,j}}{2}$ and $e_{i,j}:=b_{i,j}+c_{i,j}$. Then the 
transition matrix from the standard basis to the basis, consisting of all $a_i$ 
and all $d_{i,j}$ and $e_{i,j}$, is block-diagonal with ${k\choose 2}$ blocks 
of size $2\times 2$ and determinant one. Hence the new Gram matrix has the same
determinant as the original one. 

From the values of the Frobenius form listed above, we note that 
$(b_{i,j},x)=(c_{i,j},x)$ for all axes $x$ excluding $b_{i,j}$ and $c_{i,j}$. 
This means that each $d_{i,j}$ is orthogonal to all vectors in the new basis 
apart from itself and, possibly, $e_{i,j}$. We compute that 
$(d_{i,j},d_{i,j})=\frac{1}{4}(b_{i,j}-c_{i,j},b_{i,j}-c_{i,j})=\frac{1}{4}(2-2
(2\eta)+2)=1-\eta$, non-zero since $\eta\neq 1$. However, 
$(d_{i,j},e_{i,j})=\frac{1}{2}(b_{i,j}-c_{i,j},b_{i,j}+c_{i,j})= 
\frac{1}{2}(2+2\eta-2\eta-2)=0$. Therefore, the subspace $D$ spanned by all 
vectors $d_{i,j}$ is non-degenerate (the Gram matrix on $D$ is $(1-\eta)I$, 
where $I$ is the identity matrix) and it is orthogonal to the subspace $E$ 
spanned by all $a_i$ and all $e_{i,j}$. 

Let us now focus on the subspace $E$. Each $e_{i,j}$ is already orthogonal to 
all $a_s$ with $s\not\in\{i,j\}$. Let us amend $e_{i,j}$ so that the corrected 
vectors are also orthogonal to $a_i$ and $a_j$. Consider 
$f_{i,j}:=e_{i,j}-\al(a_i+a_j)$. Then $(f_{i,j},a_i)=(b_{i,j}+c_{i,j}-\al 
a_i-\al a_j,a_i)=\eta+\eta-\al-0$. Let us select $\al=2\eta$. Then $f_{i,j}$ 
is orthogonal to $a_i$ and, symmetrically, to $a_j$. The transition matrix is 
again of determinant one, and the subspace $W$ spanned by all $a_i$ is 
orthogonal to the subspace $F$ spanned by all $f_{i,j}$. The Gram matrix on 
$W$ is the identity matrix, so $W$ is non-degenerate. Hence the radical of $A$ 
coincides with the radical of the subspace $F$ and so we can focus on the 
latter.

Let us compute the Gram matrix on $F$ with respect to the basis consisting of 
all $f_{i,j}$. The basis vectors are in a bijection with subsets $\{i,j\}$ and 
so we have three cases to consider. First of all, if $\{i,j\}$ and $\{s,t\}$ 
are disjoint then $(f_{i,j},f_{s,t})=0$. If $\{i,j\}$ and $\{s,t\}$ meet in 
one element, say, $i=s$ then $(f_{i,j},f_{i,t})=(b_{i,j}+c_{i,j}-2\eta a_i- 
2\eta a_j,b_{i,t}+c_{i,t}-2\eta a_i-2\eta a_t)=\eta+\eta-2\eta\eta+0+\eta+
\eta-2\eta\eta+0-2\eta\eta-2\eta\eta+4\eta^2+0+0+0+0+0=
4\eta-4\eta^2=4(-\eta^2+\eta)=4\eta(1-\eta)$. Finally, $(f_{i,j},f_{i,j})= 
(b_{i,j}+c_{i,j}-2\eta a_i-2\eta a_j,b_{i,j}+c_{i,j}-2\eta a_i-2\eta 
a_j)=2+2+4\eta^2+4\eta^2+2(2\eta-2\eta\frac{\eta}{2}-2\eta\frac{\eta}{2}- 
2\eta\frac{\eta}{2}-2\eta\frac{\eta}{2}+0)=4+4\eta-8\eta^2=4(-2\eta^2+\eta+1)=
4(1-\eta)(2\eta+1)$. 

Let $m:={k\choose 2}=\frac{k(k-1)}{2}$ be the dimension of $F$. Let $I$ be 
the identity matrix of this size and $J$ be the adjacency matrix of the 
Johnson graph $J(k,2)$, the graph on $2$-element subsets of 
$\{1,2,\ldots,k\}$, where two $2$-subsets are adjacent when they meet in one 
element. Then the Gram matrix of $F$ with respect to the basis formed by the 
vectors $f_{i,j}$ is equal to $4((-2\eta^2+\eta+1)I+(-\eta^2+\eta)J)$. Over a 
field of characteristic zero, the eigenvalues of $J$ and their multiplicities 
are well-known, since $J(k,2)$ is strongly regular. They are the degree, 
$\theta_0=2(k-2)$, with multiplicity $1$, $\theta_1=k-4$ with multiplicity 
$k-1$, and $\theta_2=-2$ with multiplicity $\frac{k(k-3)}{2}$. 

\begin{proposition}
The determinant of the Gram matrix of the Frobenius form on $A=Q_k(\eta)$ is 
equal to $(2(1-\eta))^{k^2-k}(2(k-1)\eta+1)((k-2)\eta+1)^{k-1}$.
\end{proposition}

\begin{proof}
We first assume that $\F$ is of characteristic zero. Since the Gram matrix on 
$F$ coincides with $4((-2\eta^2+\eta+1)I+(-\eta^2+\eta)J)$, it has eigenvalues 
$\kappa_h=4\eta(1-\eta)\theta_h+4(1-\eta)(2\eta+1)=4(1-\eta) 
(\eta\theta_h+2\eta+1)$, $h=0,1,2$, with the corresponding multiplicities. 
This gives $\kappa_0=4(1-\eta)(2(k-1)\eta+1)$, $\kappa_1=
4(1-\eta)((k-2)\eta+1)$, and $\kappa_2=4(1-\eta)$. Since the determinant of 
a square matrix is the product of its eigenvalues taken with their 
multiplicities, we conclude that the determinant of the Gram matrix on $F$ is 
$(4(1-\eta))^m(2(k-1)\eta+1)((k-2)\eta+1)^{k-1}$. Combining this with the 
determinant of the Gram matrix on $D$, equal to $(1-\eta)^m$, and with the 
determinant of the Gram matrix on $W$, equal to $1$, we obtain the polynomial 
claimed in the proposition.

While our computation was over a field of characteristic zero, the result is a 
polynomial with integral coefficients. Therefore, it readily transfers into 
any odd characteristic and so our claim holds for any field $\F$.
\end{proof}

In this proof and in the discussion before the proposition, we glazed over the 
small cases $k=1$ and $2$, where some of the eigenvalues disappear and the 
multiplicity formula does not apply. However, by inspection, the formula for 
the determinant remains true also in these small cases.

As a consequence, we have the following statement.

\begin{theorem} \label{critical values}
The $k^2$-algebra $Q_k(\eta)$ is simple unless $2(k-1)\eta+1=0$ (equivalently, 
$\eta=-\frac{1}{2(k-1)}$) or $(k-2)\eta+1=0$ (equivalently, 
$\eta=-\frac{1}{k-2}$).
\end{theorem}

We now turn to the dimension of the radical for the special values of $\eta$. 
In zero characteristic, the multiplicity of $\kappa_0$, equal to $1$, 
(respectively, of $\kappa_1$, equal to $k-1$) gives the dimension of the 
radical when $\eta=-\frac{1}{2(k-1)}$ (respectively, $\eta=-\frac{1}{k-2}$). 
In positive characteristic, the same applies as long as 
$\kappa_0\neq\kappa_1$. This can only happen when $k=0$ (i.e., the 
characteristic of $\F$ divides $k$) and then the special value of $\eta$ is 
$\eta=\frac{1}{2}$.

We should also comment on the small cases. The above analysis fully applies 
when $k\geq 3$. When $k=1$, $A=Q_1(\eta)\cong\F$ and so is simple. When $k=2$, 
only the eigenvalue $\kappa_0$ can be zero, and so $A=Q_2(\eta)$ has a 
non-zero radical only for $\eta=-\frac{1}{2(2-1)}=-\frac{1}{2}$. Then the 
radical is $1$-dimensional (cf. Section \ref{2-generated}).

\section{Two further series}

\label{two further series}

In this final section of the paper we construct two further infinite 
series of algebras of Monster type generalizing two examples we found in 
Section \ref{3-generated}.

The $3$-transposition groups that we consider here are generalizations 
of the groups $\hat G(p)$ that we first encountered in the case of 
diagram $B_6$ in Section \ref{2-generated}. Both of these groups can be 
obtained via a more general construction due to Zara \cite{z,z1} in 
wreath products $K\wr S_n$, where $K$ has all elements of order at 
most $3$. Our two groups arise for $K=C_2$ or $C_3$.

\subsection{Series $2Q_k(\eta)$}

First, consider the semi-direct product $\hat G_n=2^n:S_n$ of the 
symmetric group $S=S_n$ and its permutational module $E=2^n$ over 
$\F_2$. The class $C$ of transpositions from $S$ extends to a larger 
conjugacy class $D$ in $\hat G_n$. We claim that $D$ is again a class of 
$3$-transpositions. Let $e_1,\ldots,e_n$ be the basis of $E$ permuted by 
$S$. Note that we retain the additive notation for the operation in $E$. 
The transposition $\sg:=(i,j)\in C$ acts as a transvection on $E$ with 
$[E,\sg]=\la e_i+e_j\ra$ and $C_E(\sg)=\la e_i+e_j\ra+\la e_s\mid s\neq 
i,j\ra$. Therefore, $|D|=2|C|=n(n-1)$; namely, the coset $E\sg$, in 
addition to $\sg$, contains the second element from $D$, 
$\tilde\sg=(e_i+e_j)\sg$. Note that this description of $D$ implies that 
$G_n=\la D\ra$ is of index $2$ in $\hat G_n$, namely, $G_n\cap E\cong 
2^{n-1}$ is the ``sum-zero'' submodule of $E$.

Manifestly, all elements of $D$ are involutions. Any two distinct 
elements from $D$ involve no more than four indices from 
$\{1,2,\ldots,n\}$. Therefore, they are contained in a subgroup 
isomorphic to the group $\hat G(2)\cong 2^3:S_4$ from Section \ref{2-generated}, 
where they are contained in the class of $3$-transpositions. Hence $D$ 
is also a class of $3$-transpositions and $(G_n,D)$ is a 
$3$-transposition group. Let $\Gm$ be the Fischer space of $(G_n,D)$ and 
$M=M_\eta(\Gm)$.

Let $n\geq 2k$ and let $\tau=(1,2)(3,4)\cdots(2k-1,2k)$ be the flip, as 
in Section \ref{symmetric group}. Let $H=\la\tau\ra$.

We start with a special case.

\begin{proposition}
If $n=2k$ then the fixed subalgebra $M_H$ is of dimension $2k^2$. Among 
the orbit vectors, there are $2k$ singles, $2(k^2-k)$ doubles, and no 
extras.
\end{proposition}

\begin{proof}
The singles are the points $(2i-1,2i)$ and $(e_{2i-1}+e_{2i})(2i-1,2i)$. 
All other orbit sums are of the form $\{(s,t)+(\bar s,\bar t)\}$ or 
$\{(e_s+e_t)(s,t)+(e_{\bar s}+e_{\bar t})(\bar s,\bar t)\}$, and they 
are all doubles. Here $s$ and $t$ are in different parts of the 
partition $P=\{\{1,2\},\{3,4\},\ldots,\{2k-1,2k\}\}$. Note that we adopt 
the bar notation from Section \ref{symmetric group} for the 
complementary elements in the parts.

Comparing with Proposition \ref{Qk}, we have twice as many singles here 
and twice as many doubles, so the claim follows.
\end{proof}

Since $M_H$ is spanned by singles and doubles, it coincides with the 
flip algebra $A(\tau)$. We use $2Q_k(\eta)$ to denote this algebra $M$. The 
$8$-dimensional example from the case $C_5$ from Section 
\ref{3-generated} is isomorphic to $2Q_2(\eta)$. 

Let us also include a statement analogous to Theorem \ref{general}. 
We note, however, that, unlike there, we cannot claim now that this 
result covers all possible flips.

\begin{proposition}
If $n=2k+m$ then $A(\tau)$ is the direct sum of $2Q_k(\eta)$ arising on the 
subset $R:=\{1,2,\ldots, 2k\}$ and the Matsuo algebra $M_\eta(\Gm')$, where 
$\Gm'$ is the Fischer space of the group $G_m\cong 2^{m-1}:S_m$ arising 
on the subset $S:=\{2k+1,\ldots,n\}$.
\end{proposition}

We skip the proof as it is quite analogous to that of Theorem 
\ref{general}.

Concerning the question of simplicity of $2Q_k(\eta)$, it is easy to see that 
the projection graph on the set of single and double axes in $2Q_k(\eta)$, 
$k\geq 2$, is connected. Hence $2Q_k(\eta)$ contains no proper ideals 
containing axes and so it is simple for all but a finite number of 
values of $\eta$. We leave determination of these special values of 
$\eta$ until another time. 

\subsection{Series $3Q_k(\eta)$}

This is similar to the first case, except now we take $\hat G_n$ to be 
the semi-direct product of $S=S_n$ and the permutational module $E=3^n$ 
over $\F_3$. Again, let $e_1,\ldots,e_n$ be the basis permuted by $S$. 
Let $D$ be the conjugacy class of $\hat G_n$ containing the class $C$ of 
transpositions from $S$. 

Taking $\sg=(i,j)$, it is easy to see that $[E,\sg]=\la e_i-e_j\ra$ and 
$C_E(\sg)=\la e_i+e_j\ra+\la e_s\mid s\neq i,j\ra$, that is, $\sg$ acts 
on $E$ as a reflection. From here we deduce that the coset $E\sg$ 
contains three elements from $D$, namely, $\sg$ and two further 
elements, $(e_i-e_j)\sg$ and $(e_j-e_i)\sg$. Thus, 
$|D|=3|C|=\frac{3}{2}n(n-1)$. Also, it follows that any two elements 
from $D$ involve no more than four indices from $\{1,2,\ldots,n\}$ and 
so they are contained in a subgroup $\hat G(3)=3^3:S_4$ as in Section 
\ref{2-generated}. Since $D\cap\hat G(3)$ is a class of 
$3$-transpositions in $\hat G(3)$, we deduce that $D$ is a class of 
$3$-transpositions. Also similarly to the first case, $G_n=\la D\ra$ is 
a proper, index $3$ subgroup of $\hat G_n$, with $G_n\cap E\cong 3^{n-1}$ 
being the ``sum-zero'' submodule in $E$.

Let $\Gm$ be the Fischer space of the $3$-transposition group $(G_n,D)$ 
and let $M=M_\eta(\Gm)$. Take $k\leq\frac{n}{2}$. Our flip $\tau$ is again 
based on $\tau_0=(1,2)(3,4)\cdots(2k-1,2k)$. However, we need a correcting 
factor to achieve the subalgebra we want. Let $z$ be the automorphism of 
$\hat G_n$, centralizing the complement $S_n$ and inverting every element 
of $E$. We set $\tau=z\tau_0$. Since $z$ and $\tau_0$ commute, $\tau$ is 
an involution.

Let $H=\la\tau\ra$.

\begin{proposition}
Suppose that $n=2k$. Then the fixed subalgebra $M_H$ is of dimension 
$3k^2$ and spanned by $3k$ singles, $3(k^2-k)$ doubles, and no 
extras.
\end{proposition}

\begin{proof}
We claim that the points $(2i-1,2i)$, $(e_{2i-1}-e_{2i})(2i-1,2i)$, and 
$(e_{2i}-e_{2i-1})(2i-1,2i)$, where $i=1,\ldots,k$, are singles. Indeed, 
$(2i-1,2i)^{\tau}=((2i-1,2i)^z)^{\tau_0}=(2i-1,2i)^{\tau_0}=(2i-1,2i)$, 
so $(2i-1,2i)$ is a single. Also, $((e_{2i-1}-e_{2i})(2i-1,2i))^{\tau}=
(((e_{2i-1}-e_{2i})(2i-1,2i))^z)^{\tau_0}=
((-e_{2i-1}+e_{2i})(2i-1,2i))^{\tau_0}=(-e_{2i}+e_{2i-1})(2i-1,2i)=
(e_{2i-1}-e_{2i})(2i-1,2i)$, so this is also a single. The calculation 
for $(e_{2i}-e_{2i-1})(2i-1,2i)$ is quite similar. 

All other orbit vectors are of the form $(s,t)+(\bar 
s,\bar t)$ or $(e_s-e_t)(s,t)+(e_{\bar t}-e_{\bar s})(\bar s,\bar t)$, 
where $s$ and $t$ are in different parts of the partition 
$P=\{\{1,2\},\{3,4\},\ldots,\{2k-1,2k\}\}$ associated with $\tau$. Note 
that the two summands in these orbit vectors have disjoint support, so 
they are all doubles, leaving no room for extras.

Note also that $(e_t-e_s)(s,t)=(e_t-e_s)(t,s)$ and so our discussion 
above indeed covers all orbit vectors.
\end{proof}

From this we see that the flip algebra $A(\tau)$ coincides with $M_H$. 
Our notation for this flip algebra is $3Q_k(\eta)$. The $12$-dimensional 
examples from the case $C_5$ from Section \ref{3-generated} is 
$3Q_2(\eta)$.

For a general $n=2k+m$ and $\tau=z\tau_0$, we again get a direct sum 
decomposition.

\begin{proposition}
Suppose that $n=2k+m$. Then $A(\tau)$ is the direct sum of $3Q_k(\eta)$ 
arising on the subset $\{1,2,\ldots, 2k\}$ and the Matsuo algebra 
$M_\eta(\Gm')$, where $\Gm'$ is the Fischer space of the group $S_m$ 
acting on the subset $\{2k+1,\ldots,n\}$.
\end{proposition}

We omit the proof. Note, however, that the Matsuo summand $M_\eta(\Gm')$ 
here is not for the group $G_m$, as in the first case, but rather for 
the group $S_m$, as in Theorem \ref{general}. This is because of the 
correcting factor $z$.

The projection graph of $3Q_k(\eta)$ is again easily seen to be 
connected, and so $3Q_k(\eta)$ is simple for all but a finite number of 
special values of $\eta$.

\Addresses
\end{document}